%% file: main_arxiv.tex
\documentclass[a4paper, twoside]{scrartcl}

\usepackage{amsfonts}
\usepackage{amssymb}
\usepackage{amsmath}
\usepackage{amsthm}
\usepackage{amstext}
\usepackage{amssymb}
\usepackage{mathtools}
\usepackage{caption}
\usepackage{subcaption}
\usepackage{float}
\usepackage{color}
\usepackage{graphicx}
\usepackage{enumerate}
\usepackage{dsfont}  
\usepackage[normalem]{ulem} 
\usepackage[authoryear,sort,round]{natbib}  
\usepackage{url}
\usepackage{hyperref}
\usepackage[noabbrev,capitalise,nosort,nameinlink]{cleveref}
\crefname{assumption}{Assumption}{Assumptions}

\usepackage{etex}
\usepackage[a4paper, top=88pt, bottom=88pt, left=72pt, right=72pt, headsep=16pt, footskip=28pt]{geometry}
\usepackage{hyperref}
\hypersetup{
	colorlinks,
	linkcolor=blue,
	citecolor=blue,
	urlcolor=blue,
}
\newcommand*{\email}[1]{\bgroup\color{blue}\href{mailto:#1}{#1}\egroup}
\usepackage{enumitem}
\setlist[enumerate]{nosep}
\setlist[itemize]{nosep}
\usepackage{mleftright} \mleftright
\renewcommand{\qedsymbol}{$\blacksquare$}
\renewenvironment{proof}[1][\proofname]{\noindent{\sffamily\bfseries #1.} }{\hfill \qedsymbol \medskip}
\usepackage{scrlayer-scrpage, xhfill}
\usepackage[labelfont={sf,bf}]{caption}
\automark[section]{section}
\setkomafont{pageheadfoot}{\normalcolor\sffamily}
\setkomafont{pagenumber}{\normalfont\normalsize\sffamily}
\clearpairofpagestyles
\let\oldtitle\title
\renewcommand{\title}[1]{\oldtitle{#1}\newcommand{\theshorttitle}{#1}}

\let\oldauthor\author
\renewcommand{\author}[1]{\oldauthor{#1}\newcommand{\theshortauthor}{#1}}
\newcommand{\shortauthor}[1]{\renewcommand{\theshortauthor}{#1}}
\cohead{\xrfill[0.525ex]{0.6pt}~\theshorttitle~\xrfill[0.525ex]{0.6pt}}
\cehead{\xrfill[0.525ex]{0.6pt}~\theshortauthor~\xrfill[0.525ex]{0.6pt}}
\newcommand{\theabstract}[1]{\par\bgroup\noindent\textbf{\textsf{Abstract.}} #1\egroup}
\newcommand{\thekeywords}[1]{\par\smallskip\bgroup\noindent\textbf{\textsf{Keywords.}}\newcommand{\and}{ $\bullet$ } #1\egroup}
\newcommand{\themsc}[1]{\par\smallskip\bgroup\noindent\textbf{\textsf{2010 Mathematics Subject Classification.}}\newcommand{\and}{ $\bullet$ } #1\egroup}
\newcommand*{\affilref}[1]{\ref{affiliation#1}}
\newcommand*{\affiliation}[3]{
	\footnotetext[#1]{\label{affiliation#2} #3}
}

\usepackage{bm}

\newtheorem{theorem}{\sffamily Theorem}
\newtheorem{lemma}{\sffamily Lemma}
\newtheorem{proposition}{\sffamily Proposition}
\newtheorem{corollary}{\sffamily Corollary}
\theoremstyle{definition}
\newtheorem{definition}{\sffamily Definition}

\newtheorem{remark}{\sffamily Remark}

\newcommand{\todo}[1]{\bgroup\color{red}#1\egroup}

\title{Bayesian Numerical Methods for Nonlinear Partial Differential Equations}

\author{%
	Junyang Wang\textsuperscript{\affilref{Newcastle}}
	\and
	Jon Cockayne\textsuperscript{\affilref{Turing}}
	\and
	Oksana Chkrebtii\textsuperscript{\affilref{OSU}}
	\and
	T.~J.\ Sullivan\textsuperscript{\affilref{Warwick}}
	\and
	Chris.\ J.\ Oates\textsuperscript{\affilref{Newcastle},\affilref{Turing}}
}
\shortauthor{J.~Wang, J.~Cockayne, O.~Chkrebtii, T.~J.~Sullivan, and C.~J.~Oates}

\begin{document}

\maketitle
\affiliation{1}{Newcastle}{School of Mathematics, Statistics \& Physics, Newcastle University, Newcastle-Upon-Tyne, NE1 7RU, United Kingdom. \email{j.wang68@ncl.ac.uk}, \email{chris.oates@ncl.ac.uk}}
\affiliation{2}{Turing}{Alan Turing Institute, British Library, 96 Euston Road, London, NW1 2DB, United Kingdom. \email{jcockayne@turing.ac.uk}}
\affiliation{3}{OSU}{Department of Statistics, Ohio State University, 1958 Neil Avenue, 429 Cockins Hall, Columbus, OH 43210-1247, United States. \email{oksana@stat.osu.edu}}
\affiliation{4}{Warwick}{Mathematics Institute and School of Engineering, The University of Warwick, Coventry, CV4 7AL, United Kingdom. \email{t.j.sullivan@warwick.ac.uk}}

\begin{abstract}
	\theabstract{\input{abstract}}
	
	\thekeywords{{approximate likelihood}\and{inverse problem}\and{Mat\'ern covariance}\and{probabilistic numerics}\and{uncertainty quantification}}



\end{abstract}

\input{introduction}

\input{methods}

\input{theory}

\input{results}

\input{conclusion}

\appendix
\section{Appendix}

\input{appendix}

\paragraph*{Acknowledgements}
\input{acknowledgements}

\input{main_arxiv.bbl}
\end{document}

%% file: abstract.tex
	The numerical solution of differential equations can be formulated as an inference problem to which formal statistical approaches can be applied.
	However, nonlinear partial differential equations (PDEs) pose substantial challenges from an inferential perspective, most notably the absence of explicit conditioning formula.
	This paper extends earlier work on linear PDEs to a general class of initial value problems specified by nonlinear PDEs, motivated by problems for which evaluations of the right-hand-side, initial conditions, or boundary conditions of the PDE have a high computational cost.
	The proposed method can be viewed as exact Bayesian inference under an approximate likelihood, which is based on discretisation of the nonlinear differential operator.
	Proof-of-concept experimental results demonstrate that meaningful probabilistic uncertainty quantification for the unknown solution of the PDE can be performed, while controlling the number of times the right-hand-side, initial and boundary conditions are evaluated.
	A suitable prior model for the solution of the PDE is identified using novel theoretical analysis of the sample path properties of Mat\'{e}rn processes, which may be of independent interest.

%% file: introduction.tex
\section{Introduction}
\label{sec: intro}


Classical numerical methods for differential equations produce an approximation to the solution of the differential equation whose error (called \textit{numerical error}) is uncertain in general.
For non-adaptive numerical methods, such as finite difference methods, the extent of numerical error is often estimated by comparing approximations at different mesh sizes \citep{strikwerda2004finite}, while for adaptive numerical methods, such as certain finite element methods, the global tolerance that was given as an algorithmic input is used as a proxy for the size of the numerical error \citep{gratsch2005posteriori}.
On the other hand, probability theory provides a natural language in which uncertainty can be expressed and, since the solution of a differential equation is unknown, it is interesting to ask whether probability theory can be applied to quantify numerical uncertainty.
This perspective, in which numerical tasks are cast as problems of statistical inference, is pursued in the field of \emph{probabilistic numerics} \citep{Larkin1972,Diaconis1988,OHagan1992,Hennig2015,Oates2019}.
A probabilistic numerical method (PNM) describes uncertainty in the quantity of interest by returning a probability distribution as its output. 
The Bayesian statistical framework is particularly natural in this context \citep{Cockayne2017}, since the output of a Bayesian PNM carries the formal semantics of \textit{posterior uncertainty} and can be unambiguously interpreted.
However, \cite{Wang2019} argued that a rigorous Bayesian treatment of ordinary differential equations (ODEs) may only be possible in a limited context (namely, when the ODE admits a solvable Lie algebra). 
This motivates the development of approximate Bayesian PNMs, which aim to approximate the differential equation in such a way that exact Bayesian inference can be performed, motivated in particular by challenging numerical tasks for which numerical uncertainty cannot be neglected. 

This paper focuses on partial differential equations (PDEs); in particular, we focus on PDEs whose governing equations must be evaluated pointwise at high computational cost.
To date, no (approximate or exact) Bayesian PNM for the numerical solution of nonlinear PDEs has been proposed.
However, the cases of nonlinear ODEs and linear PDEs have each been studied.
In \cite{Chkrebtii2016} the authors constructed an approximate Bayesian PNM for the solution of initial value problems specified by either a nonlinear ODE or a linear PDE.
The approach relied on conjugacy of Gaussian measures under linear operators
to cycle between assimilating and simulating gradient field ``data'' in parallel.
(Recall that the gradient field of an ODE can be a nonlinear function of the state variable, necessitating a linear approximation of the gradient field to obtain a conjugate Gaussian treatment.)
The case of ODEs was discussed by \cite{Cockayne2017}, who proposed a general computational strategy called \textit{numerical disintegration} (similar to a sequential Monte Carlo method for rare event simulation) that can in principle be used to condition on nonlinear functionals of the unknown solution, such as provided by the gradient field, and thereby instantiate an \textit{exact} Bayesian PNM.
However, the authors emphasised that such methods are typically impractical.
\citet{Tronarp2018,Tronarp2020} demonstrated how nonlinear filtering techniques can be used to obtain low-cost approximations to the PNM described by \cite{Cockayne2017} and applied their approach to approximate the solution of ODEs.
Further work on approximate Bayesian PNM for ODEs includes \cite{Skilling1992,Schober2014,Schober2016,Kersting2016,Kersting2018,Teymur2016,Teymur2018,Chkrebtii2019,Kraemer2020,Bosch2020}.
The focus of much of this earlier work was in the design of a PNM whose mean or mode coincides in some sense with a classical numerical method for ODEs (e.g.\ a Runge--Kutta method).
For more detail on existing PNM for ODEs the reader is referred to the recent survey of \citet{Wang2019}.
Linear elliptic PDEs were considered by \citet{sarkka2011linear,owhadi2015bayesian,Cockayne2016} and \citet{Chkrebtii2016}, and in this setting exact Gaussian conditioning can be performed.
\citet{Chkrebtii2016} also presented a nonlinear PDE (Navier--Stokes), but a pseudo-spectral projection in Fourier space was applied at the outset to transform the PDE into a system of first-order ODEs --- an approach that exploited the specific form of that PDE.
\citet{chen2021solving} extended \citet{owhadi2015bayesian} to nonlinear PDEs, focusing on \textit{maximum a posteriori} estimation as opposed to uncertainty quantification for the solution of the PDE.

This paper presents the first (approximate) Bayesian PNM for numerical uncertainty quantification in the setting of nonlinear PDEs. 
Our strategy is based on local linearisation of the nonlinear differential operator, in order to perform conjugate Gaussian updating in an approximate Bayesian framework. 
Broadly speaking, our approach is a natural generalisation of the approach taken by \citet{Chkrebtii2016} for ODEs, but with local linearisation to address the additional challenges posed by nonlinear PDEs.
The aim is to quantify numerical uncertainty with respect to the unknown solution of the PDE.
An important point is that we consider only PDEs for which evaluation of either the right-hand side or the initial or boundary conditions is associated with a high computational cost; we do not aim to numerically solve PDEs for which a standard numerical method can readily be employed to drive numerical error to a negligible level, nor do we aim to compete with standard numerical methods in terms of CPU requirement.
Such problems occur in diverse application areas, such as modelling of ice sheets, carbon and nitrogen cycles \citep{hurrell2013community}, species abundance and ecosystems \citep{fulton2010approaches}, each in response to external forcing from a meteorlogical model, or in solving PDEs that themselves depend on the solution of an auxiliary PDE, which occur both when \textit{operator splitting} methods are used \citep{macnamara2016operator} and when \textit{sensitivity equations}, expressing the rate of change of the solution of a PDE with respect to its parameters, are to be solved \citep{petzold2006sensitivity,cockayne2020probabilistic}.
These applications provide strong motivation for PNM, since typically it will not be possible to obtain an accurate approximation to the solution of the PDE and the rich, probabilistic description of numerical uncertainty provided by a PNM can be directly useful \citep[e.g.][]{oates2019bayesian}.

The remainder of the paper is structured as follows:
In \Cref{sec: methods} the proposed method is presented.
The choice of prior is driven by the mathematical considerations described in \Cref{sec: prior}.
A detailed experimental assessment is performed in \Cref{sec: results}.
Concluding remarks are contained in \Cref{sec: conc}.

%% file: methods.tex
\section{Methods} \label{sec: methods}

In \Cref{subsec: set up} we present the general form of the nonlinear PDE that we aim to solve using PNM.
The use of finite differences for local linearisation is described in \Cref{subsec: FD}.
Then, in \Cref{subsec: approach} we present our proposed approximate Bayesian PNM, discussing how computations are performed and how the associated uncertainty is calibrated.

\subsection{Set-Up and Notation} \label{subsec: set up}

For a set $S \subseteq \mathbb{R}^d$, let $C^0(S)$ denote the vector space of continuous functions $c \colon S \to \mathbb{R}$.
For two \textit{multi-indices} $\alpha, \beta \in \mathbb{N}_0^d$, we write $\alpha \leq \beta$ if $\alpha_i \leq \beta_i$ for each $i = 1,\dots,d$.
For a multi-index $\beta \in \mathbb{N}_0^d$, we let $|\beta| = \beta_1 + \dots + \beta_d$ and let $C^\beta(S) \subseteq C^0(S)$ denote those functions $c$ whose partial derivatives 
\[
\partial^\alpha c \coloneqq \partial_{z_1}^{\alpha_1} \dots \partial_{z_d}^{\alpha_d} c(z) \coloneqq \frac{ \partial^{|\alpha|} c(z) }{ \partial z_1^{\alpha_1} \dots \partial z_d^{\alpha_d} }, \qquad \alpha \leq \beta
\] 
exist and are continuous in $S$.

Let $T \in (0,\infty)$ and let $\Gamma$ be an open and bounded set in $\mathbb{R}^d$, whose boundary is denoted $\partial \Gamma$.
Let $\beta \in \mathbb{N}_0^d$ be a multi-index and consider a differential operator
\begin{align*}
D \colon C^\beta([0,T] \times \Gamma) \to C^0([0,T] \times \Gamma)
\end{align*}
and the associated initial value problem with Dirichlet boundary conditions
\begin{equation} \label{eq: somepde}
\begin{alignedat}{2}
Du(t,x) & =f(t,x), && \qquad t \in [0,T], \; x \in \Gamma  \\
u(0,x) & =g(x), && \qquad x \in \Gamma  \\
u(t,x) &= h(t,x), && \qquad t \in [0,T], \; x \in \partial \Gamma 
\end{alignedat}
\end{equation}
whose unique classical (i.e.\ strong) solution $u \in C^\beta([0,T] \times \Gamma)$ is assumed to exist\footnote{The existence of a strong solution is a nontrivial assumption, since several PDEs admit only a weak solution; see Section 1.3.2 of \cite{evans1998partial} for definitions and background. A well-known class of classical numerical methods that also presuppose the existence of a strong solution are the radial basis function methods \citep{Fornberg2015}. In \Cref{subsec: porous} we consider, empirically, the performance of the method developed in this paper when applied to a PDE for which a strong solution does not exist.}.
The task considered in this paper is to produce a probability distribution on $C^\beta([0,T] \times \Gamma)$ that (approximately) carries the semantics of Bayesian inference for $u$; i.e.\ we seek to develop an (approximate) Bayesian PNM for the numerical solution of \eqref{eq: somepde} \citep{Cockayne2017}.
In particular, we are motivated by the problems described in \Cref{sec: intro}, for which evaluation of $f$, $g$ and $h$ are associated with a high computational cost.
Such problems provide motivation for a careful quantification of uncertainty regarding the unknown solution $u$, since typically it will not be possible to obtain a sufficient number of evaluations of $f$, $g$ and $h$ in order for $u$ to be precisely identified.

\subsubsection{Why Not Emulation?}

Given that the dominant computational cost is assumed to be evaluation of $f$, $g$ and $h$, it is natural to ask whether the uncertainty regarding these functions can be quantified using a probabilistic model, such as an \textit{emulator} \citep{kennedy2001bayesian}.
This would in principle provide a straight-forward Monte Carlo solution to the problem of quantifying uncertainty in the solution $u$ of \eqref{eq: somepde}, where first one simulates an instance of $f$, $g$ and $h$ from the emulator and then applies a classical numerical method to solve \eqref{eq: somepde} to high numerical precision.
The problem with this approach is that construction of a defensible emulator is difficult; the functions $f$, $g$ and $h$ are coupled together by the nonlinear PDE in \eqref{eq: somepde} and, for example, it cannot simultaneously hold that each of $f$, $g$ and $h$ are Gaussian processes.
In fact, the challenge of ensuring that samples of $f$, $g$ and $h$ are consistent with the existence of a solution to \eqref{eq: somepde} poses a challenge that is comparable with solving the PDE itself.
This precludes a straight-forward emulation approach to \eqref{eq: somepde} and motivates our focus on PNM in the remainder, where uncertainty is quantified in the solution space of \eqref{eq: somepde}.

\subsection{Finite Difference Approximation of Differential Operators} \label{subsec: FD}

If $D$ is linear then the differential equation in \eqref{eq: somepde} is said to be \textit{linear} and one or more of the Bayesian PNM of \cite{Chkrebtii2016,Cockayne2016,Chkrebtii2019} may be applied (assuming any associated method-specific requirements are satisfied).
If $D$ is \textit{nonlinear} then at most we can express $D=P+Q$, where $P$ is linear and $Q$ is nonlinear (naturally such representations are non-unique in general).
For example, for Burgers' equation
\begin{equation*}
Du = \frac{\partial u}{\partial t} + u\frac{\partial u}{\partial x}-\varepsilon\frac{\partial^2 u}{\partial x^2}=0 \label{eq: Burger}
\end{equation*}
we have both $P= \partial_t - \varepsilon \partial_x^2$, $Q = u \partial_x$ and also the trivial $Q = D$, $P = 0$.
In this paper we aim, given a decomposition of $D$ in terms of $P$ and $Q$, to adaptively approximate $Q$ by a linear operator, in order that exact Gaussian conditioning formulae can be exploited.
Although we do not prescribe how to select $P$ and $Q$, one should bear in mind that we aim to construct a linear approximation of $Q$, meaning that a decomposition should be identified that renders $Q$ as close to linear as possible, to improve the quality of the approximation. 
The effect of different selections for $P$ and $Q$ is investigated in \Cref{subsec: porous}.

To  adaptively construct linear approximations to the nonlinear differential operator $Q$, we propose to exploit traditional finite difference formulae \citep{strikwerda2004finite}.
Note that our conceptualisation of these approximations as linear \textit{operators} for Gaussian conditioning is somewhat non-traditional. 
Define a time discretisation grid $\mathbf{t}=[t_0,t_1 \dots t_{n-1}]$, where $0=t_0<t_1<\dots<t_{n-1} \leq T$ with the increment $\delta \coloneqq t_i - t_{i-1}$ fixed.
For concreteness, consider Burgers' equation with $P= \partial_t - \varepsilon \partial_x^2$, $Q = u \partial_x$.
The following discussion is intended only to be informal.
Suppose that the unknown solution $u(t_{i-1},\cdot)$ at time $t_{i-1}$ has been approximated to accuracy $\mathcal{O}(\delta)$ by $u_{i-1}(\cdot)$, as quantified by a norm $\|\cdot\|$ on $C^\beta([0,T] \times \Gamma)$. 
Then we could adaptively build a linear approximation to $Q$ at time $t_i$ as
\begin{align}
Q_i u(t_i,x) \coloneqq u_{i-1}(x) \frac{\partial u}{\partial x} (t_i,x) . \label{eq: first lin Burger}
\end{align}
This provides an approximation $D_i = P + Q_i$ to the original differential operator $D$, at time $t_i$, with accuracy $\mathcal{O}(\delta)$.
To achieve higher order accuracy, we can use higher order approximations of $Q$.
For example, letting $\left. \frac{\partial u}{\partial t}\right|_{i-1}(x)$ denote an approximation to $\frac{\partial u}{\partial t}(t_{i-1},x)$,  we could take
\begin{align*}
Q_i u(t_i,x) \coloneqq \left[ u_{i-1}(x) + \delta \left. \frac{\partial u}{\partial t}\right|_{i-1}(x)  \right]  \frac{\partial u}{\partial x} (t_i,x) .
\end{align*}
The only requirement that we impose on finite difference approximations is that $Q_i$ uses (only) data that were gathered at earlier time points $t_{i-1}, t_{i-2}, \dots$, analogous to \textit{backward difference} formulae.
This is to ensure that the approximations $Q_i$ are well-defined before they are used in our method, which is described next.

Henceforth we assume that an appropriate representation $D = P + Q$ has been identified and an appropriate linear approximation to $Q$ has been selected.
The next section describes how probabilistic inference for $u$ can then proceed.

\subsection{Proposed Approach} \label{subsec: approach}

In this section we describe our proposed method.
Recall that we assume there exists a unique $u \in C^\beta([0,T] \times \Gamma)$ for which \eqref{eq: somepde} is satisfied.
Since \eqref{eq: somepde} represents an infinite number of constraints, it is not generally possible to recover $u$ exactly with a finite computational budget.
Our proposed method mirrors a general approach used to construct Bayesian PNM \citep{Cockayne2017}, in that we consider conditioning on only a finite number of the constraints in \eqref{eq: somepde} and reporting the remaining uncertainty as our posterior.
The case of nonlinear PDEs presents an additional challenge in that a subset of the constraints are nonlinear, and are therefore not amenable to exact Gaussian conditioning.
To circumvent this issue, we condition on linear approximations to the constraints following the ideas developed in \Cref{subsec: FD}.

\subsubsection{Prior Distribution}

The starting point of any Bayesian analysis is the elicitation of a suitable prior distribuition.
In our case, we minimally require a prior that is supported on $C^\beta([0,T] \times \Gamma)$, since we \textit{a priori} know that the solution $u$ to \eqref{eq: somepde} has this level of regularity.
Our approach is rooted in Gaussian conditioning and thus the regularity of Gaussian process sample paths must be analysed.
This analysis is somewhat technical and we therefore defer the discussion of prior elicitation to \Cref{sec: prior}.
For the remainder of this section we assume that a suitable Gaussian process prior $U \sim \mathcal{GP}(\mu,\Sigma)$ has been elicited.
Here $\mu \colon [0,T] \times \Gamma \to \mathbb{R}$, $\mu(t,x) \coloneqq \mathbb{E}[U(t,x)]$ is the \textit{mean function} and $\Sigma \colon ([0,T] \times \Gamma) \times ([0,T] \times \Gamma) \to \mathbb{R}$, $\Sigma((t,x),(t',x')) \coloneqq \mathbb{E}[(U(t,x) - \mu(t,x))(U(t',x') - \mu(t',x'))]$ is the \textit{covariance function}; see \cite{Rasmussen2006} for background.
The random variable notation $U$ serves to distinguish the true solution $u$ of \eqref{eq: somepde} from our probabilistic model for it.
The specific choices of $\mu$ and $\Sigma$ discussed in \Cref{sec: prior} have sufficient regularity for the subsequent derivations in this section to be well-defined.

\subsubsection{Initialisation}

At the outset we fix a time discretisation $\mathbf{t}=[t_0,t_1 \dots t_{n-1}]$, where $0=t_0<t_1<\dots<t_{n-1} \leq T$, and a spatial discretisation $\mathbf{x} = [x_1,x_2,\dots,x_m] \in (\Gamma \cup \partial \Gamma)^m$ where the $x_i$ are required to be distinct.
It will sometimes be convenient to interpret $\mathbf{x}$ as a set $\{x_1,\dots,x_n\}$; for instance we will write $\mathbf{x} \setminus \partial\Gamma$ to denote $\{x_1,\dots,x_n\} \setminus \partial\Gamma$.

Our first task is to condition on (or \textit{assimilate}) a finite number of constraints that encode the initial condition $u(0,x) = g(x)$, $x \in \Gamma$.
For this purpose we use the spatial discretisation $\mathbf{x}$, and condition on the data $U(0,x) = g(x)$ at each $x \in \mathbf{x} \setminus \partial \Gamma$.
(For example, if $\Gamma = [0,1]$ and $0 = x_1 < x_2 < \dots < x_{m-1} < x_m = 1$, then we condition on $U(0,x_i) = g(x_i)$ for $i = 2,\dots,m-1$.
The two boundary locations $x_1, x_m \in \partial \Gamma$ are excluded since these constraints are assimilated as part of the boundary condition, which will shortly be discussed.)
To perform conditioning, we use the following vectorised shorthand:
\begin{align*}
\bm{v}_i & \coloneqq (t_i , \mathbf{x}) \coloneqq [(t_i,x_1) , (t_i,x_2), \dots , (t_i,x_m)]^\top \in ([0,T] \times \Gamma)^m \\ 
U(\bm{v}_i) & \coloneqq [U(t_i,x_1),\dots,U(t_i,x_m)]^\top \in \mathbb{R}^{m \times 1} \\
g(\bm{v}_i) & \coloneqq [g(x_1),\dots,g(x_m)]^\top  \in \mathbb{R}^{m \times 1} \\
\Sigma((t,x), \bm{v}_i) & \coloneqq [ \Sigma((t,x),(t_i,x_1)) , \dots , \Sigma((t,x),(t_i,x_m)) ]  \in \mathbb{R}^{1 \times m} \\
\Sigma(\bm{v}_i, (t,x))  & \coloneqq \Sigma((t,x), \bm{v}_i)^\top  \in \mathbb{R}^{m \times 1}  \\
\Sigma(\bm{v}_i,\bm{v}_j) & \coloneqq \left[ \begin{array}{ccc} \Sigma((t_i,x_1),(t_j,x_1)) & \dots & \Sigma((t_i,x_1),(t_j,x_m)) \\ \vdots & & \vdots \\ \Sigma((t_i,x_m),(t_j,x_1)) & \dots & \Sigma((t_i,x_m),(t_j,x_m)) \end{array} \right]  \in \mathbb{R}^{m \times m} 
\end{align*}
Then let $\bm{a}_0 \coloneqq (t_0, \mathbf{x} \setminus \partial \Gamma)^\top$ denote the locations in $[0,T] \times \Gamma$ where the initial condition is to be assimilated.
At $\bm{a}_0$ we have the \textit{initial data} $\mathbf{y}^0 \coloneqq g(\bm{a}_0)$.
These initial data are assimilated into the Gaussian process model according to the standard conditioning formulae \citep[eq.\ 2.19;][]{Rasmussen2006}
\begin{align*}
U^0 & \coloneqq ( U|U(\bm{a}_0) = \mathbf{y}^0 ) \sim \mathcal{GP}(\mu^0,\Sigma^0) \\
\mu^0(r) & \coloneqq \mu(r) + \Sigma(r,\bm{a}_0)\Sigma(\bm{a}_0,\bm{a}_0)^{-1} 
(g(\bm{a}_0)-\mu(\bm{a_0})) \\
\Sigma^0(r,s) & \coloneqq \Sigma(r,s) - \Sigma(r,\bm{a}_0)\Sigma(\bm{a}_0,\bm{a}_0)^{-1} \Sigma(\bm{a}_0,s) 
\end{align*}
where $r, s \in [0,T] \times \Gamma$.


\subsubsection{Time Stepping}

Having assimilated the initial data, we now turn to the remaining constraints in \eqref{eq: somepde}.
Following traditional time-stepping algorithms, we propose to proceed iteratively, beginning at time $t_0$ and then advancing to $t_1$, $t_2$, and ultimately to $t_{n-1}$.
At each iteration $i$ we aim to condition on a finite number of constraints that encode the boundary condition $u(t_i,x) = h(t_i,x)$, $x \in \partial \Gamma$, and the differential equation itself $Du(t_i,x) = f(t_i,x)$, $x \in \Gamma$.
For this purpose we again use the spatial discretisation $\mathbf{x}$, and condition on the \textit{boundary data} $U(t_i,x) = h(t_i,x)$ at each $x \in \mathbf{x} \cap \partial \Gamma$ and the \textit{differential data} $DU(t_i,x) = f(t_i,x)$ at each $x \in \mathbf{x}$.
Since $D$ is nonlinear, there are no explicit formulae that can be used in general to assimilate the differential data, so instead we propose to condition on the approximate constraints $D_i U(t_i,x) = f(t_i,x)$, $x \in \mathbf{x}$ where $D_i = P + Q_i$ is an adaptively defined linear approximation to $D$, which will be problem-specific and chosen in line with the principles outlined in \Cref{subsec: FD}.

For a univariate function such as $\mu$ and a linear operator $L$, we denote $\mu_L(r) = (L\mu)(r)$.
For a bivariate function such as $\Sigma$, we denote $\Sigma_L(r,s) = L_r \Sigma(r,s)$, where $L_r$ denotes the action of $L$ on the $r$ argument.
In addition, we denote $\Sigma_{\bar{L}}(r,s) = L_s \Sigma(r,s)$ and we allow subscripts to be concatenated, such as $\Sigma_{L,L'} = (\Sigma_L)_{L'}$ for another linear operator $L'$.


Fix $i \in \{0,1,\dots,n-1\}$.
Let $\bm{b}_i = (t_i, \mathbf{x} \cap \partial \Gamma)$ denote the locations in $[0,T] \times \partial \Gamma$ where the boundary conditions at time $t_i$ are to be assimilated.
At $\bm{b}_i$ we have boundary data $h(\bm{b_i})$.
Correspondingly, we have differential data $f(\bm{v}_i)$ and we concatenate all data at time $i$ into a single vector $\mathbf{y}^i \coloneqq [h(\bm{b}_i)^\top , f(\bm{v}_i)^\top ]^\top$, so that $\mathbf{y}^i$ represents all the information on which (approximate) conditioning is to be performed. 
Upon assimilating these data we obtain
\begin{align*}
U^{i+1} &\coloneqq ( U^i|[U(\bm{b}_i), D_i U(\bm{v}_i)]=\mathbf{y}^i ) \sim \mathcal{GP}(\mu^{i+1},\Sigma^{i+1}) \\
\mu^{i+1}(r) & \coloneqq \mu^i(r)+ [\Sigma^i(r,\bm{b}_i),\Sigma_{\bar{D_i}}^i(r,\bm{v}_i)] A_i^{-1} \begin{bmatrix} 
h(\bm{b}_i)-\mu^{i}(\bm{b}_i) \\
f(\bm{v}_i)-\mu^{i}_{D_i}(\bm{v}_i)
\end{bmatrix} \\
\Sigma^{i+1}(r,s) & \coloneqq \Sigma^i(r,r')- [\Sigma^i(r,\bm{b}_i),\Sigma_{\bar{D_i}}^i(r,\bm{v}_i)] A_i^{-1} \begin{bmatrix} 
\Sigma^i(\bm{b}_i,s) \\
\Sigma^i_{D_i}(\bm{v}_i,s)
\end{bmatrix} \\
A_i & \coloneqq \begin{bmatrix} 
 \Sigma^i(\bm{b}_i,\bm{b}_i) & \Sigma^i_{\bar{D_i}}(\bm{b}_i,\bm{v}_i) \\
\Sigma^i_{{D_i}}(\bm{v}_i,\bm{b}_i) & \Sigma^i_{D_i\bar{D_i}}(\bm{v}_i,\bm{v}_i)
\end{bmatrix}
\end{align*}
The result of performing $n$ time steps of the algorithm just described is a Gaussian process $\mathcal{GP}(\mu^n, \Sigma^n)$, to which we associate the semantics of an (approximate) posterior in a Bayesian PNM for the solution of \eqref{eq: somepde}.

The Bayesian interpretation of $\mathcal{GP}(\mu^n, \Sigma^n)$ is reasonable since this distribution arises from the conditioning of the prior $\mathcal{GP}(\mu,\Sigma)$ on a finite number of constraints that are (approximately) satisfied by the solution $u$ of \eqref{eq: somepde}.
This is clarified in the following statement:

\begin{lemma} \label{lem: batch}
The stochastic process $U^n$ obtained above is identical to the distribution obtained when $U \sim \mathcal{GP}(\mu,K)$ is conditioned on the dataset
\[
\left[ \begin{array}{c} U(\bm{a}_0) \\ \left[U(\bm{b}_0),D_0U(\bm{v}_0)\right]^\top \\ \vdots \\ \left[U(\bm{b})_{n-1},D_{n}U(\bm{v}_{n-1})\right]^\top   \end{array} \right] = \left[ \begin{array}{c} \mathbf{y}^0 \\ \mathbf{y}^1 \\ \vdots \\ \mathbf{y}^{n-1} \end{array} \right].
\]
\end{lemma}
\begin{proof}
This follows immediately from the self-consistency property of Bayesian inference (invariance to the order in which data are conditioned), but for completeness we demonstrate their algebraic equivalence in \Cref{sec: derive sequential version pde}.
\end{proof}

\begin{remark}[Computational Complexity]
The computational cost of our algorithm is not competitive with that of a standard numerical method.\footnote{Technically, the computational complexity of our algorithm the same as that of a traditional numerical method that performs forward Euler increments in the temporal component and symmetric collocation in the spatial component \citep{fasshauer1999solving,Cockayne2016}.
However such methods are rarely used, with one factor for this being the computational cost.}
However, we are motivated by problems for which $f$, $g$ and $h$ are associated with a high computational cost, for which the auxiliary computation required to provide probabilistic uncertainty quantification is inconsequential.
Thus we merely remark that the iterative algorithm we presented is gated by the inversion of the matrix $A_i$ at the $i$\textsuperscript{th} time step, the size of which is $O(m)$, independent of $i$, and therefore the complexity of predicting the final state $u(T,\cdot)$ of the PDE by performing $n$ iterations of the above algorithm is $O(nm^3)$.
For comparison, direct Gaussian conditioning on the information in \Cref{lem: batch} would incur a higher computational cost of $O(n^3 m^3)$, but would provide the joint distribution over the solution $u(t,\cdot)$ at all times $t \in [0,T]$.
Although we do not pursue it in this paper, in the latter case the grid structure present in $\mathbf{t}$ and $\mathbf{x}$ could be exploited to mitigate the $O(n^3m^3)$ cost; for example, a compactly supported covariance model $\Sigma$ would reduce the cost by a constant factor \citep{gneiting2002compactly}, or if the preconditions of \cite{schafer2017compression} are satisfied then their approach would reduce the cost to $O(nm \log(nm) \log^{d+1}(nm/\epsilon))$ at the expense of introducing an error of $O(\epsilon)$.
See also the recent work of \cite{de2021high}.
\end{remark}

\begin{remark}
The posterior mean $\mu^{i+1}$ can be interpreted as a particular instance of a radial basis method \citep{Fornberg2015}, as a consequence of the representer theorem for kernel interpolants \citep{scholkopf2001generalized}.
For brevity we do not explore this connection further, but we note that a similar connection was explored in detail in \cite{Cockayne2016}.
\end{remark}

\subsubsection{Calibration of Uncertainty}

The principal advantage of PNM over classical numerical methods is that they provide probabilistic quantification of uncertainty, in our case expressed in the Bayesian framework, which can be integrated along with other sources of uncertainty to facilitate inferences and decision-making in a real-world context.
In order for our posterior distribution to faithfully reflect the scale of uncertainty about the solution of \eqref{eq: somepde}, we must allow the hyper-parameters of the prior model to adapt to the dataset.
However, we do not wish to sacrifice the sequential nature of our algorithm and thus we seek an approach to hyper-parameter estimation that operates in real-time as the algorithm is performed.

To achieve this we focus on a covariance model $\Sigma(\cdot,\cdot;\sigma)$ with a scalar hyper-parameter denoted $\sigma > 0$, which is assumed to satisfy $\Sigma(\cdot,\cdot;\sigma)=\sigma^2 \Sigma(\cdot,\cdot;1)$.
Such a $\sigma$ is sometimes called a \textit{scale} or \textit{amplitude} hyper-parameter of the covariance model.
From \Cref{lem: batch} it follows that $\sigma$ directly controls the spread of the posterior and it is therefore essential that $\sigma$ is estimated from data in order that the uncertainty reported by the posterior can be meaningful.
To estimate $\sigma$, we propose to maximise the predictive likelihood of the ``differential data'' $f(\bm{v}_i)$, given the information collected up to iteration $i-1$, for $i \in \{0,\dots,n-1\}$, which can be considered as an \textit{empirical Bayes} approach based on just those factors in the likelihood that correspond to the differential data.
The reasons for focussing on the differential data (as opposed to also including the initial and boundary data) are twofold; first, the differential data constitutes the vast majority of the dataset, and second, this simplifies the computational implementation, described next.

At iteration $i$, the predictive likelihood for $U_{D_i}(\bm{v}_i)$ is $\mathcal{N}( \mu^{i}_{D_i}(\bm{v}_i),\Sigma^{i}(\bm{v}_i,\bm{v}_i;\sigma))$, and the observed differential data are $f(\bm{v}_i)$. 
Thus we select $\sigma$ to maximise the full predictive likelihood of the differential data
\begin{equation}
\prod_{i=0}^{n-1} \mathcal{N}(f(\bm{v}_i); \mu^{i}_{D_i}(\bm{v}_i),\Sigma^{i}_{D_i}(\bm{v}_i,\bm{v}_i;\sigma)) . \label{eq: pred lik}
\end{equation}
Crucially, the linear operators $D_i$ that we constructed do not directly depend on $\sigma$, and it is a standard property of Gaussian conditioning that $\Sigma_{D_i}^i(\bm{v}_i,\bm{v}_i ; \sigma) = \sigma^2 \Sigma_{D_i}(\bm{v}_i,\bm{v}_i;1)$.
These facts permit a simple closed form expression for the maximiser $\hat{\sigma}$ of \eqref{eq: pred lik}, namely
\begin{equation}
\hat{\sigma}^2 = \frac{1}{n} \sum_{i=0}^{n-1}  \left\| (\Sigma^{i}_{D_i}(\bm{v}_i,\bm{v}_i);1)^{-\frac{1}{2}}(f(\bm{v}_i)- \mu^{i}_{D_i}(\bm{v}_i))  \right\|^2 \label{eq: hat sigma}
\end{equation}
where $M^{-1/2}$ denotes an inverse matrix square root; $(M^{1/2})^2 = M$.
Furthermore, it is clear from \eqref{eq: hat sigma} that in practice one can simply run our proposed algorithm with the prior covariance model $\Sigma(\cdot,\cdot ; 1)$ and then report the posterior covariance $\hat{\sigma}^2 \Sigma^{i+1}(\cdot,\cdot ; 1)$, so that hyper-parameter estimation is performed in real-time without sacrificing the sequential nature of the algorithm.

Closed form expressions such as \eqref{eq: hat sigma} are not typically available for other hyper-parameters that may be included in the covariance model, and we therefore assume in the sequel that any other hyper-parameters have been expert-elicited.
This limits the applicability of our method to situations where some prior expert insight can be provided.
However, we note that data-driven estimation of the amplitude parameter $\sigma$ is able to compensate to a degree for mis-specification of other parameters in the covariance model.

\subsubsection{Relation to Earlier Work}

Here we summarise how the method just proposed relates to existing literature on Bayesian PNM and beyond.

The sequential updating procedure that we have proposed is similar to that of \cite{Chkrebtii2016} in the special case of a linear PDE.
It is not identical in these circumstances though, for two reasons:
First, \cite{Chkrebtii2016} incorporated the initial condition $u(0,x) = g(x)$, $x \in \Gamma$, into the prior model, whereas we explicitly conditioned on initial data $g(\bm{a}_0)$ during the initialisation step of the method.
This direct encoding of the initial condition in \cite{Chkrebtii2016} relies on $g$ being analytically tractable in order that a suitable prior can be derived by hand.
Our treatment of $g$ as a black-box function from which initial data are provided is therefore more general.
Second, in \cite{Chkrebtii2016} the authors advocated the use of an explicit \textit{measurement error} model, whereas our conditioning formula assume that the differential data $\mathbf{y}^i$ are exact measurements of $U$, as clarified in \Cref{lem: batch}.
For linear PDEs this assumption is correct, but it is an approximation in the case a nonlinear PDE.
Our decision not to employ a measurement error model here is due to the fact that the scale of the measurement error cannot easily be estimated in an online manner as part of a sequential algorithm, without further approximations being introduced.

To limit scope, the adaptive selection of the $t_i$ or $x_j$ was not considered, but we refer the reader to \cite{Chkrebtii2019} for an example of how this can be achieved using Bayesian PNM.
Note, however, that adaptive selection of a time grid may be problematic when evaluation of either $f$ or $h$ is associated with a high computational cost, since the possibility of taking many small time steps relinquishes control of the computational budget.
For this reason, non-adaptive methods may be preferred in this context, since the run-time of the PNM can be provided up-front.

The choice of linearisation $Q_i$ was left as an input to the proposed method, with some guidelines (only) provided in \Cref{subsec: FD}.
This can be contrasted with recent work for ODEs in \cite{Tronarp2018,Tronarp2020,Bosch2020}, where first order Taylor series were used to automatically linearise a nonlinear gradient field.
It would be possible to also consider the use of Taylor series methods for nonlinear PDEs.
However, their use assumes that the gradient field is analytically tractable and can be differentiated, while in the present paper we are motivated by situations in which $f$ is a black-box that can (only) be point-wise evaluated.
The use of linearisations in PNM was also explored \citet{chen2021solving}, in the \textit{maximum a posteriori} estimation context.

The combination of local linearisation and Gaussian process conditioning was also studied by \cite{Raissi2018}, who considered initial value problems specified by PDEs, where the initial condition was random and the goal was to approximate the implied distribution over the solution space of the PDE.
The authors observed that if the initial condition was a Gaussian process, then approximate conjugate Gaussian computation is possible when a finite difference approximation to the differential operator was employed.
This provided a one-pass, cost-efficient alternative to the Monte Carlo approach of repeatedly sampling an initial condition and then applying a classical numerical method.
Our work bears a superficial similarity to \cite{Raissi2018} and related work on \textit{physics-informed} Gaussian process regression \citep[e.g.][]{wheeler2014mechanistic,Wang2016,jidling2017linearly,Chen2020}, in that finite difference approximations enable approximate Gaussian conditioning to be performed.
However, these authors are addressing a \textit{fundamentally different} problem to that addressed in the present paper; we aim to quantify numerical uncertainty for a single (i.e.\ non-random) PDE.
Accordingly, in this paper we emphasise issues that are critical to the performance of PNM, such as ensuring that the posterior is supported on a set of functions whose regularity matches that known to be possessed by the solution of the PDE (\Cref{sec: prior}), and explicitly assessing the quality of the credible sets provided by our PNM (\Cref{sec: results}).

%% file: theory.tex
\section{Prior Construction} \label{sec: prior}

This section is dedicated to presenting a prior construction that ensures samples generated from the prior are elements of $C^\beta([0,T] \times \Gamma)$, the set in which a solution to \eqref{eq: somepde} is sought.
First, in \Cref{subsec: prior properties} we introduce the technical notions of \textit{sample continuity} and \textit{sample differentiability}, clarifying what properties of the prior are required to hold.
These sample-path properties are distinct from mean-square properties, the latter being more commonly studied.
Then, in \Cref{subsec: Matern} we formally prove that the required properties holds for a particular Mat\'{e}rn tensor product, which we then advocate as a default choice for our PNM.
These results may also be of independent interest.

\subsection{Mathematical Properties for the Prior} \label{subsec: prior properties}

This paper is concerned with the strong solution of \eqref{eq: somepde}, which is an element of $C^\beta([0,T] \times \Gamma)$.
It is therefore logical to construct a prior distribution whose samples also belong to this set. 
In particular, if the true solution has $\beta_i$ derivatives in the variable $z_i$ (for instance because the PDE features a term $\partial_{z_i}^{\beta_i}u$, where we have set $z = (t,x)$), it would be appropriate to construct a prior (and hence a posterior) whose samples also have $\beta_i$ derivatives in the variable $z_i$. 

To make this discussion precise, we make explicit a probability space $(\Omega, \mathcal{F}, \mathbb{P})$ and recall the fundamental definitions of sample continuity and sample differentiability for a random field $X \colon I \times \Omega \to \mathbb{R}$ defined on an open, pathwise-connected set $I \subseteq \mathbb{R}^d$ (i.e.\ $I$ is an interval when $d=1$):

\begin{definition}[Sample Continuity]
$X$ is said to be \emph{sample continuous} if, for $\mathbb{P}$-almost all $\omega \in \Omega$, the sample path $X(\cdot,\omega)$ is continuous (everywhere) in $I$.
\end{definition}


\begin{definition}[Sample Differentiability]
Let $v^1,\dots,v^p \in \mathbb{R}^d$ be a sequence of directions and $v=(v^1,\dots,v^p)$. Then $X$ is said to be \emph{sample partial differentiable in the sequence of directions $v$} if for $\mathbb{P}$-almost all $\omega \in \Omega$, the following limit exists for all $z \in I$
\[
\mathcal{D}^pX(z,v,\omega)=\lim_{h_1 \to 0} \dots \lim_{h_p \to 0} \frac{\Delta^pX(z,v,h,\omega)}{\prod_{i=1}^p h_i}  < \infty
\]
where
\[
\Delta^pX(z,v,h,\omega) \coloneqq \sum_{r \in \{0,1\}^p} (-1)^{p - \sum_{i=1}^p r_i} X\left( z + \sum_{i=1}^p r_i h_i v^i ,\omega \right) .
\]
\end{definition}

\noindent The limits above are taken sequentially from left to right. In the discussions that follow, we take $v^i \in \{e_1,e_2,\dots,e_d\}$, the standard Cartesian unit basis vectors of $\mathbb{R}^d$, in which case the usual partial derivatives are retrieved, and we use the shorthand $\mathcal{D}^pX(z,v,\omega)=\partial^{\alpha}X(z,\omega)$ to denote sample partial derivatives, where $\alpha = (\alpha_1,\dots,\alpha_p)$, $|\alpha| =p$, and $\alpha_i$ denotes the number of times the variable $z_i$ is differentiated. 



A similar property, which is more easily studied than sample continuity (resp.\ sample differentiability), is mean-square continuity (resp.\ mean-square differentiability).
This property is recalled next, since we will make use of mean-square properties en route to establishing sample path properties in  \Cref{subsec: Matern}.

\begin{definition}[Mean-Square Continuity]
$X$ is said to be \emph{mean-square continuous at $z \in I$} if
$$
\mathbb{E}\left[ X(z,\omega)^2 \right] < \infty, \qquad \lim_{z' \to z} \mathbb{E} \left[ (X(z',\omega)-X(z,\omega))^2 \right] = 0 . 
$$
\end{definition}

\begin{definition}[Mean-Square Differentiability]
Let $v^1,\dots,v^p \in \mathbb{R}^d$ be a sequence of directions and $v=(v^1,\dots,v^p)$. Then $X$ is said to be \emph{mean-square partial differentiable at $z \in I$ in the sequence of directions $v$} if there exists a finite random field $\omega \mapsto \mathcal{D}_{\textsc{ms}}^p X(z,v,\omega)$ such that
\[
\lim_{h_1 \to 0} \dots \lim_{h_p \to 0} \mathbb{E}\left[ \left(\frac{\Delta^p X(z,v,h,\omega)}{\prod_{i=1}^p h_i} -\mathcal{D}_{\textsc{ms}}^p X(z,v,\omega)\right)^2\right] =0
\]
is well-defined.
\end{definition}

For a mean-square differentiable Gaussian processes $X$, with mean function $\mu \in C^\alpha(I)$ and covariance function $\Sigma \in C^{(\alpha,\alpha)}(I \times I)$, one has
\begin{align*}
\partial_{\textsc{ms}}^\alpha X \sim \mathcal{GP}(\partial^\alpha \mu , \partial^\alpha \bar{\partial}^\alpha \Sigma )
\end{align*}
where we use the shorthand 
$\mathcal{D}^p_{\textsc{ms}}X(z,v,\omega)=\partial_{\textsc{ms}}^{\alpha}X(z,\omega)$ to denote mean-square partial derivatives, where again the $v^i$ are unit vectors parallel to the coordinate axes, $|\alpha| =p$, and $\alpha_i$ denotes the number of times the variable $z_i$ is differentiated.\footnote{The shorthand notation surpresses the order in which derivatives are taken, and can therefore only be applied in situations where partial derivatives are continuous, to ensure that their order can be interchanged without affecting the result.
In the following, the notation $\partial_{\textsc{ms}}^\alpha X$ is used only for Gaussian processes with $\partial^\alpha \bar{\partial}^\alpha \Sigma \in C(I \times I)$.}
See \citet[][Section 2.6]{Stein1999}.
If $X$ is mean-square continuous (resp.\ mean-square differentiable for all $\alpha \leq \beta$) at all $z \in I$, then we say simply that $X$ is \emph{mean-square continuous} (resp.\ \emph{order $\beta$ mean-square differentiable} ).
In contrast to sample path properties, mean-square properties are often straight-forward to establish.
In particular, if $X$ is weakly stationary with autocovariance function $\Sigma(z) = \Sigma(z,0)$, then
\begin{equation}
\mathbb{E}\left[ (X(z,\omega)-X(z',\omega))^2 \right]=2(\Sigma(0)-\Sigma(z-z')) , 
\label{eq: stationarycontinuity}
\end{equation}
meaning that $X$ is mean-square continuous whenever its autocovariance function $\Sigma$ is continuous at $0$ \citep[][Section 2.4]{Stein1999}.

\subsection{Mat\'{e}rn Tensor Product} \label{subsec: Matern}

Our aim in this section is to establish sample path properties for a particular choice of prior, namely a Gaussian process with tensor product Mat\'{e}rn covariance, to ensure that prior (and posterior) samples are contained in $C^\beta([0,T] \times \Gamma)$.
Surprisingly, we are unable to find explicit results in the literature for the sample path properties of commonly used covariance models; this is likely due to the comparative technical difficulty in establishing sample path properties compared to mean-square properties.
Our aim in this section is, first, to furnish a gap in the literature by rigorously establishing the sample differentiability properties of Gaussian processes defined by the Mat\'{e}rn covariance function and, second, to put forward a default prior for our PNM that takes values in $C^\beta([0,T] \times \Gamma)$.

\begin{definition}[Mat\'{e}rn Covariance]
Let $\nu=p+\frac{1}{2}$ where $p \in  \mathbb{N}$.
The Mat\'{e}rn covariance function is defined, for $z,z' \in \mathbb{R}$, as
\begin{equation}
K_{\nu}(z,z')=K_{\nu}(z-z')=\sigma^2 \exp \left( -\frac{|z-z'|}{\rho} \right) \frac{p!}{(2p)!} \sum_{k=0}^{p} \frac{(2p-k)!}{(p-k)!k!} \left( \frac{2}{\rho} \right)^{k}|z-z'|^{k} .
\label{eq: Maternkernel}
\end{equation}
\end{definition}


The proof of the following result can be found in \Cref{app: Maternmsdifferentiable}.

\begin{proposition}[Mean-Square Differentiability of Univariate Mat\'{e}rn]\label{thm: Maternmsdifferentiable}
Let $I \subseteq \mathbb{R}$ be an open set and let $\mu \in C^p(I)$.
Then any process $X \sim \mathcal{GP}(\mu,K_{\nu})$, with $K_\nu$ as in \eqref{eq: Maternkernel} with $\nu = p + \frac{1}{2}$, is order $p$ mean-square differentiable. 
Furthermore, $\partial_{\textsc{ms}}^{p}X$ is mean-square continuous.
\end{proposition}



Following a general approach outlined in \cite{Potthoff2010}, and focussing initially on the univariate case, our first step toward establishing sample differentiability is to establish sample continuity of the mean-square derivatives.
Recall that, for two stochastic processes $X,\tilde{X}$ on a domain $I$, we say $\tilde{X}$ is a \emph{modification of $X$} if, for every $z \in I$, $\mathbb{P}(X(z,\omega)=\tilde{X}(z,\omega))=1$. 
A modification of a stochastic process does not change its mean square properties, but sample path properties need not be invariant to modification.\footnote{To build intuition into the role of modifications, let $X \colon [0,1] \times \Omega \to \mathbb{R}$ be a sample continuous stochastic process and consider the process $\tilde{X}(z,\omega) \coloneqq X(z,\omega) + 1[z =Z(\omega)]$ where $Z \sim \mathcal{U}(0,1)$, independent of $X$. Then $\tilde{X}$ is a modification of $X$ whose finite dimensional distributions (and hence mean square properties) are identical to those of $X$, but $\tilde{X}$ is almost surely \textit{not} sample continuous. In such circumstances it is convenient (and standard practice) to work with the sample continuous process, $X$, as opposed to $\tilde{X}$. }
For Gaussian processes, which are characterised up to modifications by their finite dimensional distributions, it is standard practice to work with continuous modifications when they exist \citep[see for example][]{DUDLEY1967290,Marcus1972}.
The proof of the following result can be found in \Cref{app: Maternderivativesamplecont}.

\begin{proposition}\label{thm: Maternderivativesamplecont}
Let $X$ be as in \Cref{thm: Maternmsdifferentiable}.
Then $\partial_{\textsc{ms}}^{i}X$ has a modification that is sample continuous for all $0\leq i \leq p$.
\end{proposition}

The second step is to leverage a fundamental result on the sample path properties of stochastic processes.

\begin{theorem}[Criterion for Sample Differentiability; Theorem 3.2 of \cite{Potthoff2010}]\label{thm: Sampledifferentiabilitycriterion}
Let $I \subseteq \mathbb{R}^d$ be an open, pathwise connected set, and consider a random field $X \colon I \times \Omega \to \mathbb{R}$ such that $\mathbb{E}[X(z,\omega)^2] < \infty$ for all $z \in I$. 
Suppose $X$ is first order mean-square differentiable, with mean-square partial derivatives $\mathcal{D}_{\textsc{ms}}^1 X(\cdot,e_k,\omega)$, $1\leq k \leq d$, themselves being mean-square continuous and having modifications that are sample continuous.
Then $X$ has a modification $\tilde{X}$ that is first order sample partial differentiable, with partial derivatives $\mathcal{D}^1  \tilde{X}(\cdot,e_k,\omega)$, $1\leq k \leq d$, themselves being sample continuous and satisfying, almost surely, $\mathcal{D}^1  \tilde{X}(\cdot,e_k,\omega) =\mathcal{D}_{\textsc{ms}}^1 X(\cdot,e_k,\omega)$, $1\leq k \leq d$.

\end{theorem}

\noindent Since continuity of partial derivatives implies differentiability, the conclusion of \Cref{thm: Sampledifferentiabilitycriterion} implies that $X$ is first order sample differentiable.

Iterative application of \Cref{thm: Sampledifferentiabilitycriterion} to higher order derivatives provides the following, whose proof can be found in \Cref{app: Sampledifferentiabilitynthordercriterion}.

\begin{corollary}\label{thm: Sampledifferentiabilitynthordercriterion}



Fix $p \in \mathbb{N}$.
Let $I \subseteq \mathbb{R}^d$ be an open, pathwise connected set and consider $X \sim \mathcal{GP}(\mu,\Sigma)$ with $\mu \in C^{p}(I)$ and $\Sigma \in C^{(p,p)}(I \times I)$, so that $X$ has mean-square partial derivatives $\partial_{\textsc{ms}}^\beta X$, $\beta  \in \mathbb{N}_0^d, |\beta| \leq p$.
Suppose $\partial_{\textsc{ms}}^\beta X$ is mean-square continuous and sample continuous for all $|\beta| \leq p$. 
Then $X$ has continuous sample partial derivatives $\partial^\beta X$, and they satisfy $\partial^\beta X = \partial_{\textsc{ms}}^\beta X$ almost surely, for all $|\beta| \leq p$.

\end{corollary}

\noindent This provides a strategy to establish sample properties of Mat\'{e}rn processes, such as the following:

\begin{corollary}[Sample Differentiability of Univariate Mat\'{e}rn] 
\label{cor: diff univ mat}
Let $I \subseteq \mathbb{R}$ be an open interval and let $X \sim \mathcal{GP}(\mu,K_{\nu})$, with $\mu \in C^p(I)$ and $K_\nu$ as in \eqref{eq: Maternkernel}.
Then there exists a modification $\tilde{X}$ of $X$ such that $\mathbb{P}(\tilde{X} \in C^p(I)) = 1$.
\end{corollary}

\begin{proof}
By \Cref{thm: Maternmsdifferentiable}, $X$ is order $p$ mean-square differentiable and $\partial_{\textsc{ms}}^{i}X$ is mean-square continuous for $0\leq i \leq p$. 
By \Cref{thm: Maternderivativesamplecont}, we may work with a modification $\tilde{X}$ of $X$ such that $\partial_{\textsc{ms}}^{i}\tilde{X} (= \partial_{\textsc{ms}}^{i} X)$ is sample continuous for each $0\leq i \leq p$. 
One can directly verify that $K_\nu \in C^{(p,p)}(I \times I)$; see the calculations in \Cref{app: Maternmsdifferentiable}. 
The result then follows from \Cref{thm: Sampledifferentiabilitynthordercriterion}.
\end{proof}


\noindent \Cref{cor: diff univ mat} is stronger than existing results in the literature, the most relevant of which is \citet[][Theorem 5]{Scheuerer2010}, who showed that samples from $\mathcal{GP}(\mu,K_{\nu + \epsilon})$ are $C^p(I)$ for any $\epsilon > 0$.

Finally, we present a multivariate version of the previous result, which will be exploited in the experiments that we perform in \Cref{sec: results}.
Importantly, we allow for different smoothness in the different variables, which is necessary to properly capture the regularity of solutions to PDEs.
The proof of this result is contained in \Cref{app: product of Matern}.

\begin{theorem}[Sample Differentiability of Mat\'{e}rn Tensor Product] \label{thm: product of Matern}


Let $I = (a_1,b_1) \times \dots \times (a_d,b_d)$ be a bounded hyper-rectangle in $\mathbb{R}^d$.
Fix $\beta \in \mathbb{N}_0^d$.
Let $\mu \in C^\beta(I)$ be bounded\footnote{The requirement that $\mu$ be bounded, which did not appear in the corresponding univariate result (\Cref{cor: diff univ mat}) is likely an artefact of our proof strategy, rather than a necessary condition for \Cref{thm: product of Matern} to hold.} in $I$, and consider a covariance function $\Sigma \colon I \times I \to \mathbb{R}$ of the form
\begin{equation*}
\Sigma(z , z') = \prod_{i=1}^{d} K_{\nu_i}(z_i-z_i')
\label{eq: Maternproductcovariance}
\end{equation*}
where $z=(z_1,z_2,\dots,z_d)$, $z'=(z_1',z_2',\dots,z_d')$ and $\nu_i = \beta_i+\frac{1}{2}$ for each $i = 1,\dots,d$. 
Then a Gaussian process of the form $X \sim \mathcal{GP}(\mu,\Sigma)$ has a modification $\tilde{X}$ that satisfies $\mathbb{P}(\tilde{X} \in C^\beta(I)) = 1$.
\end{theorem}

%% file: results.tex
\section{Experimental Assessment} \label{sec: results}

In this section, proof-of-concept numerical studies of three different initial value problems are presented. 
The first and simplest case is a homogeneous Burger's equation, a PDE with one nonlinear term and a solution that is known to be smooth. 
The second case is a porous medium equation, with two nonlinear terms appearing in the PDE and a solution that is known to be piecewise smooth, so that a classical solution does not exist and our modelling assumptions are violated. 
The third case returns to Burger's equation but now with forcing, to simulate a scenario where the right hand side $f$ is a black box function that may be evaluated at a high computational cost. 
All three experiments are \textit{synthetic}, in the sense that the functions $f$, $g$ and $h$, which in our motivating task are considered to be black boxes associated with a high computational cost, are in actual fact simple analytic expressions, enabling a thorough empirical assessment to be performed.

In order to assess the empirical performance of our algorithm, two distinct performance measures were employed.
The first of these aims to assess the accuracy of the posterior mean, which is analogous to how classical numerical methods are assessed.
For this purpose the $L^\infty$ error was considered:
\begin{equation}\label{eq: linfinityerror}
E_\infty \coloneqq \sup_{ t \in [0,T] , x \in \Gamma } \left| \mu^n(t,x) - u(t,x) \right|
\end{equation}
In practice the value of \eqref{eq: linfinityerror} is approximated by taking the maximum over the grid $\mathbf{t} \times \mathbf{x}$ on which the data $\mathbf{y}^0, \dots, \mathbf{y}^{n-1}$ were obtained. 
Accuracy that is comparable to a classical numerical method is of course desirable, but it is not our goal to compete with classical numerical methods in terms of $L^\infty$ error.
The second statistic that we consider assesses whether the distributional output from our PNM is \textit{calibrated}, in the sense that the scale of the Gaussian posterior is comparable with the difference between the posterior mean $\mu^n$ and the true solution $u$ of \eqref{eq: somepde}:
\begin{equation}\label{eq: zscore}
Z \coloneqq \sup_{ t \in [0,T] , x \in \Gamma }  \cfrac{\left| \mu^n(t,x) -u(t,x)\right|}{ \hat{\sigma} \Sigma^n(t,x)^{1/2} }
\end{equation}
This performance measure will be called a \textit{$Z$-score}, in analogy with traditional terminology from statistics.
For the purpose of this exploratory work, values of $Z$ that are orders of magnitude smaller than 1 are interpreted as indicating that the distributional output from the PNM is under-confident, while values that are orders of magnitude greater than 1 indicate that the PNM is over-confident.
A PNM that is neither under nor over confident is said to be \textit{calibrated} \citep[precise definitions of the term ``calibrated'' can be found in][but the results we present are straight-forward to interpret using the informal approach just described]{Karvonen2020,Cockayne2021}.
Our goal in this work is to develop an approximately Bayesian PNM for nonlinear PDEs that is both accurate and calibrated.
Again, in practice the supremum in \eqref{eq: zscore} is approximated by the maximum over the $\mathbf{t} \times \mathbf{x}$ grid.

For all experiments below, we consider uniform temporal and spatial grids of respective sizes $n = 2^i + 1$, $m = 2^j + 1$, where $i,j \in \{2,3,4,5,6,7\}$. 
This ensures that the grid points at which data are obtained are strictly nested as either the temporal exponent $i$ or the spatial exponent $j$ are increased. 
The prior mean $\mu(t,x) = 0$ for all $t \in [0,T]$, $x \in \Gamma$, will be used throughout.

\subsection{Homogeneous Burger's Equation} \label{subsec: homog burger}

Our first example is the homogeneous \textit{Burger's equation}
\begin{equation*}
\frac{\partial u}{\partial t} + u\frac{\partial u}{\partial x} -\alpha\frac{\partial^2 u}{\partial x^2}=0, \qquad t \in [0,T], \; x \in [0,L] \label{eq: Burgerhomogeneous}
\end{equation*}
with initial and boundary conditions
\begin{equation*}
\begin{alignedat}{2}
u(0,x) & = 2 \alpha \Big( \frac{a k \mathrm{sin}(kx)}{b + a \mathrm{cos}(kx)} \Big), && \qquad x \in [0,L] \\
u(t,0) = u(t,2\pi) & = 0, && \qquad t \in [0,T]
\end{alignedat}
\end{equation*}
and, for our experiments, $\alpha=0.02$, $a=1$, $b=2$, $k=1$, $T = 30$ and $L = 2 \pi$. 
These initial and boundary conditions were chosen because they permit a closed-form solution 
\begin{equation} \label{eq: Burgeranalytic}
u(t,x)=2 \alpha \Big( \frac{ ak \mathrm{exp}(-\alpha k^2 t)\mathrm{sin}(kx)}{b + a \mathrm{exp}(-\alpha k^2 t) \mathrm{cos}(kx)}  \Big)
\end{equation}
that can be used as a ground truth for our assessment.

To linearise the differential operator Burger's equation we consider approximations of the form in \eqref{eq: first lin Burger}, i.e.
\begin{align*}
Q_i u(t_i,x) \coloneqq u_{i-1}(x) \frac{\partial u}{\partial x} (t_i,x) 
\end{align*}
where $u_{i-1}(x)$ was taken equal to the predictive mean $\mu^{i-1}(t_i,x)$ arising from the Gaussian process approximation $U^{i-1}$.

\paragraph{Default Prior:}

\begin{figure}[t!]
\centering
\includegraphics[width = 0.75\textwidth]{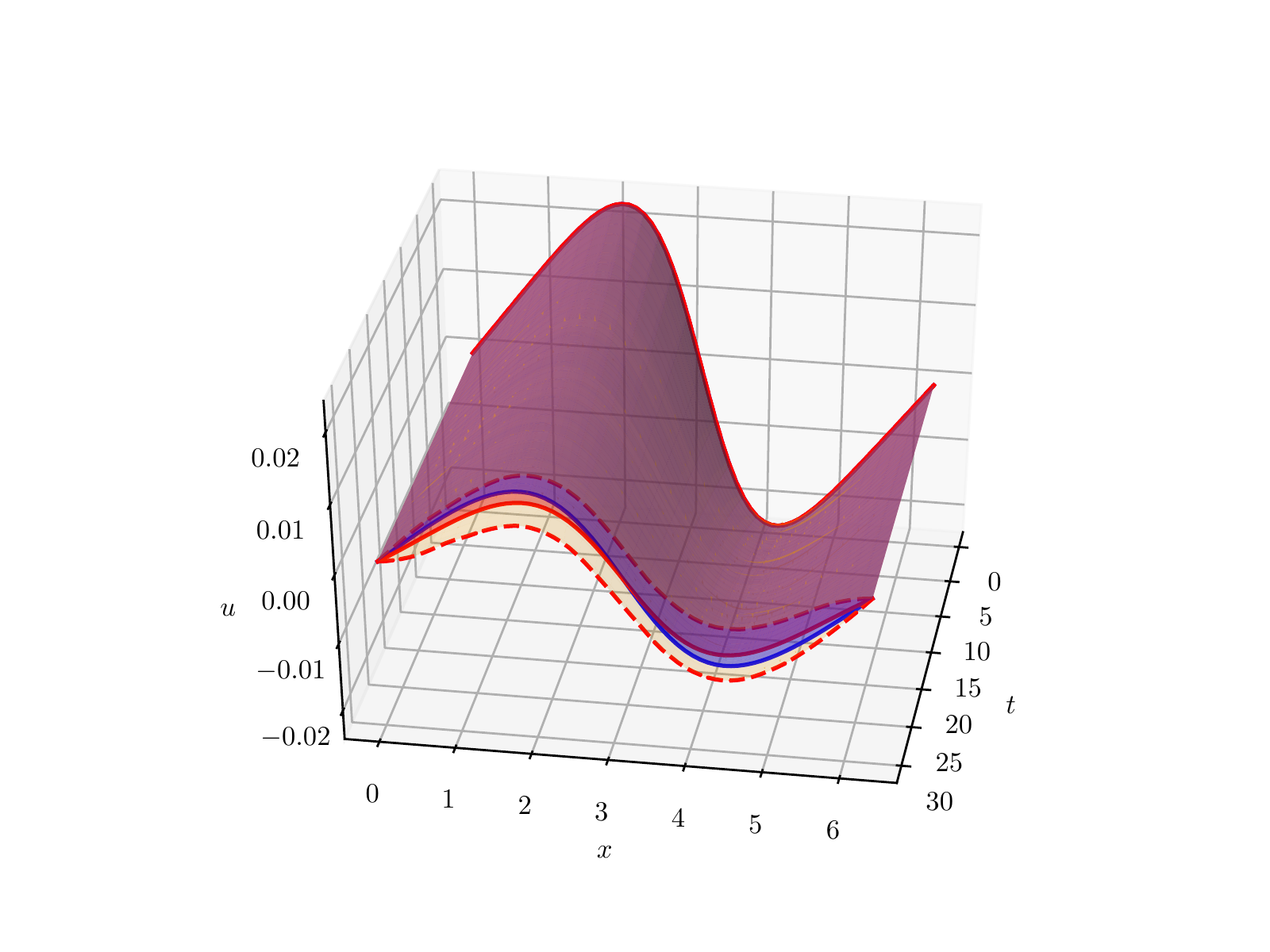}
\caption{Homogeneous Burger's equation: 
For each point $(t,x)$ in the domain we plot:
the analytic solution $u(t,x)$ (blue), the posterior mean $\mu^n(t,x)$ (red) from the proposed probabilistic numerical method, and 0.025 and 0.975 quantiles 
of the posterior distribution at each point (orange).
Here the default prior was used, with a spatial grid of size $m = 65$ and a temporal grid of size $n = 65$.
}
\label{fig: burgernoforce3d}
\end{figure}

Burger's equation has a first order temporal derivative term and a second order spatial derivative term, so following the discussion in \Cref{sec: prior} we consider as a default a Gaussian process prior with covariance function $\Sigma$ that is a product between a Mat\'{e}rn $3/2$ kernel $K_{3/2}(t,t')$ for the temporal component, and a Mat\'{e}rn $5/2$ kernel $K_{5/2}(x,x')$ for the spatial component:
\begin{equation}\label{eq: Maternfullkernel}
\Sigma((t,x),(t',x'))=K_{3/2}(t,t' ; \rho_1, \sigma_1) K_{5/2}(x,x' ; \rho_2 , \sigma_2)
\end{equation} 
The notation in \eqref{eq: Maternfullkernel} makes explicit the dependence on the amplitude hyper-parameters $\sigma_1, \sigma_2$ and the length-scale hyper-parameters $\rho_1, \rho_2$; note that only the product $\sigma \coloneqq \sigma_1 \sigma_2$ of the two amplitide parameters is required to be specified.
For the experiments below $\sigma$ was estimated as per \eqref{eq: hat sigma}, while the length-scale parameters were fixed at values $\rho_1=6$, $\rho_2=3$ (not optimised; these were selected based on a \emph{post-hoc} visual check of the credible sets in \Cref{fig: burgernoforce3d}). 
From \Cref{thm: product of Matern}, this construction ensures that the prior is supported on $C^{(1,2)}([0,T] \times [0,L])$.
Typical output from our PNM equipped with the default prior is presented in \Cref{fig: burgernoforce3d}.

\paragraph{An Alternative Prior:}

The Mat\'{e}rn covariance models assume only the minimal amount of smoothness required for the PDE to be well-defined.
However, in this assessment the ground truth $u$ is available \eqref{eq: Burgeranalytic} and is seen to be infinitely differentiable in $(0,T] \times [0,2\pi]$.
It is therefore interesting to explore whether a prior that encodes additional smoothness can improve on the default prior in \eqref{eq: Maternfullkernel}.
A prototypical example of such a prior is 
\begin{equation*}
\Sigma((t,x),(t',x')) = C(t,t' ; \rho_3, \sigma_3) C(x,x' ; \rho_4, \sigma_4)
\end{equation*}
where 
\begin{equation*}\label{eq: Rationalquadraticx}
C(z,z' ; \rho, \sigma) \coloneqq \sigma \left( 1 + \frac{(z-z')^2}{\rho^2} \right)^{-1}
\end{equation*}
is the \textit{rational quadratic} covariance model.
For the experiments below $\sigma$ was estimated as per \eqref{eq: hat sigma}, while the length-scale parameters were fixed at values $\rho_3=\sqrt{3}$, $\rho_4=\sqrt{3}$.

\begin{figure}[t!]
\centering
\includegraphics[width = 0.48\textwidth]{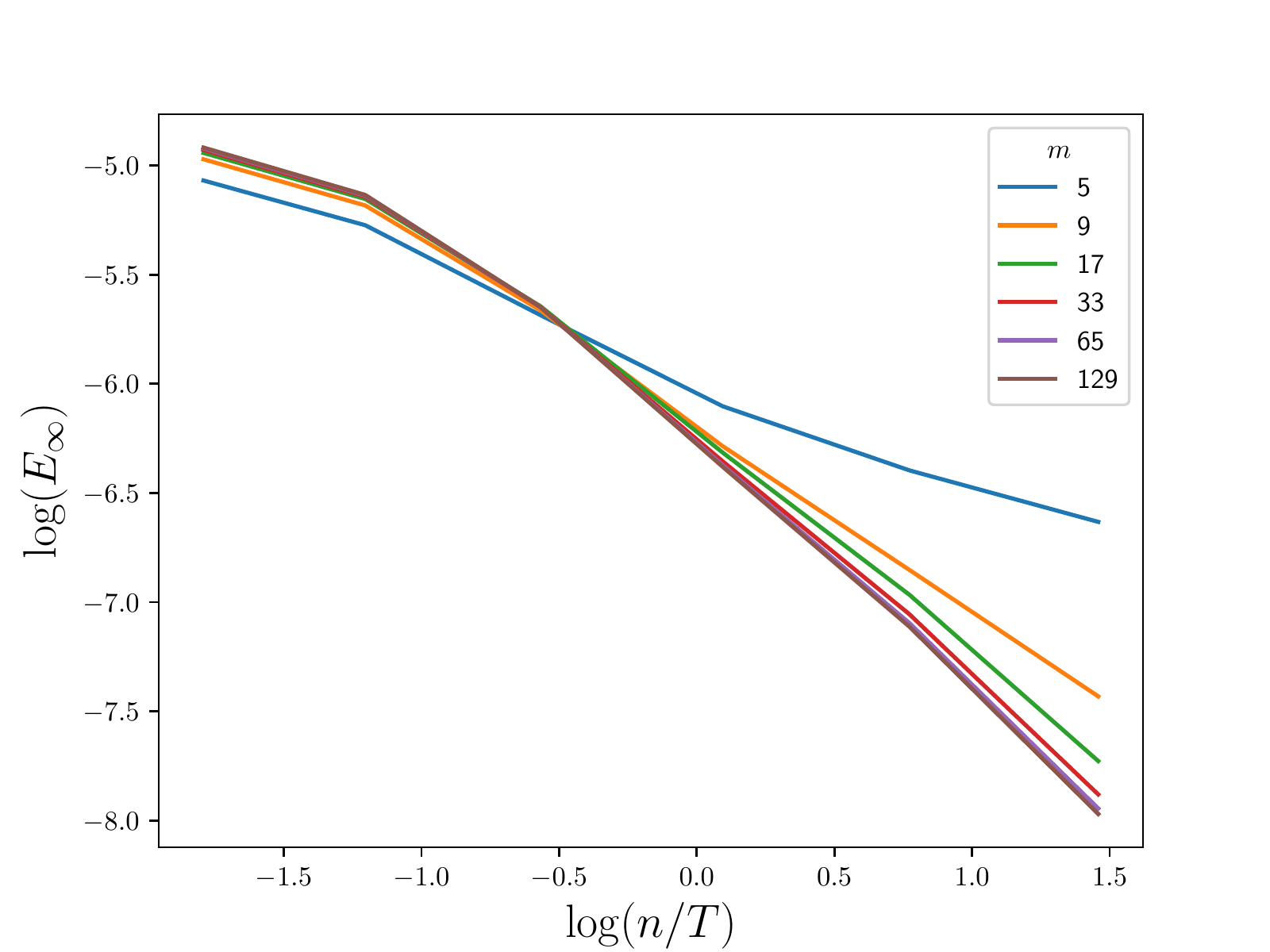}
\includegraphics[width = 0.48\textwidth]{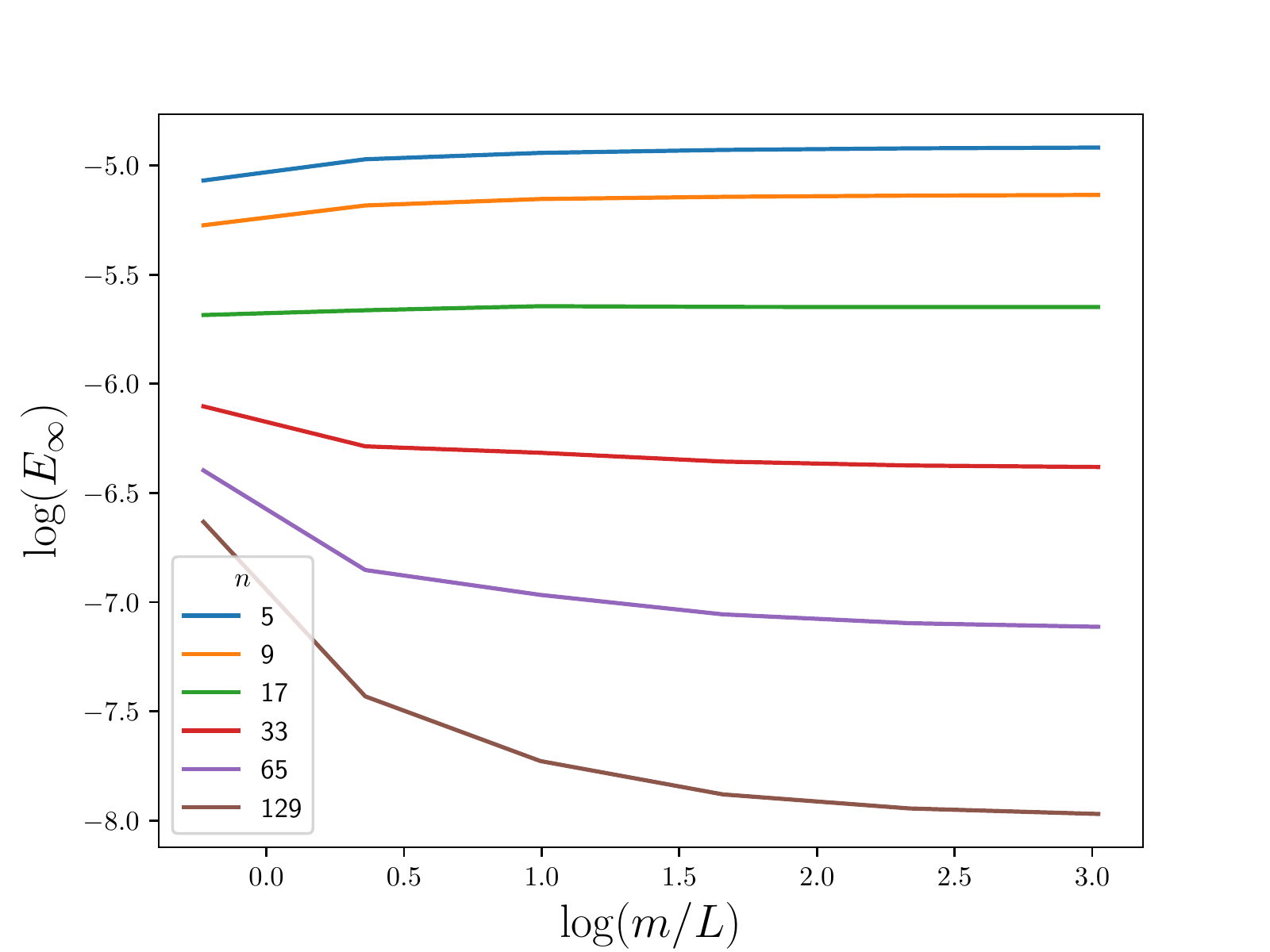}
\includegraphics[width = 0.48\textwidth]{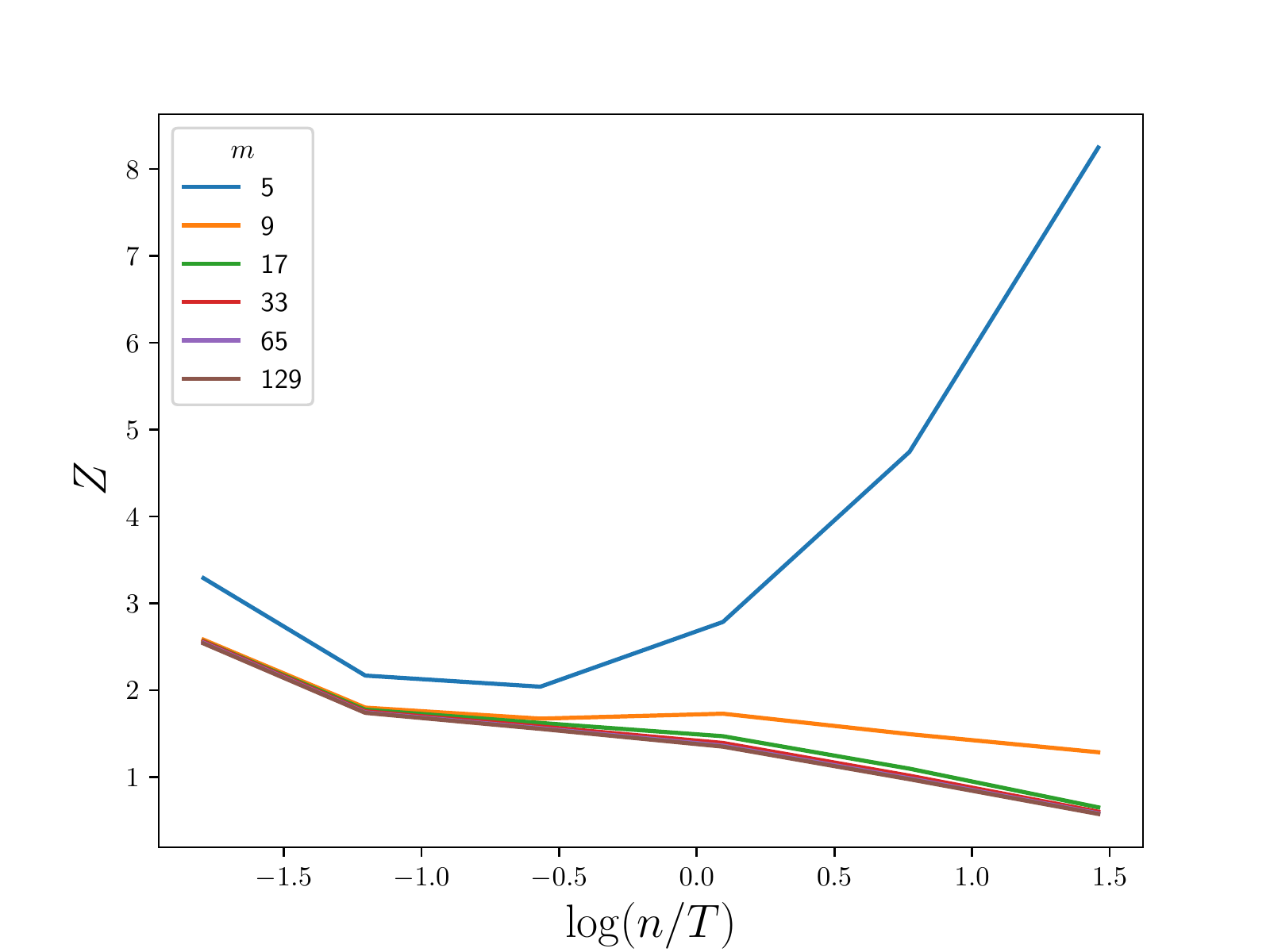}
\includegraphics[width = 0.48\textwidth]{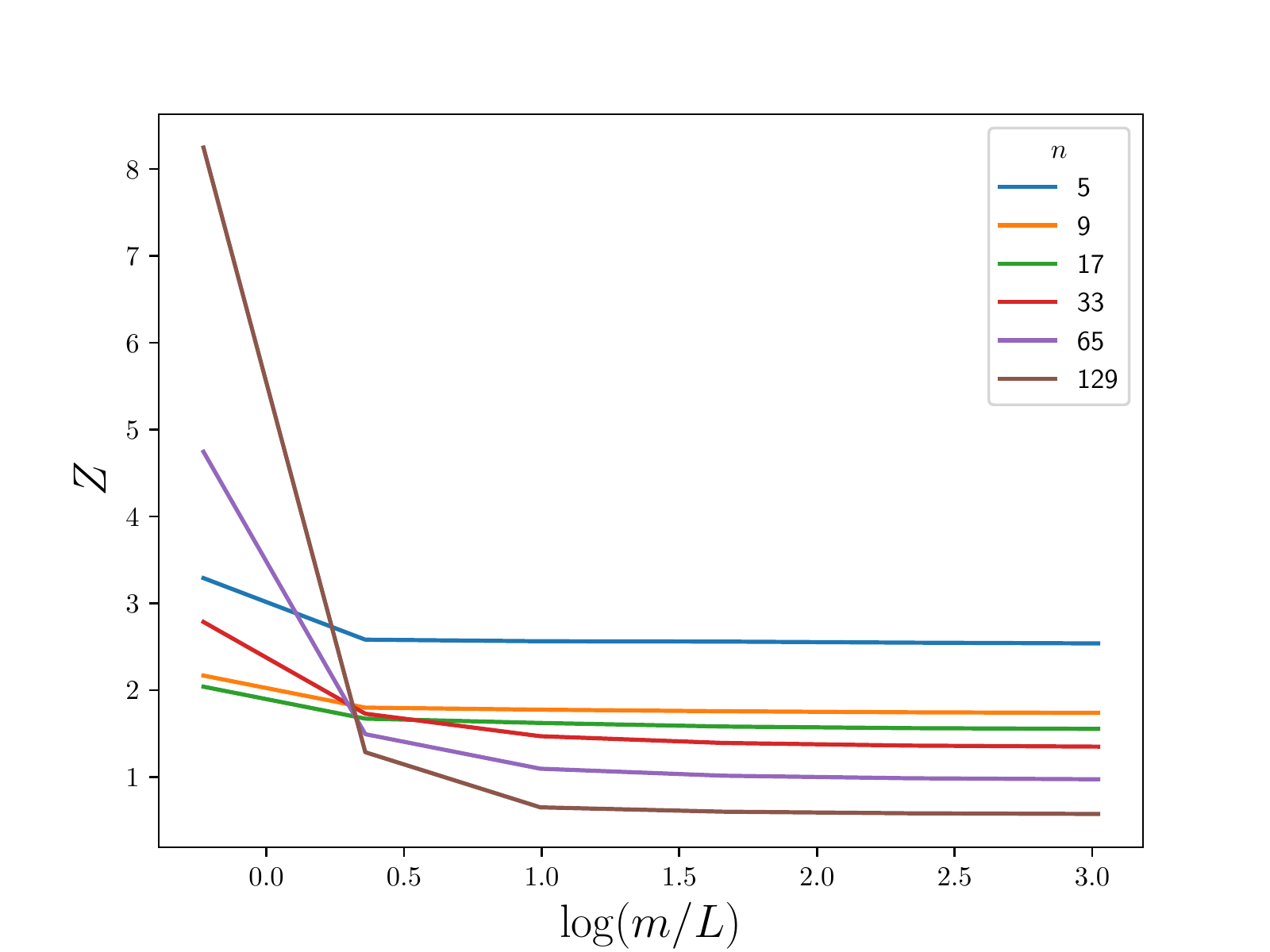}
\caption{Homogeneous Burger's equation, default prior: 
For each pair $(n,m)$ of temporal ($n$) and spatial ($m$) grid sizes considered, we plot:
(top left) the error $E_\infty$ for fixed $m$ and varying $n$;
(top right) the error $E_\infty$ for fixed $n$ and varying $m$;
(bottom left) the $Z$-score for fixed $m$ and varying $n$;
(bottom right) the $Z$-score for fixed $n$ and varying $m$.
}
\label{fig: burgernoforceerrorplots}
\end{figure}

\begin{figure}[t!]
\centering
\includegraphics[width = 0.48\textwidth]{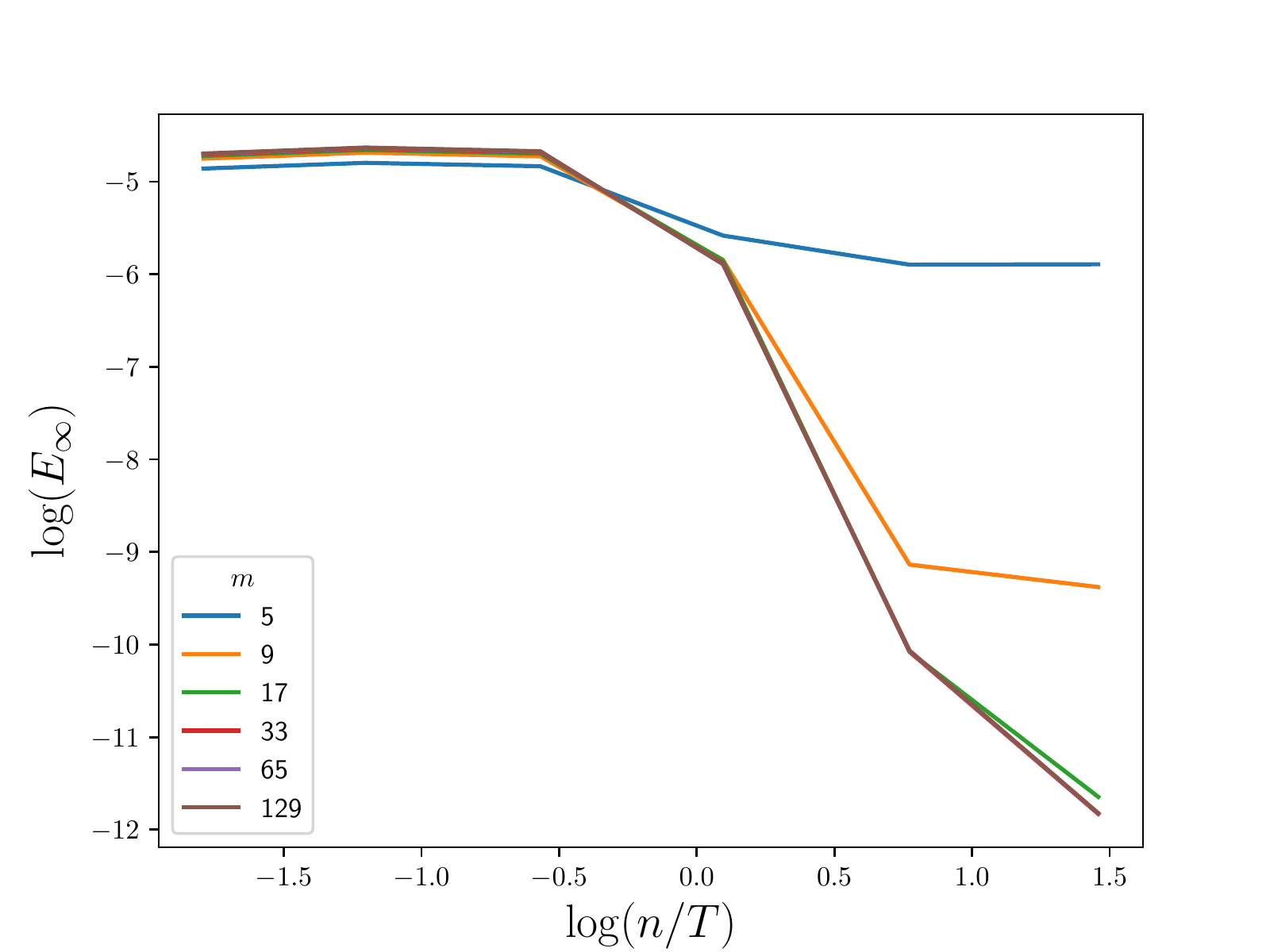}
\includegraphics[width = 0.48\textwidth]{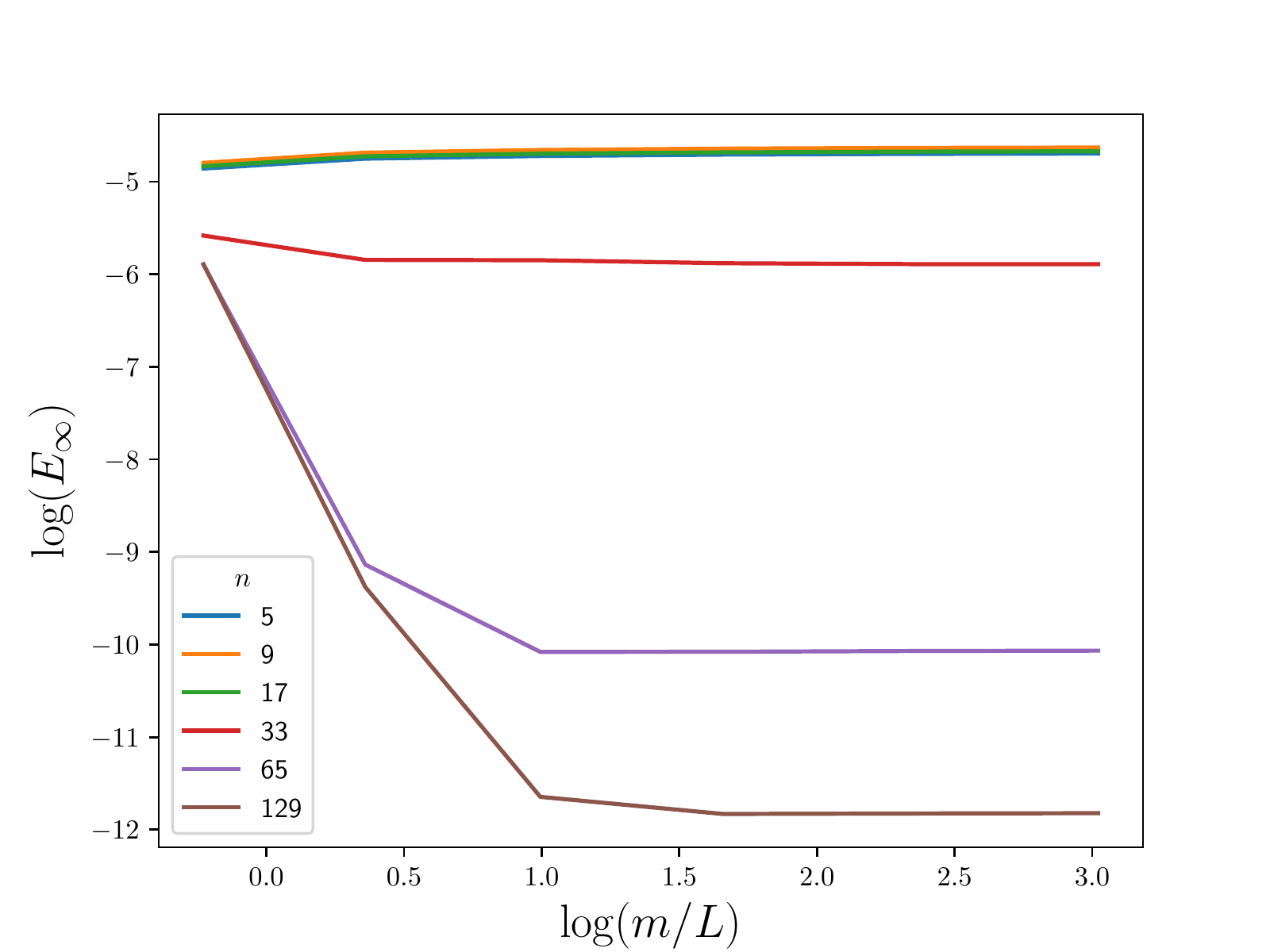}
\includegraphics[width = 0.48\textwidth]{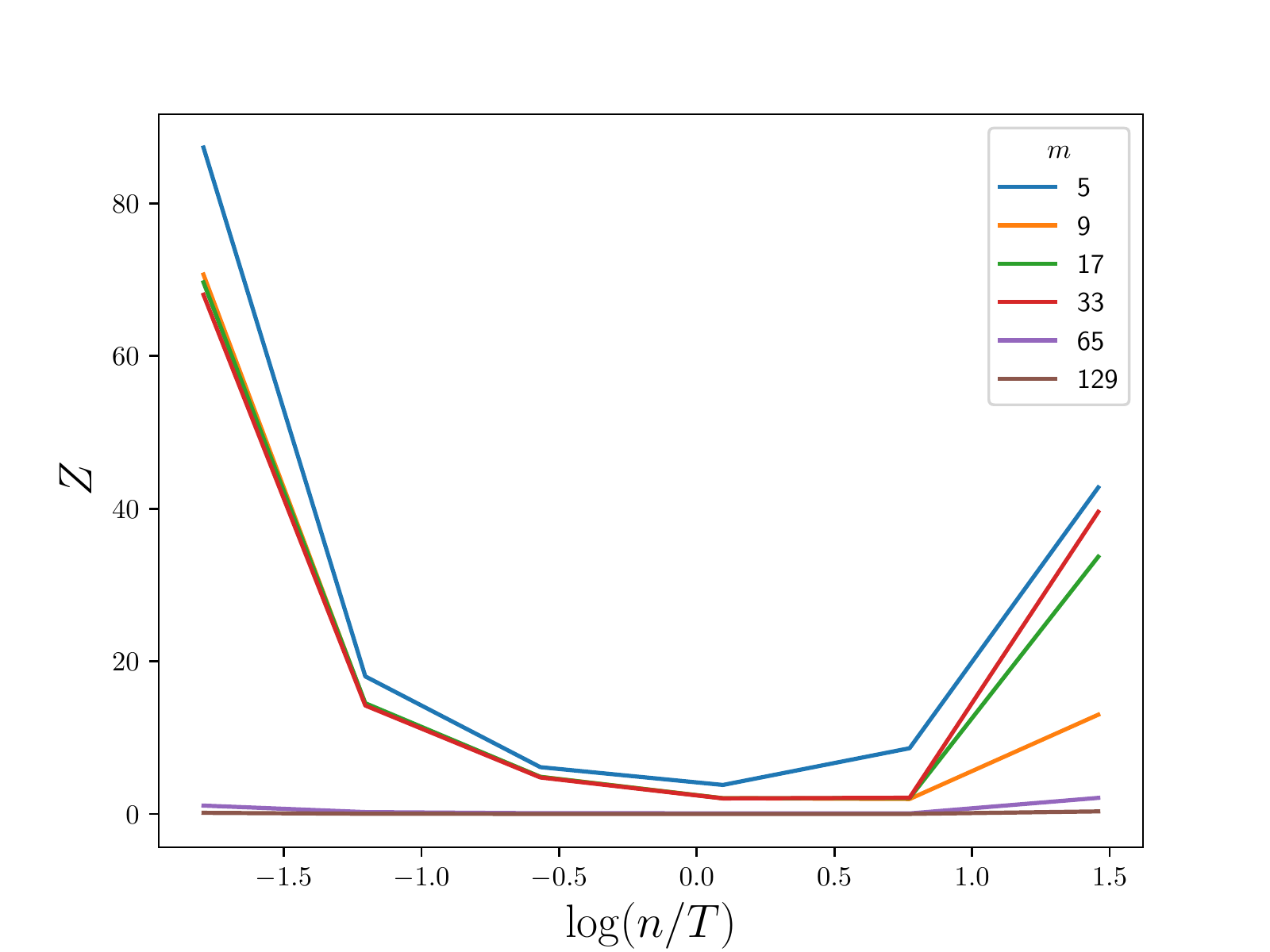}
\includegraphics[width = 0.48\textwidth]{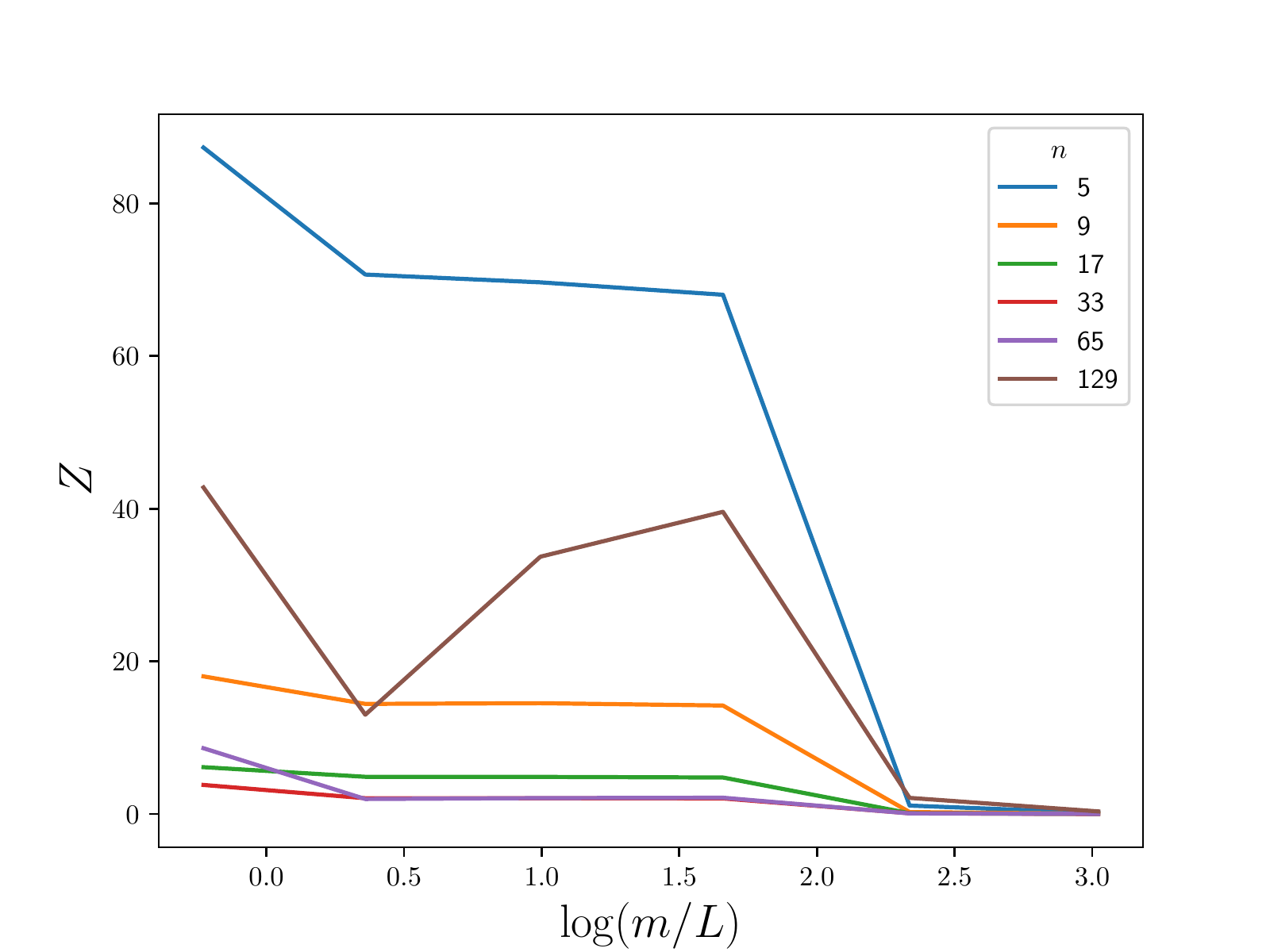}
\caption{Homogeneous Burger's equation, alternative prior: 
For each pair $(n,m)$ of temporal ($n$) and spatial ($m$) grid sizes considered, we plot:
(top left) the error $E_\infty$ for fixed $m$ and varying $n$;
(top right) the error $E_\infty$ for fixed $n$ and varying $m$;
(bottom left) the $Z$-score for fixed $m$ and varying $n$;
(bottom right) the $Z$-score for fixed $n$ and varying $m$.
}
\label{fig: burgernoforcerationalerrorplots}
\end{figure}

\paragraph{Results:}
The error $E_\infty$ was computed at 36 combinations of temporal and spatial grid sizes $(n,m)$ and results for the default prior are displayed in the top row of \Cref{fig: burgernoforceerrorplots}.
It can be seen that the error $E_\infty$ is mostly determined, in this example, by the finite length $n$ of the temporal grid rather than the length $m$ of the spatial grid.
The slope of the curves in \Cref{fig: burgernoforceerrorplots} is consistent with a convergence rate of $O(n^{-1})$ for the error $E_\infty$ when spatial discretisation is neglected. 
The $Z$-scores asociated with the default prior (bottom row of \Cref{fig: burgernoforceerrorplots}) appear to be bounded as $(n,m)$ are increased, tending toward 0 but taking values of order $1$ for all regimes, except for the smallest value ($m=5$) of the spatial grid. 
This provides evidence that our proposed PNM, equipped with the default prior, is either calibrated or slightly under-confident but, crucially from a statistical perspective, it is not over-confident.

Equivalent results for the alternative covariance model are presented for the error $E_\infty$ in the top row of \Cref{fig: burgernoforcerationalerrorplots} and for the $Z$-score in the bottom row of \Cref{fig: burgernoforcerationalerrorplots}.
Here the error $E_\infty$ is again gated by the size $n$ of the temporal grid and decreases at a faster rate compared to when the default prior was used. 
This is perhaps expected because the rational quadratic covariance model reflects the true smoothness of the solution better than the Mat\'{e}rn model. 
However, the $Z$-scores associated with the alternative covariance model are considerably higher, appearing to grow rapidly as $n \to \infty$ with $m$ fixed.
This suggests that the alternative covariance model is inappropriate, causing our PNM to be over-confident.
We speculate this may be because $u$ does not belong to the support of the rational quadratic covariance model, but note that the support of a Gaussian process can be difficult to characterise \citep{karvonen2021small}. 
These results support our proposed strategy for prior selection in \Cref{sec: prior}.

\subsection{Porous Medium Equation}
\label{subsec: porous}

Our second example is the \textit{porous medium equation}
\begin{equation}
\frac{\partial u}{\partial t} - \frac{\partial^2 (u^k)}{\partial x^2}=0 , \qquad t \in [t_0,t_0+T], \; x \in [-L/2,L/2], \label{eq: porousmediumunexpanded}
\end{equation}
which is more challenging compared to Burger's equation because the solution is only piecewise smooth, meaning that a strong solution does not exist and our modelling assumptions are violated.
Furthermore, there are two distinct nonlinearities in the differential operator, allowing us to explore the impact of the choice of linearisation on the performance of the PNM.
For our experiment we fix $k=2$, so that the porous medium equation becomes
\begin{equation*}
\frac{\partial u}{\partial t} - 2\left(\frac{\partial u}{\partial x}\right)^2 -2u\frac{\partial^2 u}{\partial x^2}=0 , \qquad t \in [t_0,t_0 + T], \; x \in [-L/2,L/2], \label{eq: porousmedium}
\end{equation*}
and we consider the initial and boundary conditions
\begin{equation*}
\begin{alignedat}{2}
u(t_0,x) & = t_0^{-1/3} \max(0,1-x^2/(12 t_0^{2/3})), && \qquad x \in [-L/2,L/2] \\
u(t,-L/2) = u(t,L/2) & = 0, && \qquad t \in [t_0,t_0+T]
\end{alignedat} \label{eq: porousconditions}
\end{equation*}
with $t_0 = 2$, $T = 8$, $L/2 = 10$.
These initial and boundary conditions were chosen because they permit a (unique) closed-form solution, due to \cite{Barenblatt1952}:
\begin{align*}
u(t,x)=\max \left( 0, \frac{1}{t^{1/3}} \left( 1-\frac{1}{12}\frac{x^2}{t^{2/3}} \right) \right)
\end{align*}
The solution is therefore only piecewise smooth, with discontinuous first derivatives at $x^2 = 12 t^{2/3}$, which are inside of the domain $[-L/2,L/2]$ for all $t \in [t_0,t_0+T]$.

\paragraph{Prior:}
Henceforth we consider the default prior advocated in \Cref{sec: prior} and \Cref{subsec: homog burger}, with amplitude $\sigma$ estimated using maximum likelihood and length-scale parameters fixed at values $\rho_1=1$, $\rho_2=2$ (not optimised; based on a simple \emph{post-hoc} visual check).

\paragraph{Choice of Linearisation:}
The differential operator here contains the nonlinear component $Qu = (\partial_xu)^2 + u \partial_x^2 u$ that must be linearised.
The first term $(\partial_xu)^2$ appeared also in Burger's equation, and we linearise this term in an identical way to that used in \Cref{subsec: homog burger}.
The second term $u \partial_x^2u$ can be linearised in at least two distinct ways, fixing either $u$ or $\partial_x^2 u$ to suitable constant values adaptively based on quantities that have been pre-computed.
Thus we consider the two linearisations
\begin{align*}
Q_i^{(1)} u(t_i,x) & \coloneqq \mu^{i-1}(t_i,x) \frac{\partial u}{\partial x} (t_i,x) + \mu^{i-1}(t_i,x) \frac{\partial^2 u}{\partial x^2}(t_i,x) \\
Q_i^{(2)} u(t_i,x) & \coloneqq \mu^{i-1}(t_i,x) \frac{\partial u}{\partial x} (t_i,x) + u(t_i,x) \frac{\partial^2 \mu^{i-1}}{\partial x^2}(t_i,x) 
\end{align*}
where we recall that $\mu^{i-1}(t_i,x)$ is the predictive mean arising from the Gaussian process approximation $U^{i-1}$.
Through simulation we aim to discover which (if either) linearisation is more appropriate for use in our PNM.

\paragraph{Conservation of Mass:}
In addition to admitting multiple linearisations, we consider the porous medium equation because when $k > 1$ it exhibits a \textit{conservation law}, which is typical of many nonlinear PDEs that are physically-motivated.
Specifically, integrating \eqref{eq: porousmediumunexpanded} with respect to $x$ gives
\[
\frac{\mathrm{d}}{\mathrm{d} t}\int_{-L}^{L} u(t,x) \,\mathrm{d}x - \frac{\partial (u^k)}{\partial x}\Big|_{-L}^{L}=0
\]
and, from the fact that $u = 0$ for all $x^2 \geq 12 t^{2/3}$, it follows that $\partial_x (u^k) = ku^{k-1} \partial_x u = 0$ for all $x^2 \geq 12 t^{2/3}$ and thus $\int_{-L}^L u(t,x) \mathrm{d}x$ is $t$-invariant.
A desirable property of a numerial method is that it respects conservations law of this kind; as exemplified by the \textit{finite volume} methods \citep{leveque2002finite} and \textit{symplectic integrators} \citep{sanz1992symplectic}.
Interestingly, it is quite straight-forward to enforce this conservation law in our PNM by adding additional linear constraints to the system in \Cref{lem: batch}.
Namely, we add the linear constaints
\begin{align*}
\int_{-L}^L u(t_i,x) \mathrm{d}x & = \int_{-L}^L u(t_0,x) \mathrm{d}x = 4(3^{\frac{1}{2}}-3^{-\frac{1}{2}})
\end{align*}
at each point $i \in \{1,\dots,n-1\}$ on the temporal grid.
The performance of our PNM both with and without conservation of mass will be considered.

\begin{figure}[t!]
\centering
\includegraphics[width = 0.48\textwidth]{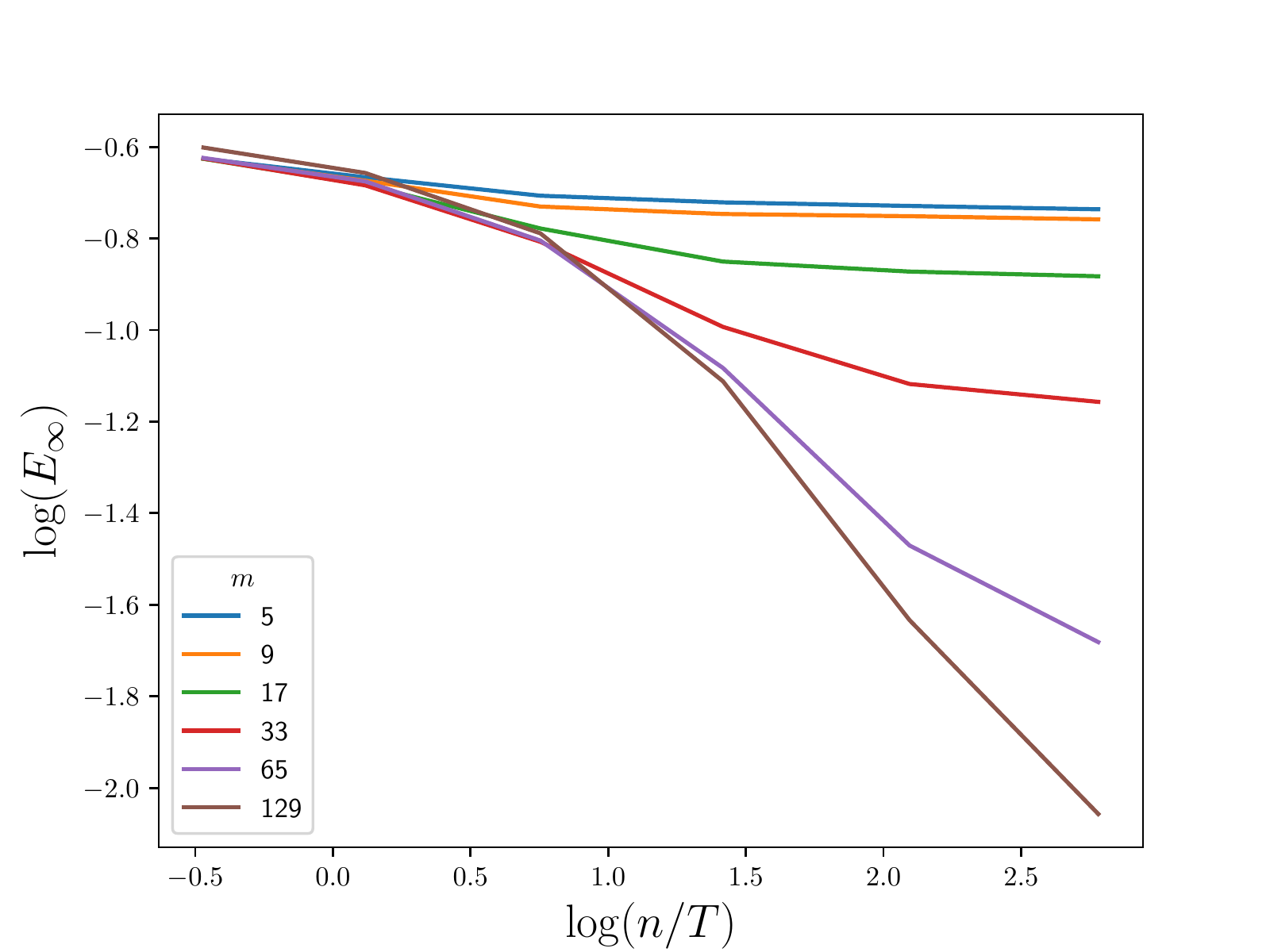}
\includegraphics[width = 0.48\textwidth]{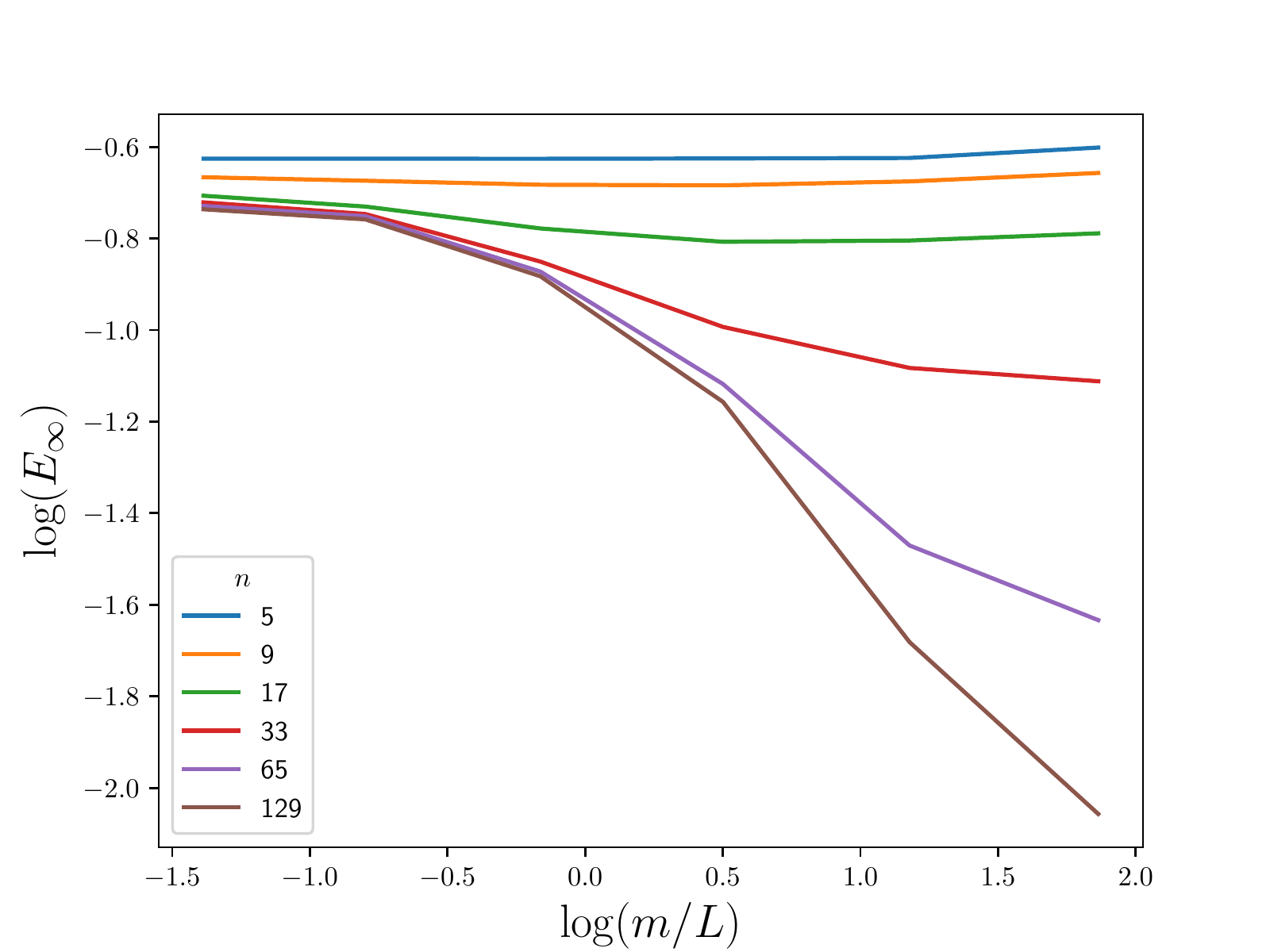}
\includegraphics[width = 0.48\textwidth]{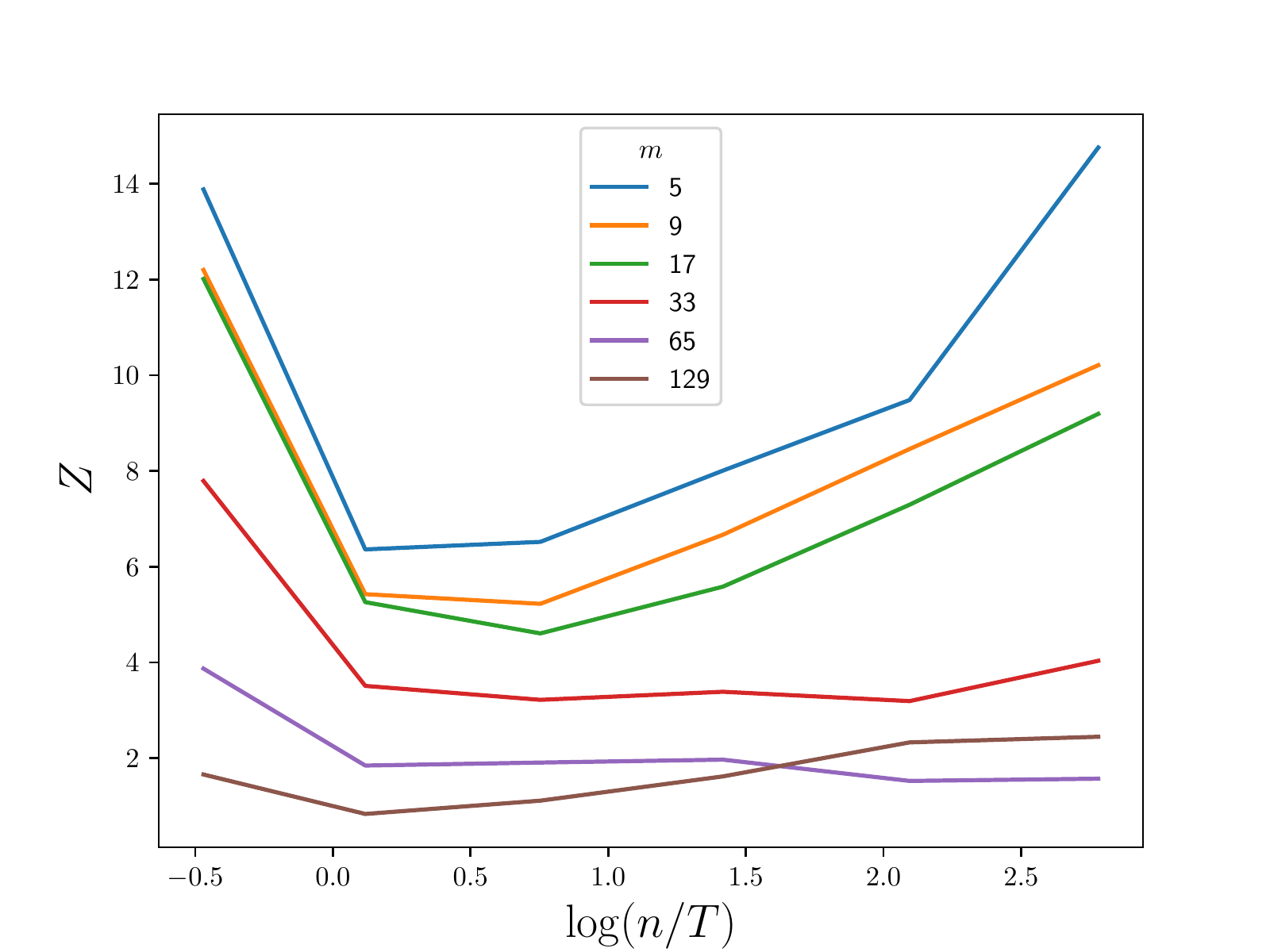}
\includegraphics[width = 0.48\textwidth]{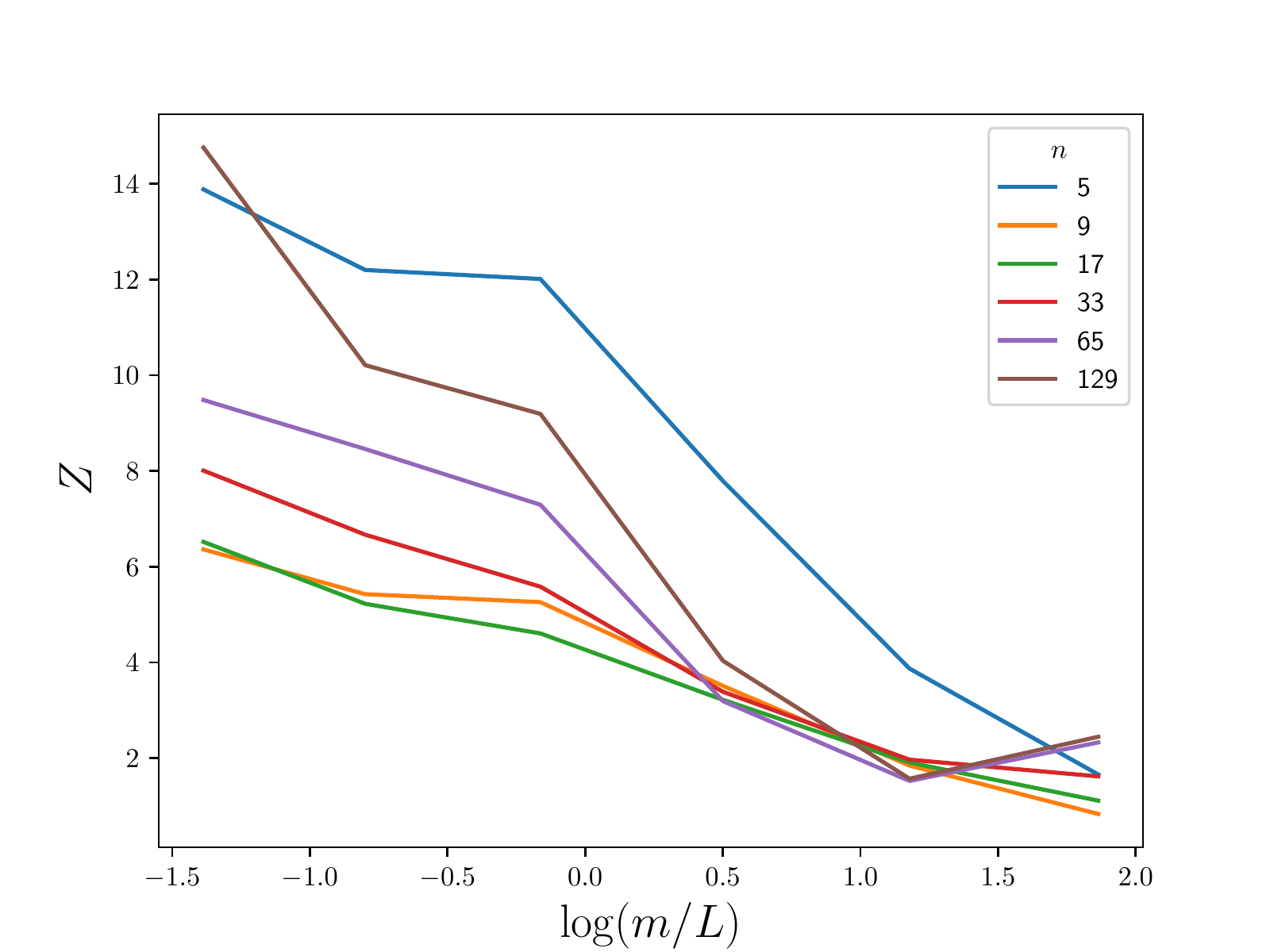}
\caption{
Porous medium equation, with linearisation $Q^{(1)}$:
For each pair $(n,m)$ of temporal ($n$) and spatial ($m$) grid sizes considered, we plot:
(top left) the error $E_\infty$ for fixed $m$ and varying $n$;
(top right) the error $E_\infty$ for fixed $n$ and varying $m$;
(bottom left) the $Z$-score for fixed $m$ and varying $n$;
(bottom right) the $Z$-score for fixed $n$ and varying $m$.
}
\label{fig: porouserrorplots}
\end{figure}

\begin{figure}[t!]
\centering
\includegraphics[width = 0.48\textwidth]{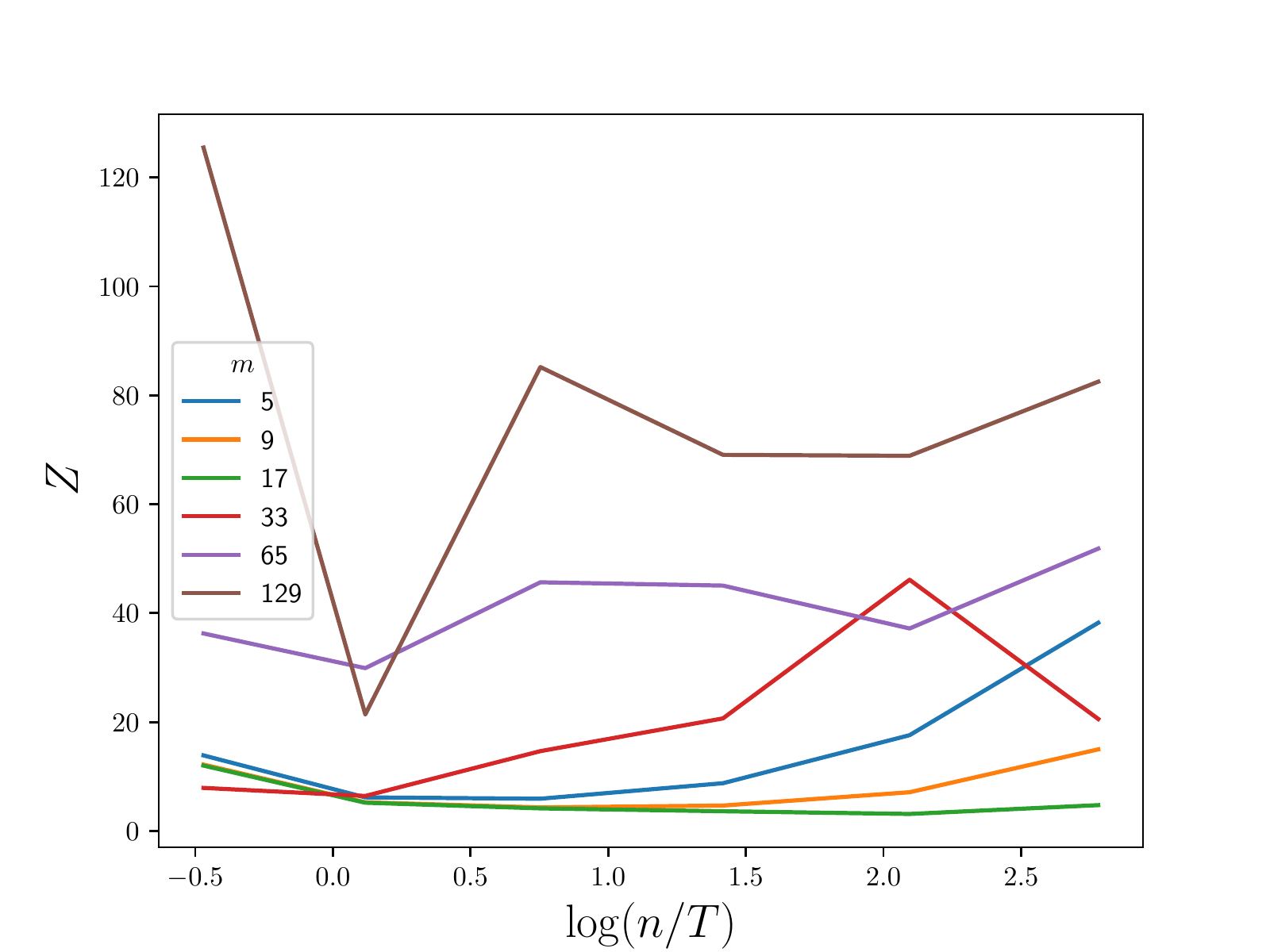}
\includegraphics[width = 0.48\textwidth]{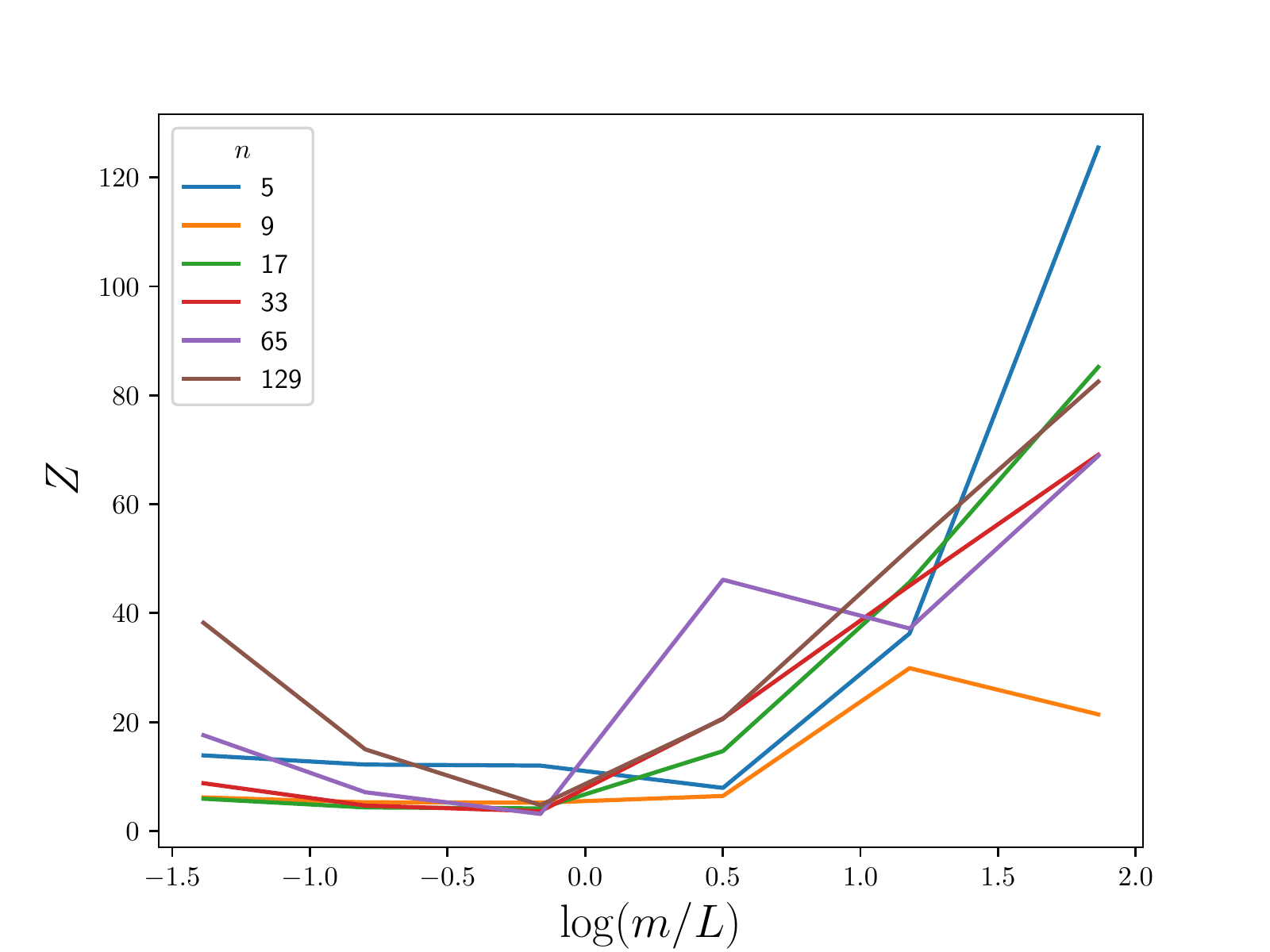}
\caption{
Porous medium equation, with linearisation $Q^{(2)}$:
For each pair $(n,m)$ of temporal ($n$) and spatial ($m$) grid sizes considered, we plot:
(left) the $Z$-score for fixed $m$ and varying $n$;
(right) the $Z$-score for fixed $n$ and varying $m$.
}
\label{fig: porousalterrorplots}
\end{figure}

\begin{figure}[t!]
\centering
\includegraphics[width = 0.48\textwidth]{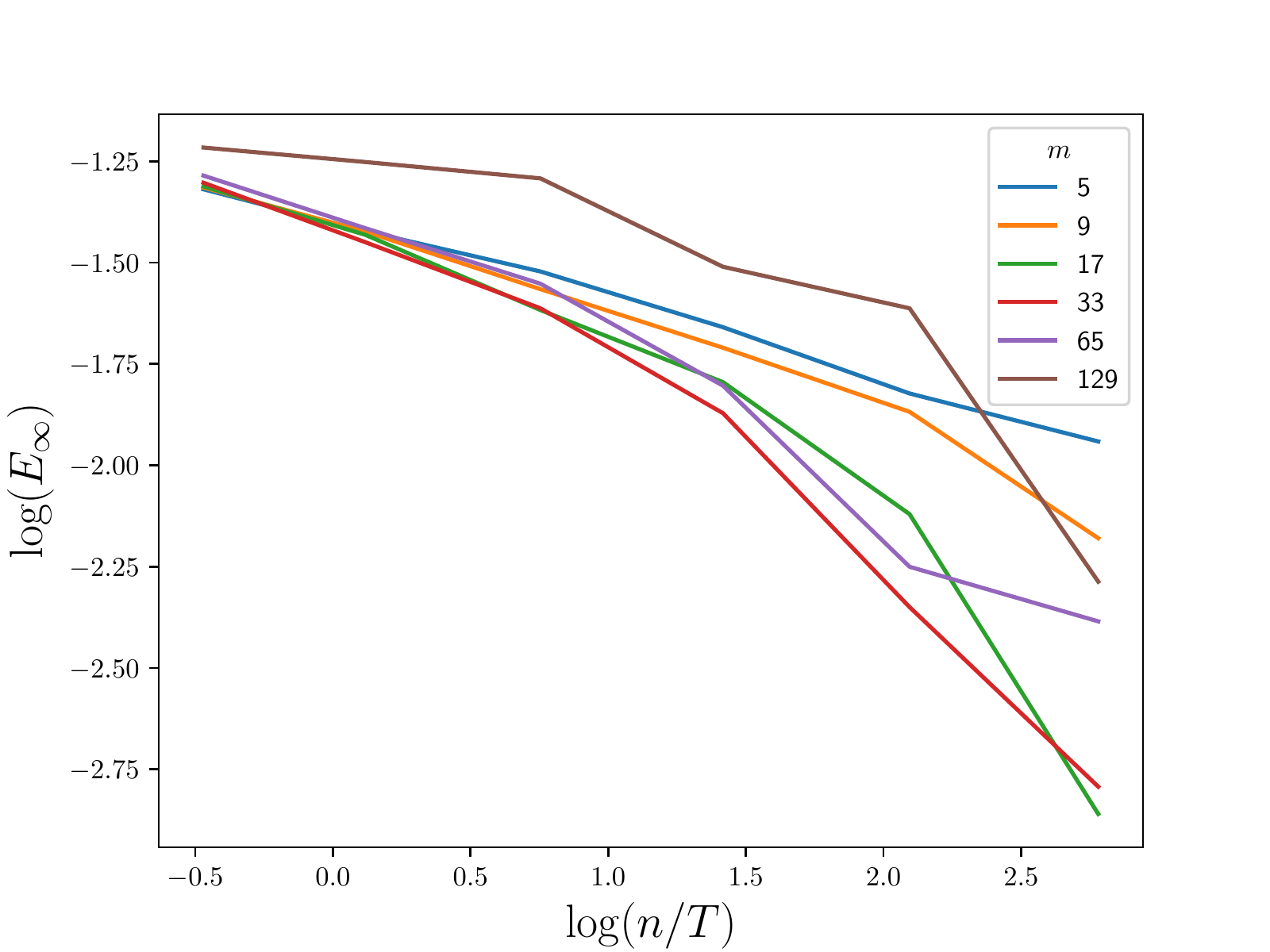}
\includegraphics[width = 0.48\textwidth]{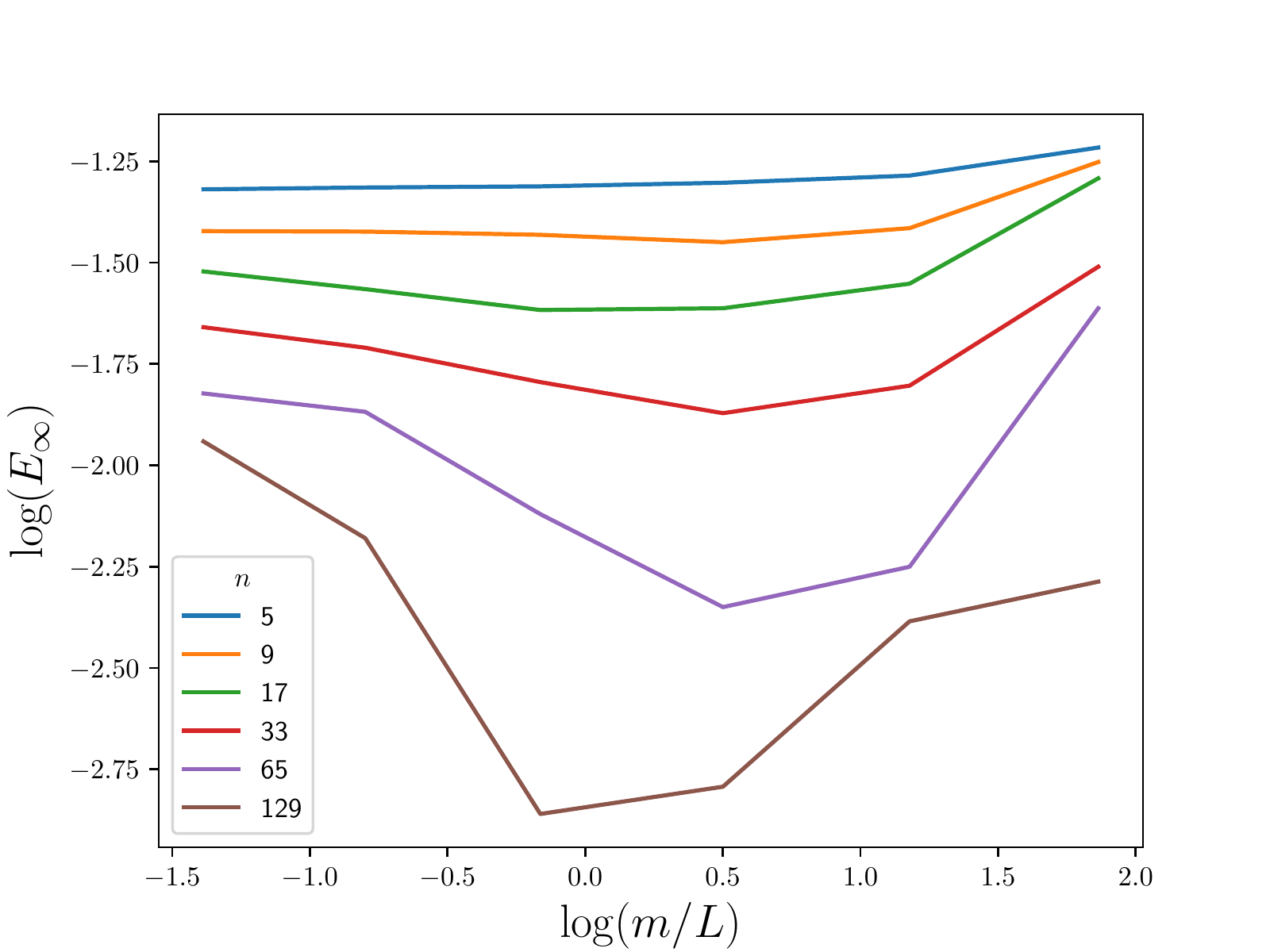}
\includegraphics[width = 0.48\textwidth]{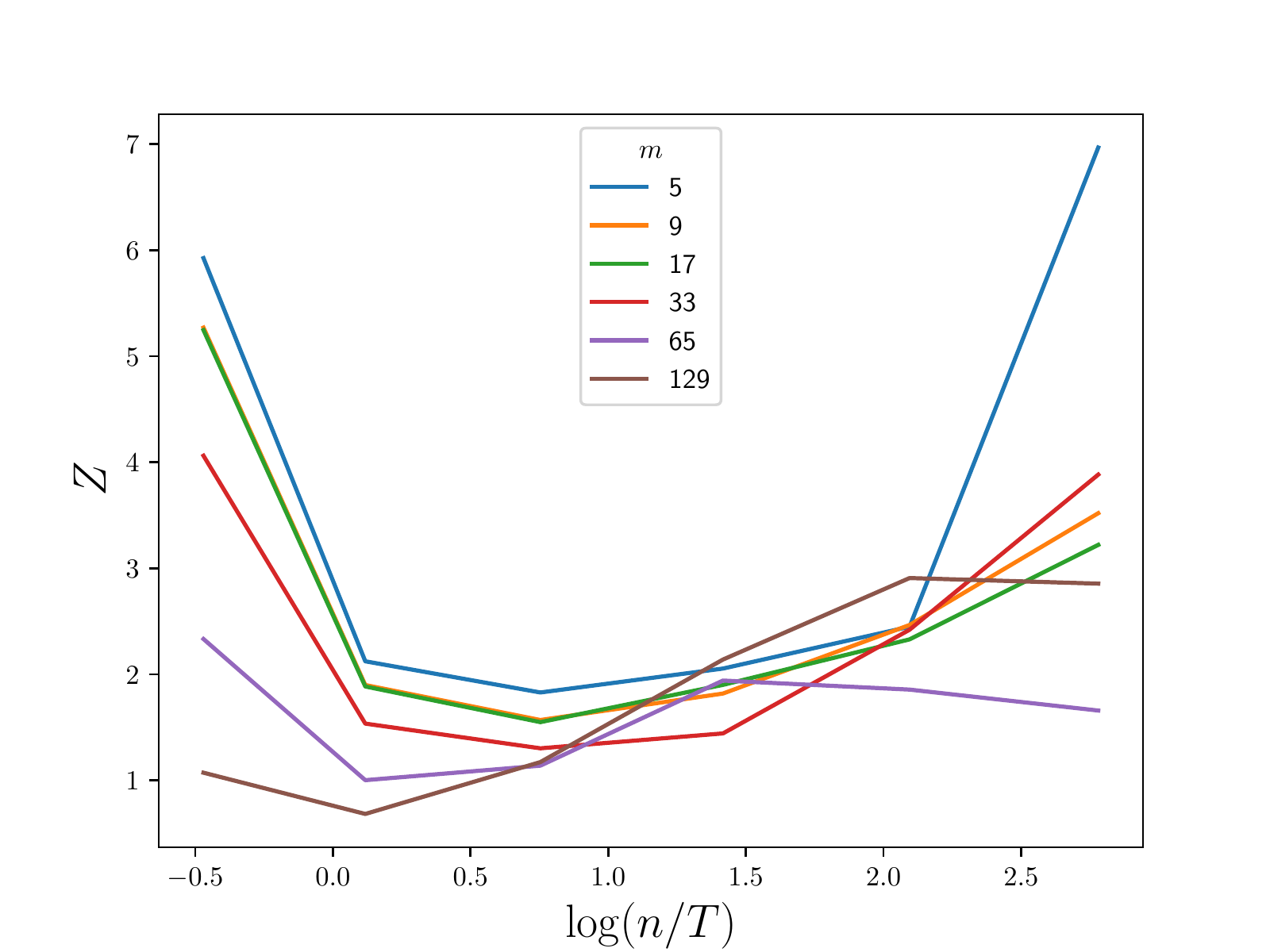}
\includegraphics[width = 0.48\textwidth]{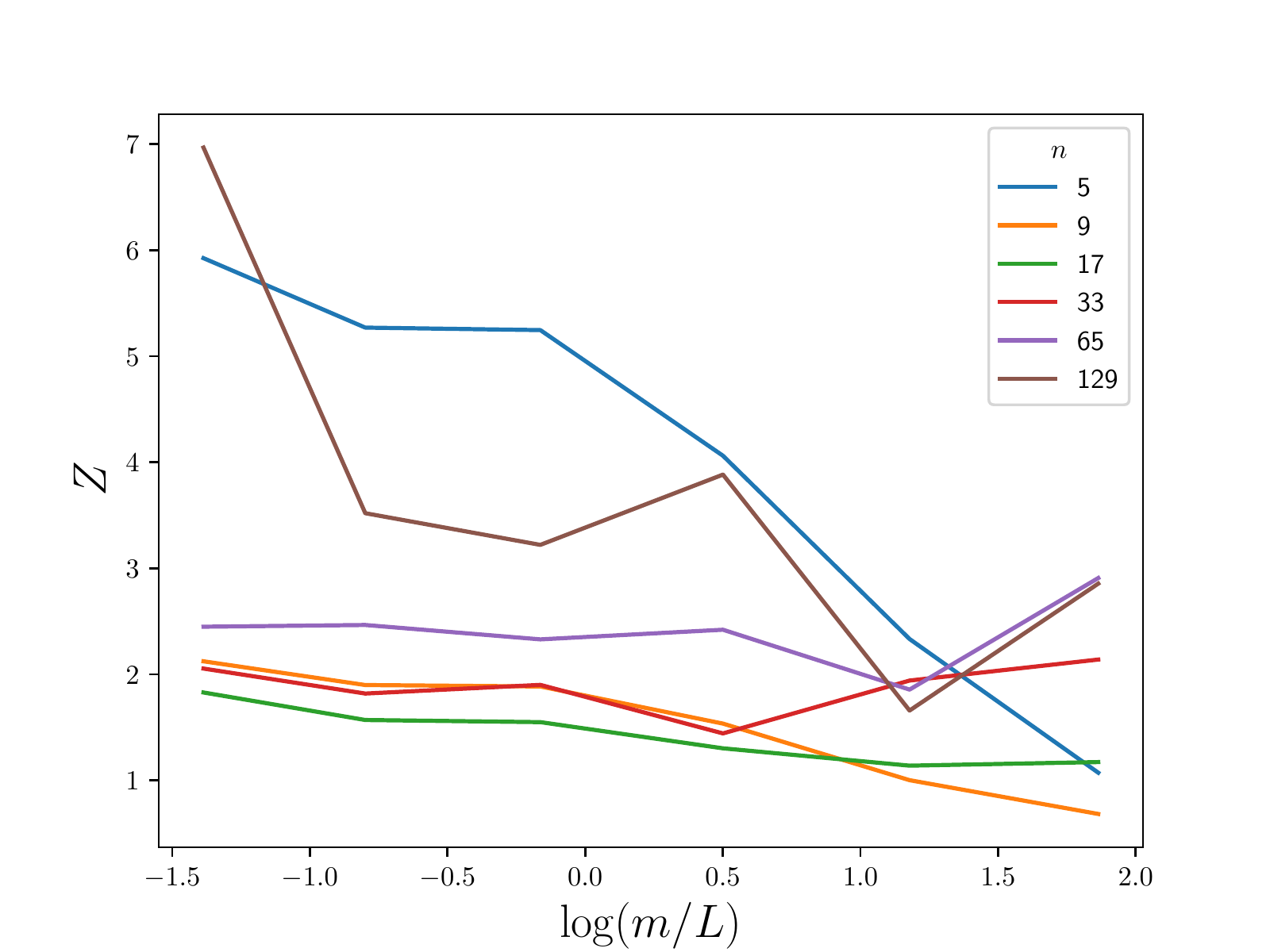}
\caption{
Porous medium equation, with mass conserved:
For each pair $(n,m)$ of temporal ($n$) and spatial ($m$) grid sizes considered, we plot:
(top left) the error $E_\infty$ for fixed $m$ and varying $n$;
(top right) the error $E_\infty$ for fixed $n$ and varying $m$;
(bottom left) the $Z$-score for fixed $m$ and varying $n$;
(bottom right) the $Z$-score for fixed $n$ and varying $m$.
}
\label{fig: porousintegralerrorplots}
\end{figure}

\paragraph{Results:}
Empirical results based on the linearisation $Q^{(1)}$ (without conservation of mass) are contained in \Cref{fig: porouserrorplots}.
In this case (and in contrast to our results for Burger's equation in \Cref{subsec: homog burger}), the error $E_\infty$ is seen to be gated by the smaller of the finite length $n$ of the temporal grid and the length $m$ of the spatial grid.
The $Z$-score values appear to be of order 1 as $(n,m)$ are simultaneously increased, but are slightly higher than for Burger's equation, which may reflect the fact that the solution to the porous medium equation is only piecewise smooth.
For increasing $n$ with $m$ fixed the PNM appears to become over-confident, while for increasing $m$ with $n$ fixed the PNM appears to become under-confident; a conservative choice would therefore be to take $m \geq n$.
Interestingly, this behaviour of the $Z$-scores is similar to that observed for the rational quadratic covariance model in \Cref{subsec: homog burger}, and may reflect the fact that in both cases the solution $u$ is outside the support of the covariance model.

Next we compared the performance of the linearisation $Q^{(1)}$ with the linearisation $Q^{(2)}$.
The error $E_\infty$ associated to $Q^{(2)}$ (not shown) was larger than the error of $Q^{(1)}$, and the $Z$-scores for $Q^{(2)}$ are displayed in \Cref{fig: porousalterrorplots}.
Our objective is to quantify numerical uncertainty, so it is essential that output from the PNM is calibrated.
Unfortunately, it can be seen that the $Z$-scores associated with $Q^{(2)}$ are unsatisfactory; for large $m$ the scores are two orders of magintude larger than 1, indicating that the PNM is over-confident.
The failure of $Q^{(2)}$ to provide calibrated output is likely due to the fact that approximation of the second order derivative term $\partial_x^2 u$ is more challenging compared to approximation of the solution $u$, since $\partial_x^2$ is less regular than $u$ and since our initial and boundary data relate directly to $u$ itself.

Finally we considered inclusion of the conservation law into the PNM.
For this purpose we used the best-performing linearisation $Q^{(1)}$.
The errors $E_\infty$ and $Z$-scores are shown in \Cref{fig: porousintegralerrorplots} and can be compared to the equivalent results without the conservation law applied, in \Cref{fig: porouserrorplots}.
It can be seen that the error $E_\infty$ is lower when the conservation law is applied, and moreover the $Z$-scores are slightly reduced, remaining order 1.
These results agree with the intuition that incorporating additional physical constraints, when they are known, can have a positive impact on the performance of our PNM.

\subsection{Forced Burger's Equation}

Our final experiment concerns a nonlinear PDE whose right hand side $f$ is considered to be a black box, associated with a substantial computational cost.
To avoid confounding due to the choice of differential operator, we consider again the differential operator from Burger's equation
\begin{equation}
\frac{\partial u}{\partial t} + u\frac{\partial u}{\partial x} -\alpha\frac{\partial^2 u}{\partial x^2} = f(t,x), \qquad t \in [0,T], \; x \in [0,L], \label{eq: Burgerforced}
\end{equation}
for which the behaviour of our PNM was studied in \Cref{subsec: homog burger} (in the case $f = 0$).
The initial and boundary conditions are 
\begin{equation*}
\begin{alignedat}{2}
u(0,x) & = 0, && \qquad x \in [0,L] \\
u(t,0) = u(t,L) & = 0, && \qquad t \in [0,T]
\end{alignedat}
\end{equation*}
and we set $\alpha=1$, $T = 30$, $L=1$.
The aims of this experiment are two-fold: 
Our first aim is to evaluate the performance of our PNM when the function $f$ is non-trivial (e.g.\ involving oscillatory behaviour), to understand whether the output from our PNM remains calibrated or not.
Recall that our experiments are \textit{synthetic}, meaning that the black box $f$ is in actual fact an analytic expression, in this case 
\begin{equation*}
f(t,x) \coloneqq 10 \sin (6\pi x) \cos (3\pi t)+2 \left|\sin(3\pi x) \cos(6\pi t) \right|,
\end{equation*}
enabling a thorough assessment to be performed.
This forcing term is deliberately chosen to have some non-smoothness (from the absolute value function) and oscillatory behaviour, as might be encountered in output from a complex computer model.
Our second aim is to compare the accuracy of our PNM against a classical numerical method whose computational budget (as quantified by the number of times $f$ is evaluated) is identical to our PNM.

The solution to \eqref{eq: Burgerforced} does not admit a closed form, so for our ground truth we used a numerical solution computed using the \texttt{MATLAB} function \texttt{pdepe}, which implements \cite{skeel1990method} based on a uniform spatial grid of size $512$ and an adaptively selected temporal grid. 
Our PNM was implemented with the same linearisation used in \Cref{subsec: homog burger}.

\paragraph{Prior:}
Again we consider the default prior advocated in \Cref{sec: prior} and \Cref{subsec: homog burger}, with amplitude $\sigma$ estimated using maximum likelihood and length-scale parameters fixed at values $\rho_1=0.5$, $\rho_2=0.5$ (not optimised; based on a simple \emph{post-hoc} visual check).

\paragraph{Crank--Nicolson Benchmark:}
In this scenario, where the black box function $f$ is associated with a high computational cost, non-adaptive numerical methods are preferred, to control the total computational cost.
From a classical perspective, the finite difference methods are natural candidates for the numerical solution of \eqref{eq: Burgerforced}. Finite difference methods are classified into explicit and implicit schemes. Explicit schemes are much easier to solve but typically require certain conditions to be met for numerical stability \citep[][Table 5.3.1]{Thomas1998}. For example, for the 2D heat equation, it is required that $\Delta t/(\alpha\Delta x^2)+\Delta t/(\alpha\Delta y^2) \leq 1/2$, where $\alpha$ is the diffusivity constant, $\Delta t$ the time resolution, and $\Delta x, \Delta y$ the spatial resolutions  \citep[][page 158]{Thomas1998}. Such conditions, which require the spatial resolution to be much finer than the time resolution, may be difficult to establish when $f$ is a black box or when manual selection of the solution grid is not possible. Implicit methods in general are more difficult to solve, but stability is often guaranteed. For example, for the same 2D heat equation problem, the \textit{Crank--Nicolson} scheme (a second order, implicit method) is unconditionally stable \citep[][page 159]{Thomas1998}.  
For these reasons, we considered the \textit{Crank--Nicolson} finite difference method \citep{Crank1947} as a classical numerical method that is well-suited to the task at hand. In the implementation of Crank-Nicholson on the inhomogeneous Burger's equation, the nonlinear term is approximated via a \textit{lag nonlinear term}  \cite[][page 140]{Thomas1998}.
The same regular temporal grid $\mathbf{t}$ and regular spatial grid $\mathbf{x}$ were employed in both Crank--Nicolson and our PNM, so that the computational costs for both methods (as quantified in terms of the number of evaluations of $f$) are identical.

\begin{figure}[t!]
\centering
\includegraphics[width = 0.48\textwidth]{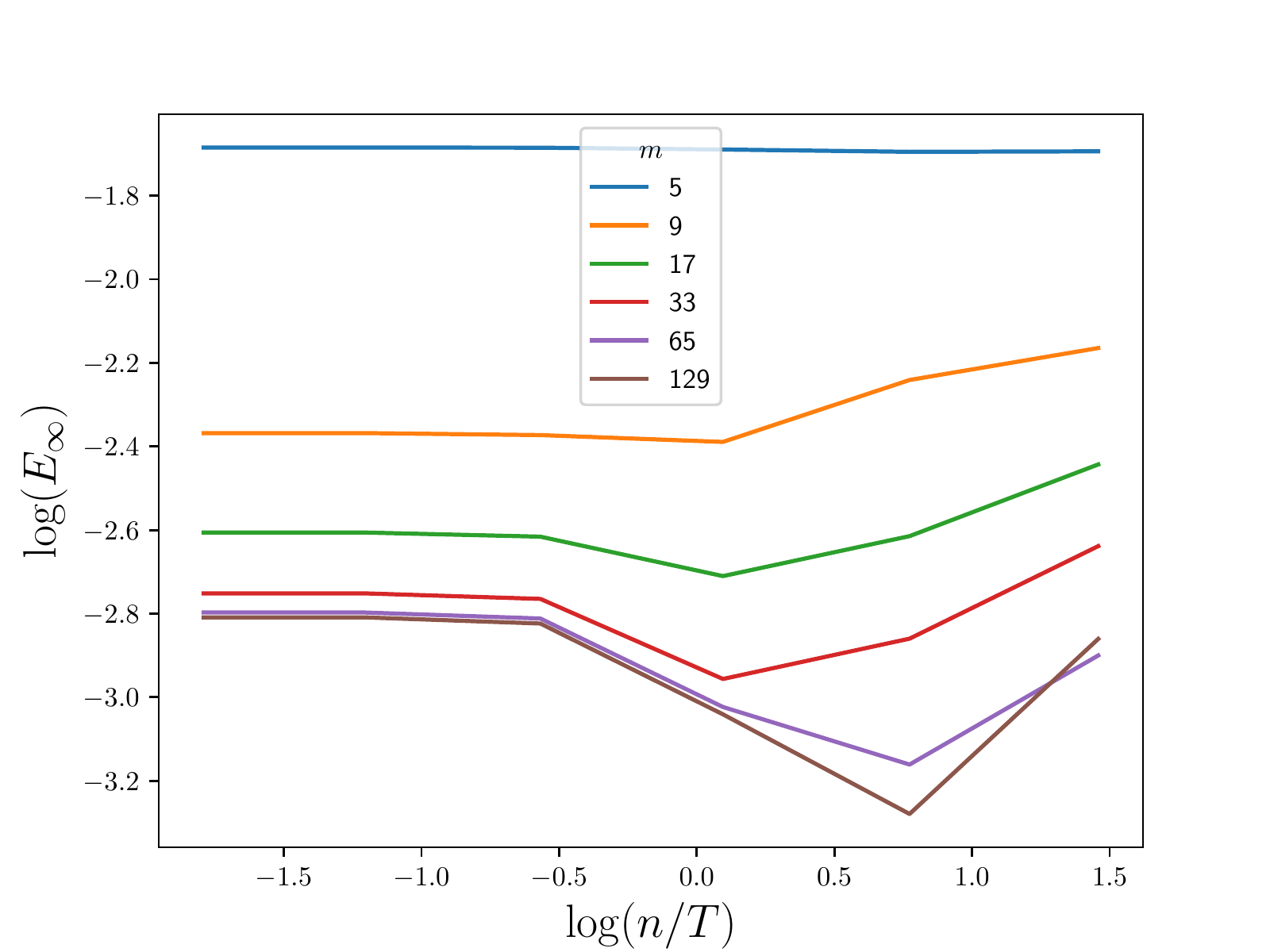}
\includegraphics[width = 0.48\textwidth]{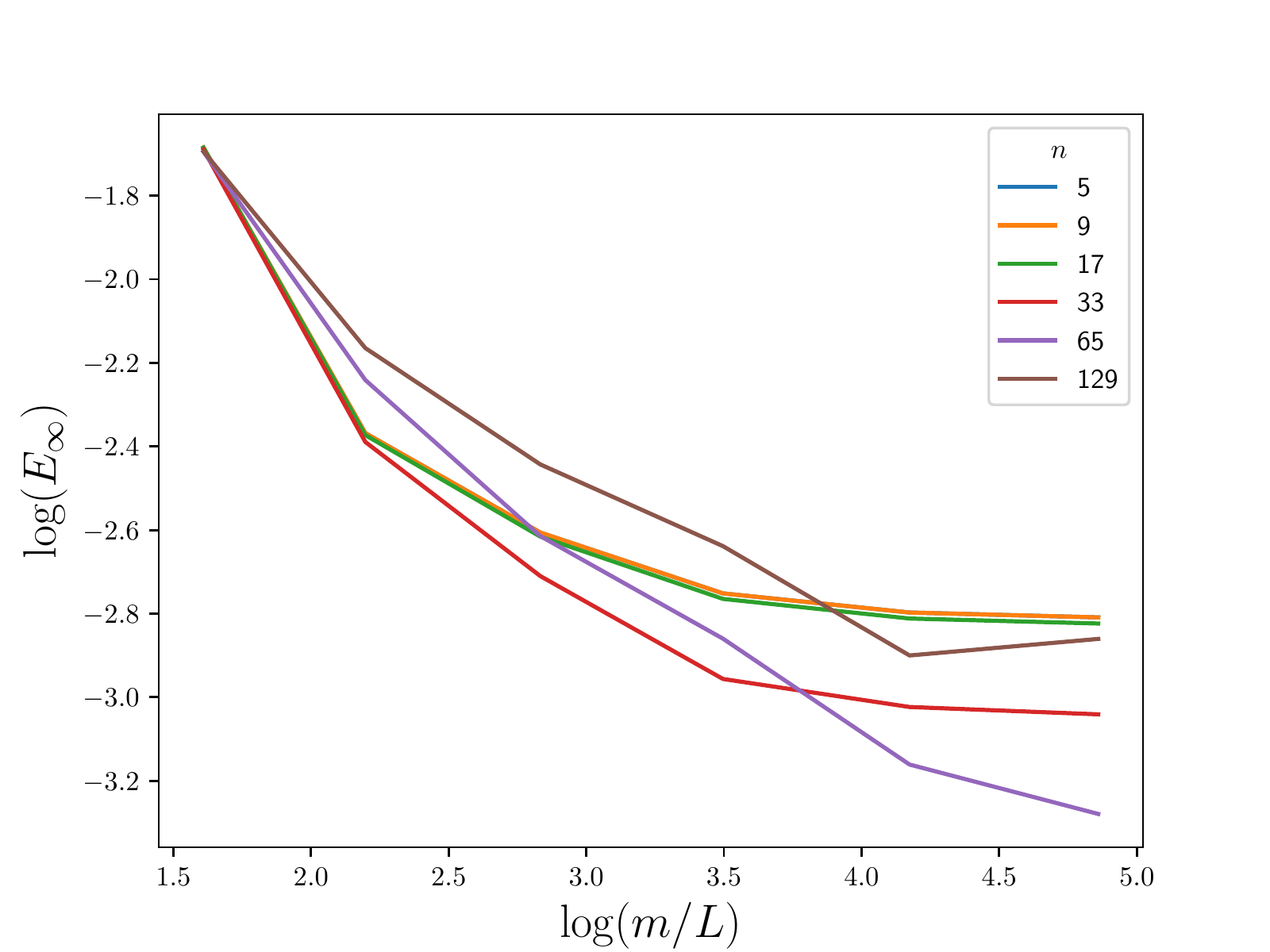}
\includegraphics[width = 0.48\textwidth]{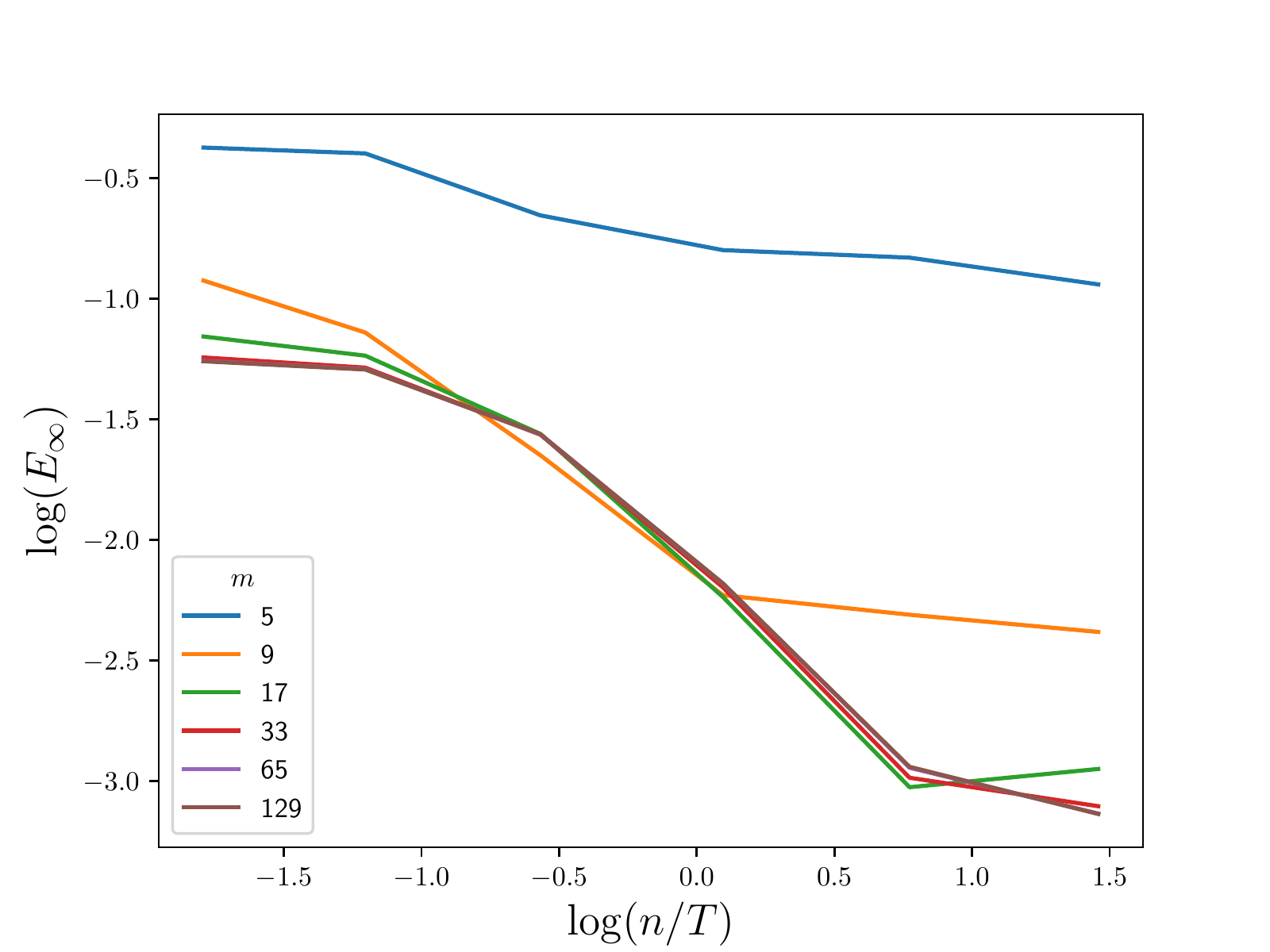}
\includegraphics[width = 0.48\textwidth]{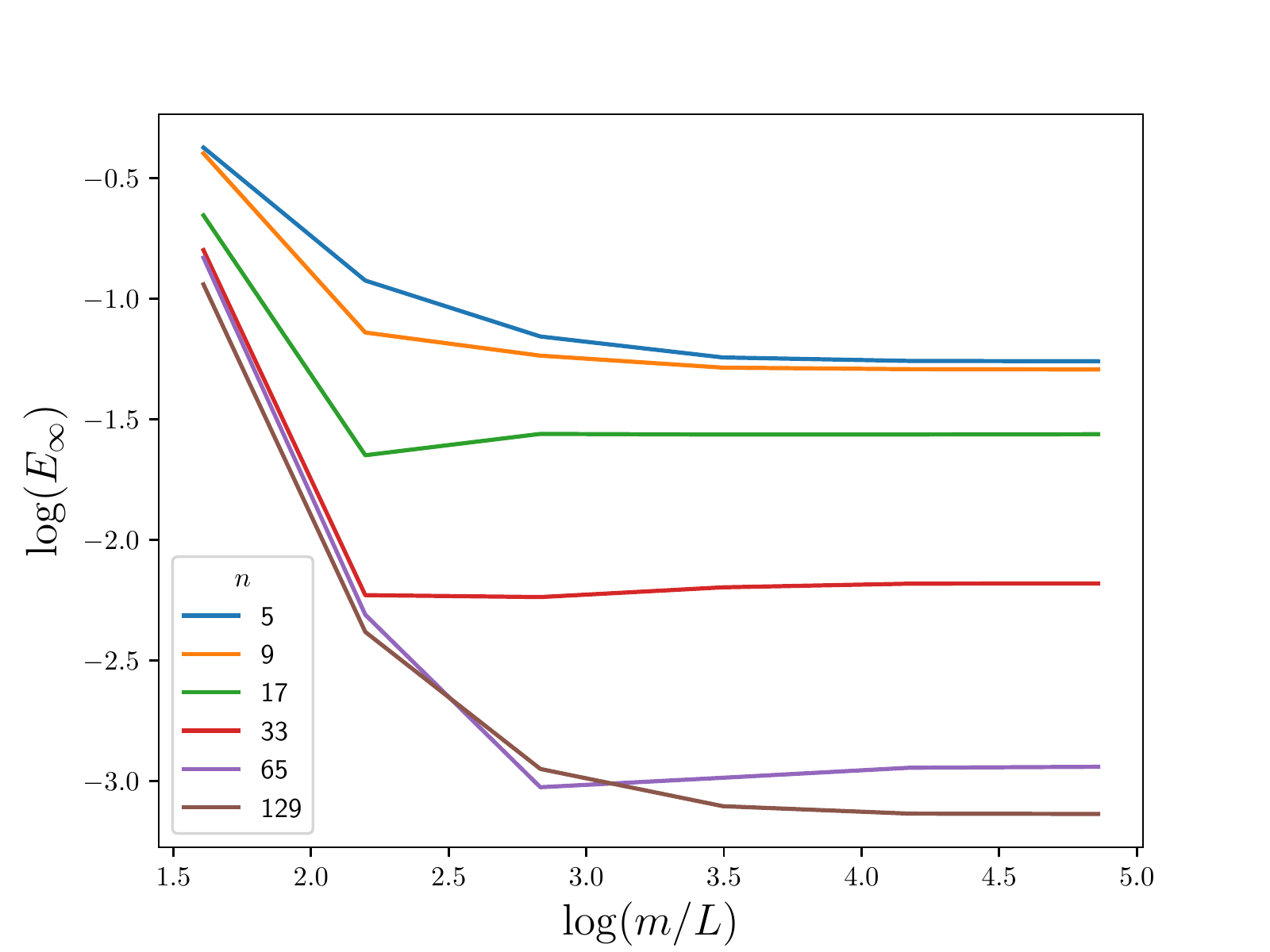}
\includegraphics[width = 0.48\textwidth]{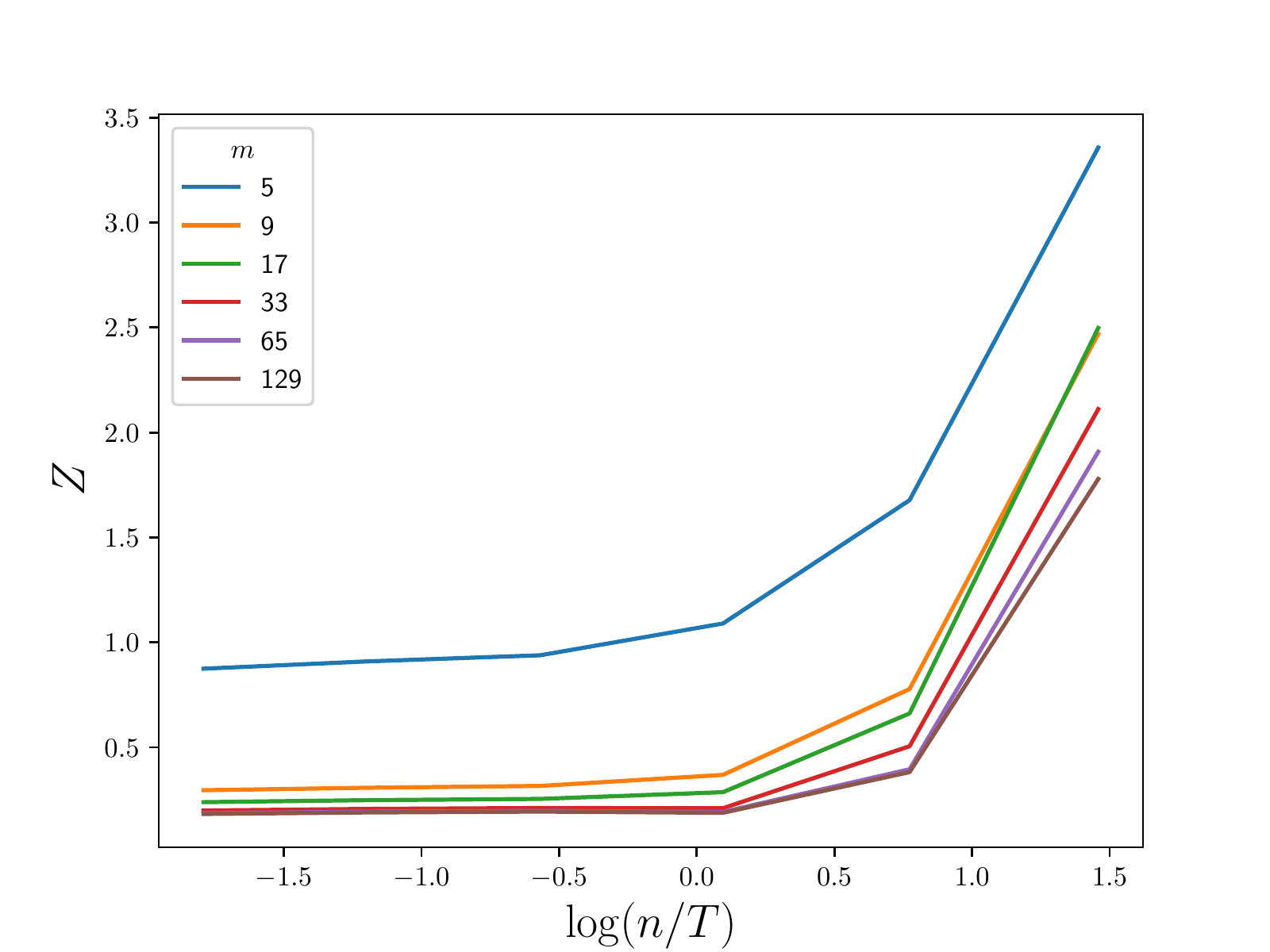}
\includegraphics[width = 0.48\textwidth]{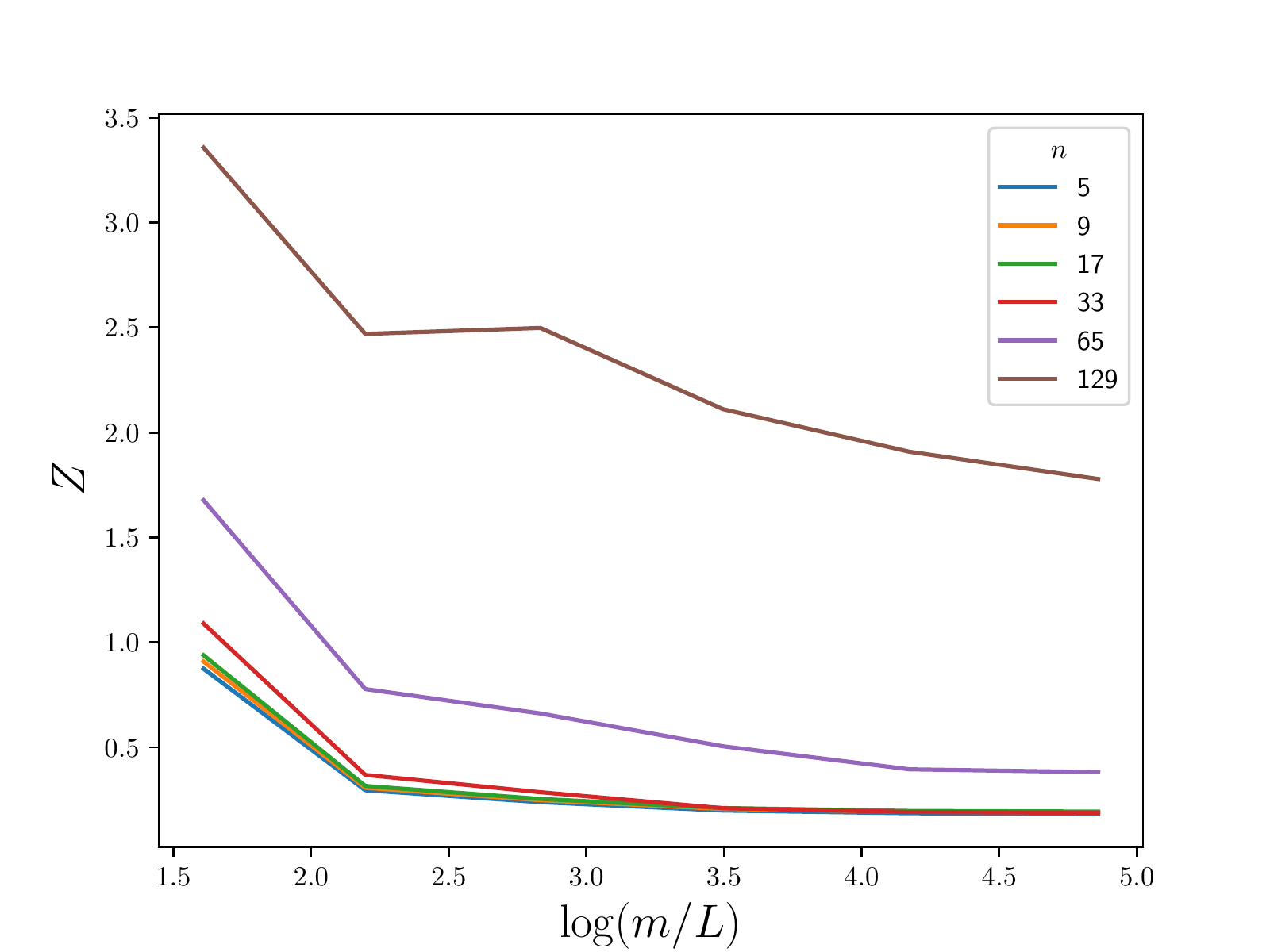}
\caption{
Forced Burger's equation:
For each pair $(n,m)$ of temporal ($n$) and spatial ($m$) grid sizes considered, we plot:
(top left) the error $E_\infty$ for fixed $m$ and varying $n$ for our PNM;
(top right) the error $E_\infty$ for fixed $n$ and varying $m$ for our PNM;
(middle left) the error $E_\infty$ for fixed $m$ and varying $n$, Crank--Nicolson method;
(middle right) the error $E_\infty$ for fixed $n$ and varying $m$, Crank--Nicolson method;
(bottom left) the $Z$-score for fixed $m$ and varying $n$;
(bottom right) the $Z$-score for fixed $n$ and varying $m$.
}
\label{fig: burgerforcederrorplots}
\end{figure}


\paragraph{Results:}
The error $E_\infty$ and $Z$-scores for our PNM are displayed in \Cref{fig: burgerforcederrorplots}.
The error $E_\infty$ is seen to be gated by the size $m$ of the spatial grid and descreases as $(n,m)$ are simultaneously increased.
The $Z$-score values appear to be of order 1 as $(n,m)$ are simultaneously increased, but for increasing $n$ with $m$ fixed the PNM appears to become over-confident; a conservative choice would be to take $m \geq n$, which is also what we concluded from the porous medium equation.
These results suggest the output from our PNM is reasonably well-calibrated.
Finally, we considered the accuracy of our PNM compared to the Crank--Nicolson benchmark.
The error $E_\infty$ for Crank--Nicolson is jointly displayed in \Cref{fig: burgerforcederrorplots}, below the error plots for our PNM, and interestingly, it is generally larger than the error obtained with our PNM.
This provides reassurance that our PNM is as accurate as could reasonably be expected.

%% file: conclusion.tex
\section{Conclusion} \label{sec: conc}

This paper addressed an important and under-studied problem in numerical analysis; the numerical solution of a PDE under severe restrictions on evaluation of the initial, boundary and/or forcing terms $f$, $g$ and $h$ in \eqref{eq: somepde}.
Such restructions occur when $f$, $g$ and/or $h$ are associated with a computational cost, such as being output from a computationally intensive computer model \citep{fulton2010approaches,hurrell2013community} or arising as the solution to an auxiliary PDE \citep{macnamara2016operator,cockayne2020probabilistic}.
In many such cases it is not possible to obtain an accurate approximation of the solution of the PDE, and at best one can hope to describe trajectories that are compatible with the limited information available on the PDE.
To provide a principled resolution, in this paper we cast the numerical solution of a nonlinear PDE as an inference problem within the Bayesian framework and proposed a probabilistic numerical method (PNM) to infer the unknown solution of the PDE.
This approach enables formal quantification of \textit{numerical uncertainty}, in such settings where the solution of the PDE cannot be easily approximated.

Our contribution extends an active line of research into the development of PNM for a range of challenging numerical tasks \citep[see the survey in][]{Hennig2015}.
A common feature of these tasks is that their difficulty justifies the use of sophisticated statistical machinery, such as Gaussian processes, that themselves may be associated with a computational cost.
The PNM developed in this paper has a complexity $O(nm^3)$ to approximate the final state of the PDE, or $O(n^3m^3)$ to approximate the full solution trajectory of the PDE.
This renders our PNM computationally intensive -- potentially orders of magnitude slower than a classical numerical method -- but such increase in cost can be justified when the demands of evaluating $f$, $g$ and $h$ exceed those of running the PNM \citep[for example, when evaluation of $f$ requires simulation from a climate model][]{fulton2010approaches,hurrell2013community}.

Further work will be required to establish our approach as a general purpose numerical tool for nonlinear PDEs:
First, the non-unique partitioning of the differential operator $D$ into linear and nonlinear components, $P$ and $Q$, together with the non-unique linearisation of $Q$, necessitates some expert input.
This is analogus to the selection of a suitable numerical method in the classical setting, but the classical literature has benefitted from decades of research and extensive practical guidence is now available in that context.
Here we took a first step to automation by rigorously establishing sample path properties of a Mat\'{e}rn tensor product covariance model $\Sigma$, along with presenting a closed-form maximum likelihood estimator for the amplitude of $\Sigma$.
The user is left to provide suitable length-scale parameter(s), which is roughly analogous to requiring the user to specify a mesh density in a finite element method (an accepted reality in that context).
Second, an extensive empirical assessment will be required to systematically assess the performance of the method; our focus in the present paper was methodology and theory, providing only an experimental proof-of-concept.
In particular, it will be important to assess diagnostics for \textit{failure} of the method; it seems plausible that statistically-motivated diagnostics, such as held-out predictive likelihood, could be used to indicate the quality of the output from the PNM.
Finally, we acknowledge that the problem we considered in \eqref{eq: somepde} represents only one class of nonlinear PDEs and further work will be required to develop PNM for other classes of PDEs, such as boundary value problems and PDEs defined on more general domains.


%% file: appendix.tex
This supplement contains proofs for all novel theoretical results in the main text.

\subsection{Proof of \Cref{lem: batch}} \label{sec: derive sequential version pde}

\begin{proof}[Proof of \Cref{lem: batch}]

Our starting point is the stochastic process defined in \Cref{lem: batch}, and from this we aim to derive the iterative formulation in the main text.
The derivation requires several items of notation to be introduced.
First, let 
$$
L^{i+1}(r) :=
\begin{bmatrix} 
\Sigma(r,\bm{a}_0) \\
\Sigma(r,\bm{b}_0) \\
\Sigma_{D_0}(r,\bm{v}_0) \\
\vdots \\
\Sigma(r,\bm{b}_i) \\
\Sigma_{\bar{D}_i}(r,\bm{v}_i)
\end{bmatrix}^\top ,
\; \;
R^{i+1}(r') :=
\begin{bmatrix} 
\Sigma(\bm{a}_0,r') \\
\Sigma(\bm{b}_0,r') \\
\Sigma_{D_0}(\bm{v}_0,r') \\
\vdots \\
\Sigma(\bm{b}_i,r') \\
\Sigma_{D_i}(\bm{v}_i,r')
\end{bmatrix} ,
\; \;
F_{i+1} :=\begin{bmatrix} 
g(\bm{a}_0)-m(\bm{a}_0)\\
h(\bm{b}_0)-m(\bm{b}_0) \\
f(\bm{v}_0)-m_{D_0}(\bm{v}_0) \\
\vdots \\
h(\bm{b}_i)-m(\bm{b}_i) \\
f(\bm{v}_i)-m_{D_i}(\bm{v}_i)
\end{bmatrix}
$$
and
$$
M_{i+1} :=
\begin{bmatrix} 
\Sigma(\bm{a}_0,\bm{a}_0) & \Sigma(\bm{a}_0,\bm{b}_0) &\Sigma_{\bar{D_0}}(\bm{a}_0,\bm{v}_0) & \dots & \Sigma(\bm{a}_0,\bm{b}_i) &\Sigma_{\bar{D_i}}(\bm{a}_0,\bm{v}_i)  \\
\Sigma(\bm{b}_0,\bm{a}_0) & \Sigma(\bm{b}_0,\bm{b}_0) &\Sigma_{\bar{D_0}}(\bm{b}_0,\bm{v}_0) & \dots & \Sigma(\bm{b}_0,\bm{b}_i) &\Sigma_{\bar{D_i}}(\bm{b}_0,\bm{v}_i) \\
\Sigma_{{D_0}}(\bm{v}_0,\bm{a}_0)& \Sigma_{{D_0}}(\bm{v}_0,\bm{b}_0) &\Sigma_{D_0\bar{D_0}}(\bm{v}_0,\bm{v}_0) & \dots & \Sigma_{{D_i}}(\bm{v}_0,\bm{b}_i) &\Sigma_{D_i\bar{D_i}}(\bm{v}_0,\bm{v}_i) \\
\vdots & \vdots & \vdots & \dots & \vdots & \vdots  \\
\Sigma(\bm{b}_i,\bm{a}_0) & \Sigma(\bm{b}_i,\bm{b}_0) & \Sigma_{\bar{D_i}}(\bm{b}_i,\bm{v}_0) & \dots & \Sigma(\bm{b}_i,\bm{b}_i) & \Sigma_{\bar{D_i}}(\bm{b}_i,\bm{v}_i) \\ 
\Sigma_{{D_i}}(\bm{v}_i,\bm{a}_0) & \Sigma_{{D_i}}(\bm{v}_i,\bm{b}_0) & \Sigma_{D_i\bar{D_i}}(\bm{v}_i,\bm{v}_0) & \dots & \Sigma_{{D_i}}(\bm{v}_i,\bm{b}_i) & \Sigma_{D_i\bar{D_i}}(\bm{v}_i,\bm{v}_i) 
\end{bmatrix} .
$$
The mean and covariance of $U^{i+1}$, as defined in \Cref{lem: batch}, are equal to
\begin{align*}
\mu^{i+1}(r) & = \mu(r)+L^{i+1}(r) (M_{i+1})^{-1}F_{i+1} \\
\Sigma^{i+1}(r,r') & = \Sigma(r,r')-L^{i+1}(r)(M_{i+1})^{-1}R^{i+1}(r')
\end{align*}
Similarly, introduce the notation
$$
L^{i+1}_{D_i}(r) :=
\begin{bmatrix} 
\Sigma_{D_0}(r,\bm{a}_0) \\
\Sigma_{D_0}(r,\bm{b}_0) \\
\Sigma_{D_0\bar{D_0}}(r,\bm{v}_0) \\
\vdots \\
\Sigma_{D_i}(r,\bm{b}_i) \\
\Sigma_{D_i\bar{D_i}}(r,\bm{v}_i)
\end{bmatrix}^\top ,
\qquad
R^{i+1}_{D_i}(r') :=
\begin{bmatrix} 
\Sigma_{\bar{D_0}}(\bm{a}_0,r') \\
\Sigma_{\bar{D_0}}(\bm{b}_0,r') \\
\Sigma_{D_0\bar{D_0}}(\bm{v}_0,r') \\
\vdots \\
\Sigma_{\bar{D_i}}(\bm{b}_i,r') \\
\Sigma_{D_i\bar{D_i}}(\bm{v}_i,r')
\end{bmatrix} ,
$$
so that the application of $D_i$ to $\mu^{i+1}$ and $\Sigma^{i+1}$ may be expressed as
\begin{align*}
\mu^{i+1}_{D_i}(r) &=\mu_{D_i}(r)+L^{i+1}_{D_i}(r)(M_{i+1})^{-1}F_{i+1} \\
\Sigma^{i+1}_{D_i}(r,r') & =\Sigma_{D_i}(r,r')-L^{i+1}_{D_i}(r) (M_{i+1})^{-1} R^{i+1}(r') \\
\Sigma^{i+1}_{\bar{D_i}}(r,r') & =\Sigma_{\bar{D_i}}(r,r')-L^{i+1}(r) (M_{i+1})^{-1} R^{i+1}_{D_i}(r') \\
\Sigma^{i+1}_{D_i\bar{D_i}}(r,r') & =\Sigma_{D_i\bar{D_i}}(r,r')-L^{i+1}_{D_i}(r) (M_{i+1})^{-1} R^{i+1}_{D_i}(r')
\end{align*}
Notice that we have a recursive partitioning of $M_{i+1}$ into blocks of the form
$$M_{i+1}=
\begin{bmatrix} 
M_i & \beta_i \\
\gamma_i & \delta_i
\end{bmatrix}
$$
where 
$$
\beta_i := 
\begin{bmatrix}
R^{i}(\bm{b}_i) & R^{i}_{D_i}(\bm{v}_i)
\end{bmatrix} ,
\; \;
\gamma_i :=
\begin{bmatrix} 
L^{i}(\bm{b}_i) \\
L^{i}_{D_i}(\bm{v}_i)
\end{bmatrix}, 
\; \;
\delta_i :=\begin{bmatrix} 
\Sigma(\bm{b}_i,\bm{b}_i) & \Sigma_{\bar{D_i}}(\bm{b}_i,\bm{v}_i) \\
\Sigma_{{D_i}}(\bm{v}_i,\bm{b}_i) & \Sigma_{D_i\bar{D_i}}(\bm{v}_i,\bm{v}_i) 
\end{bmatrix} .
$$
Thus we may use the block matrix inversion formula to deduce that
$$M_{i+1}^{-1}=
\begin{bmatrix} 
M_i^{-1}(I+\beta_i(\delta_i-\gamma_iM_i^{-1}\beta_i)^{-1}\gamma_iM_i^{-1})& -M_i^{-1}\beta_i(\delta_i-\gamma_iM_i^{-1}\beta_i)^{-1} \\
-(\delta_i-\gamma_iM_i^{-1}\beta_i)^{-1}\gamma_iM_i^{-1} & (\delta_i-\gamma_iM_i^{-1}\beta_i)^{-1}
\end{bmatrix}
$$
Setting $A_i := \delta_i-\gamma_iM_i^{-1}\beta_i$, we observe that
\begin{align*}
A_i=
\begin{bmatrix} 
 \Sigma^i(\bm{b}_i,\bm{b}_i) & \Sigma^i_{\bar{D_i}}(\bm{b}_i,\bm{v}_i) \\
\Sigma^i_{{D_i}}(\bm{v}_i,\bm{b}_i) & \Sigma^i_{D_i\bar{D_i}}(\bm{v}_i,\bm{v}_i)
\end{bmatrix} ,
\end{align*}
so our definition of $A_i$ coincides with that in the main text, and enables us to simplify $M_{i+1}^{-1}$ into
$$M_{i+1}^{-1}=
\begin{bmatrix} 
M_i^{-1}(I+\beta_iA_i^{-1}\gamma_iM_i^{-1})& -M_i^{-1}\beta_iA_i^{-1} \\
-A_i^{-1}\gamma_iM_i^{-1} & A_i^{-1}
\end{bmatrix} .
$$
Therefore we have that
\begin{align*}
\mu^{i+1}(r)&=\mu(r)+L^{i+1}(r) (M_{i+1})^{-1}F_{i+1} \\
&=\mu(r)+L^{i}(r)M_i^{-1} F_{i} +L^{i}(r)M_i^{-1}\beta_iA_i^{-1}\gamma_iM_i^{-1}F_{i} - [\Sigma(r,\bm{b}_i),\Sigma_{\bar{D_i}}(r,\bm{v}_i)]A_i^{-1}\gamma_iM_i^{-1}F_{i} \\
& \quad - L^{i}(r)M_i^{-1}\beta_iA_i^{-1}
\begin{bmatrix} 
h(\bm{b}_i)-\mu(\bm{b}_i) \\
f(\bm{v}_i)-\mu_{D_i}(\bm{v}_i)
\end{bmatrix} +[\Sigma(r,\bm{b}_i),\Sigma_{\bar{D_i}}(r,\bm{v}_i)]A_i^{-1}
\begin{bmatrix} 
h(\bm{b}_i)-\mu(\bm{b}_i) \\
f(\bm{v}_i)-\mu_{D_i}(\bm{v}_i)
\end{bmatrix} \\
&= \mu^i(r)+\{L^{i}(r)M_i^{-1}\beta_i- [\Sigma(r,\bm{b}_i),\Sigma_{\bar{D_i}}(r,\bm{v}_i)]\}A_i^{-1}\gamma_iM_i^{-1}F_{i} \\
& \quad +  \{[\Sigma(r,\bm{b}_i),\Sigma_{\bar{D_i}}(r,\bm{v}_i)] -  L^{i}(r)M_i^{-1}\beta_i\}
A_i^{-1}
\begin{bmatrix} 
h(\bm{b}_i)-\mu(\bm{b}_i) \\
f(\bm{v}_i)-\mu_{D_i}(\bm{v}_i)
\end{bmatrix}
\end{align*}
which can be simplified by noting that
\begin{align*}
L^{i}(r)M_i^{-1}\beta_i- [\Sigma(r,\bm{b}_i),\Sigma_{\bar{D_i}}(r,\bm{v}_i)]&=
L^{i}(r)M_i^{-1}
\begin{bmatrix}
R^{i}(\bm{b}_i) & R^{i}_{D_i}(\bm{v}_i)
\end{bmatrix}
-[\Sigma(r,\bm{b}_i),\Sigma_{\bar{D_i}}(r,\bm{v}_i)]
\\
&=[-\Sigma^i(r,\bm{b}_i),-\Sigma_{\bar{D_i}}^i(r,\bm{v}_i)] \\
\gamma_iM_i^{-1}F_{i} &=
\begin{bmatrix} 
L^{i}(\bm{b}_i) \\
L^{i}_{D_i}(\bm{v}_i)
\end{bmatrix}M_i^{-1}F_{i} 
=\begin{bmatrix} 
\mu^{i}(\bm{b}_i)-\mu(\bm{b}_i) \\
\mu^{i}_{D_i}(\bm{v}_i)-\mu_{D_i}(\bm{v}_i)
\end{bmatrix}
\end{align*}
to produce the iterative formulation for the mean in the main text:
\begin{align*}
\mu^{i+1}(r) &= \mu^i(r)+ [\Sigma^i(r,\bm{b}_i),\Sigma_{\bar{D_i}}^i(r,\bm{v}_i)]A_i^{-1} \left( \begin{bmatrix} 
h(\bm{b}_i)-\mu(\bm{b}_i) \\
f(\bm{v}_i)-\mu_{D_i}(\bm{v}_i)
\end{bmatrix}
-\begin{bmatrix} 
\mu^{i}(\bm{b}_i)-\mu(\bm{b}_i) \\
\mu^{i}_{D_i}(\bm{v}_i)-\mu_{D_i}(\bm{v}_i)
\end{bmatrix} \right) \\
&= \mu^i(r)+ [\Sigma^i(r,\bm{b}_i),\Sigma_{\bar{D_i}}^i(r,\bm{v}_i)]A_i^{-1}\begin{bmatrix} 
h(\bm{b}_i)-\mu^{i}(\bm{b}_i) \\
f(\bm{v}_i)-\mu^{i}_{D_i}(\bm{v}_i)
\end{bmatrix}
\end{align*}

For the covariance we have
\begin{align*}
\Sigma^{i+1}(r,r')&=\Sigma(r,r')-L^{i+1}(r)(M_{i+1})^{-1}R^{i+1}(r') \\
&=\Sigma(r,r')-L^{i}(r) M_i^{-1}R^{i}(r')-L^{i}(r)M_i^{-1}\beta_iA_i^{-1}\gamma_iM_i^{-1}R^{i}(r') 
\\
& \quad + [\Sigma(r,\bm{b}_i),\Sigma_{\bar{D_i}}(r,\bm{v}_i)]A_i^{-1}\gamma_iM_i^{-1}R^{i}(r') 
\\
& \quad + L^{i}(r)M_i^{-1}\beta_iA_i^{-1}\begin{bmatrix} 
\Sigma(\bm{b}_i,r') \\
\Sigma_{D_i}(\bm{v}_i,r')
\end{bmatrix}
- [\Sigma(r,\bm{b}_i),\Sigma_{\bar{D_i}}(r,\bm{v}_i)]A_i^{-1}\begin{bmatrix} 
\Sigma(\bm{b}_i,r') \\
\Sigma_{D_i}(\bm{v}_i,r')
\end{bmatrix} \\
&= \Sigma^i(r,r')-\{L^{i}(r)M_i^{-1}\beta_i - [\Sigma(r,\bm{b}_i),\Sigma_{\bar{D_i}}(r,\bm{v}_i)]\}A_i^{-1}\gamma_iM_i^{-1}R^{i}(r') \\
& \quad -  \{[\Sigma(r,\bm{b}_i),\Sigma_{\bar{D_i}}(r,\bm{v}_i)] -  L^{i}(r)M_i^{-1}\beta_i\}A_i^{-1}\begin{bmatrix} 
\Sigma(\bm{b}_i,r') \\
\Sigma_{D_i}(\bm{v}_i,r')
\end{bmatrix} .
\end{align*}
This can be simplified by noting that
\begin{align*}
L^{i}(r)M_i^{-1}\beta_i- [\Sigma(r,\bm{b}_i),\Sigma_{\bar{D_i}}(r,\bm{v}_i)]&=
L^{i}(r)M_i^{-1}
\begin{bmatrix}
R^{i}(\bm{b}_i) & R^{i}_{D_i}(\bm{v}_i)
\end{bmatrix}
-[\Sigma(r,\bm{b}_i),\Sigma_{\bar{D_i}}(r,\bm{v}_i)]
\\
&=[-\Sigma^i(r,\bm{b}_i),-\Sigma_{\bar{D_i}}^i(r,\bm{v}_i)] \\
\gamma_iM_i^{-1}R^{i}(r') &=\begin{bmatrix} 
L^{i}(\bm{b}_i) \\
L^{i}_{D_i}(\bm{v}_i)
\end{bmatrix}M_i^{-1}R^{i}(r') 
=
\begin{bmatrix} 
-\Sigma^i(\bm{b}_i,r')+\Sigma(\bm{b}_i,r') \\
-\Sigma^i_{D_i}(\bm{v}_i,r')+\Sigma_{D_i}(\bm{v}_i,r')
\end{bmatrix}
\end{align*}
to obtain
\begin{align*}
\Sigma^{i+1}(r,r') &= \Sigma^i(r,r')- [\Sigma^i(r,\bm{b}_i),\Sigma_{\bar{D_i}}^i(r,\bm{v}_i)]A_i^{-1} \left( \begin{bmatrix} 
\Sigma^i(\bm{b}_i,r')-\Sigma(\bm{b}_i,r') \\
\Sigma^i_{D_i}(\bm{v}_i,r')-\Sigma_{D_i}(\bm{v}_i,r')
\end{bmatrix}
+
\begin{bmatrix} 
\Sigma(\bm{b}_i,r') \\
\Sigma_{D_i}(\bm{v}_i,r')
\end{bmatrix} \right) \\
&= \Sigma^i(r,r')- [\Sigma^i(r,\bm{b}_i),\Sigma_{\bar{D_i}}^i(r,\bm{v}_i)]A_i^{-1}\begin{bmatrix} 
\Sigma^i(\bm{b}_i,r') \\
\Sigma^i_{D_i}(\bm{v}_i,r')
\end{bmatrix} ,
\end{align*}
identical to the iterative formulation in the main text.
\end{proof}

\subsection{Proof of \Cref{thm: Maternmsdifferentiable}} \label{app: Maternmsdifferentiable}

\begin{proof}[Proof of \Cref{thm: Maternmsdifferentiable}]

The assumed regularity of $\mu$, being an element of $C^p(I)$, implies that $X$ and $X-\mu$ have identical differentiability properties up to order $p$, and we therefore assume $\mu = 0$ for simplicity in the remainder.

The mean square differentiability of the Mat\'{e}rn covariance function has been well-documented. In particular, because of the stationarity of the Mat\'{e}rn covariance function, $X$ is order $p$ mean-square differentiable if and only if $K_{\nu}^{(p,p)}(0,0)=(-1)^{p}K_{\nu}^{(2p)}(0)$ exists and is finite, and the Mat\'{e}rn covariance function with parameter $\nu$ is $2p$ times differentiable if and only if $\nu>p$; see Section 2 of \cite{Stein1999}.
This establishes existence of the mean-square derivative $\partial_{\textsc{ms}}^{(p)}X$.

It remains to prove that $\partial_{\textsc{ms}}^{(p)}X$ is mean-square continuous. 
From the discussion in \eqref{eq: stationarycontinuity}, $\partial_{\textsc{ms}}^{(p)}X$ is mean-square continuous if and only if its autocovariance function, $K_{\nu}^{(2p)}$, is continuous at 0.
Let $h=z-z'$ and $K_{\nu}(h)=f(h)g(h)$ where
$$
f(h)=\sigma^2 \exp\left(-\frac{|h|}{\rho} \right) \frac{p!}{(2p)!}, \qquad 
g(h)=\sum_{k=0}^{p} \frac{(2p-k)!}{(p-k)!k!} \left(\frac{2}{\rho} \right)^{k}|h|^{k}
$$
so that, by Leibniz's generalised product rule, for $m \in \mathbb{N}_0$,
$$
K_{\nu}^{(m)}(h)=\sum_{n=0}^{m} {m\choose n}f^{(m-n)}(h)g^{(n)}(h).
$$ 
One can verify that, for $n \in \{0,1,\dots,m\}$,
\begin{align*}
f^{(m-n)}(h) & = \left\{ \begin{array}{ll} 
\frac{(-1)^{m-n}}{\rho^{m-n}}\sigma^2 \exp(-\frac{h}{\rho}) \frac{p!}{(2p)!} & h > 0 \\
\frac{1}{\rho^{2p-n}}\sigma^2 \exp(-\frac{-h}{\rho}) \frac{p!}{(2p)!} & h < 0
 \end{array} \right. \\
 g^{(n)}(h) & = \left\{ \begin{array}{ll}
	\sum_{k=0}^{p-n} \frac{(2p-n-k)!}{(p-n-k)!k!}(\frac{2}{\rho})^{n+k}h^{k} & h > 0 , n\leq p \\
	(-1)^{n}\sum_{k=0}^{p-n} \frac{(2p-n-k)!}{(p-n-k)!k!}(\frac{2}{\rho})^{n+k}(-h)^{k} & h < 0, n\leq p \\
	0 & h \neq 0, n>p 
\end{array} \right.
\end{align*}
from which it follows that

%
%
%
%
%
%
%

\begin{align*}
\lim_{h \downarrow 0} K_{\nu}^{(2p)}(h) &= \frac{p!}{(2p)!} \sigma^2 \sum_{n=0}^{p} {2p\choose n} \frac{(2p-n)!}{(p-n)!}\frac{(-1)^{2p-n}}{\rho^{2p-n}} \left( \frac{2}{\rho} \right)^{n} \\
&=  \frac{\sigma^2}{\rho^{2p}} \sum_{n=0}^{p} {p\choose n} (-1)^{2p-n}2^{n} 
=(-1)^{p} \frac{\sigma^2}{\rho^{2p}} (2-1)^p
=(-1)^{p} \frac{\sigma^2}{\rho^{2p}}
\end{align*}
and an analogous calculation shows
\begin{align*}
\lim_{h \uparrow 0} K_{\nu}^{(2p)}(h) & = (-1)^{p} \frac{\sigma^2}{\rho^{2p}} .
\end{align*}
%
%
%
%
%
%
%
%
%
Finally, we must check that the value $K_{\nu}^{(2p)}(0)$ agrees with the two limits just derived. 
$K_{\nu}^{(2p-1)}(h)$ is continuously differentiable, so $K_{\nu}^{(2p-1)}(0)=0$ because it is an odd function (as it is an odd derivative of an even function $K_{\nu}$). 
Thus we have that
\begin{align*}
K_{\nu}^{(2p)}(0)&=\lim_{h \to 0} \cfrac{K_{\nu}^{(2p-1)}(h)-K_{\nu}^{(2p-1)}(0)}{h} 
= \lim_{h \to 0} \cfrac{K_{\nu}^{(2p-1)}(h)}{h} \\
&=\lim_{h \to 0}\frac{p!}{(2p)!} \sigma^2 \exp \left( -\frac{|h|}{\rho} \right) \sum_{n=0}^{p-1} {2p-1\choose n} \frac{(2p-n-1)!}{(p-n-1)!}\frac{(-1)^{2p-1-n}}{\rho^{2p-1-n}} \left(\frac{2}{\rho} \right)^{n+1} +O(h) \\
&=(-1)^p\sigma^2 \sum_{n=0}^{p-1} {p-1\choose n}\frac{(-1)^{p-1-n}2^{n}}{\rho^{2p}}=(-1)^p\frac{\sigma}{\rho^{2p}}(2-1)^{p-1}=(-1)^p\frac{\sigma}{\rho^{2p}}
\end{align*}
as required.
\end{proof}

\subsection{Proof of \Cref{thm: Maternderivativesamplecont}} \label{app: Maternderivativesamplecont}

As in the proof of \Cref{thm: Maternmsdifferentiable}, the assumed regularity of $\mu$, being an element of $C^p(I)$, implies that $X$ and $X-\mu$ have identical differentiability properties up to $p$th order, and we may therefore assume $\mu = 0$.

Our main tool is the \emph{Kolmogorov continuity theorem} (see for example \cite{Kunita1997}, Section 1.4):

\begin{theorem}[Kolmogorov's Continuity Theorem]\label{thm: Kolmogorovcontinuity}
Let $I \subseteq \mathbb{R}^d$ be an open set, and let $\mathbf{z}$ be a dense subset of $I$. 
Let $X:I \times \Omega \rightarrow \mathbb{R}$ be a random field. If there exists constants $\alpha, \beta, C > 0$ such that
\begin{equation*}
\mathbb{E} \left[ |X(z,\omega)-X(z',\omega)|^{\alpha} \right] \leq C\|z-z'\|^{d+\beta}
\label{eq: Kolmogorovinequality}
\end{equation*}
for all $z,z' \in \mathbf{z}$, then there exists a modification of $X$ that is sample continuous.
\end{theorem}

\begin{lemma}\label{lemma: 2ndordermomentsonly}
Let $I \subseteq \mathbb{R}^d$ be an open set and suppose that a positive definite function $\Sigma : I \times I \rightarrow \mathbb{R}$ satisfies
\begin{equation*}
\Sigma(z,z)+\Sigma(z',z')-2\Sigma(z,z') \leq C\|z-z'\|^\gamma
\label{eq: Kolmogorovsimpler}
\end{equation*}
for some $\gamma, C \in (0,\infty)$ and all $z,z' \in I$. 
Let $X \sim \mathcal{GP}(0,\Sigma)$.
Then there exists a modification of $X$ that is sample continuous. 
\end{lemma}

\begin{proof}[Proof of \Cref{lemma: 2ndordermomentsonly}]
Notice that 
$$
\Sigma(z,z)+\Sigma(z',z')-2\Sigma(z,z') = \mathbb{E}[(X(z,\omega)-X(z',\omega))^{2}] ,
$$
so if $\gamma>d$ then the required result follows from Kolmogorov's continuity theorem (\Cref{thm: Kolmogorovcontinuity}, with $\alpha = 2$). 
If not, then we can consider higher order moments via Isserlis's theorem: 
For $n \in \mathbb{N}$,
$$
\mathbb{E}[(X(z,\omega)-X(z',\omega))^{2n}] =
	\frac{(2n)!}{2^n n!} (\mathbb{E}[(X(z,\omega)-X(z',\omega))^{2}])^{n}
$$
and thus, with any $n > d / \gamma$, we have
$$
\mathbb{E}[(X(z,\omega)-X(z',\omega))^{2n}] \leq \frac{(2n)!}{2^n n!}C^{n} \|z-z'\|^{\gamma n} 
\leq \tilde{C} \|z-z'\|^{d+\beta}
$$
with $\tilde{C} = \frac{(2n)!}{2^n n!}C^{n}$ and $\beta = \gamma n - d$.
The result then follows from Kolmogorov's continuity theorem (\Cref{thm: Kolmogorovcontinuity}, with $\alpha = 2n$).
\end{proof}

\begin{proof}[Proof of \Cref{thm: Maternderivativesamplecont}]

Our aim is to show that $\partial_{\textsc{ms}}^{(i)}X$ satisfies the preconditions of \Cref{lemma: 2ndordermomentsonly}.
This process has covariance function $K_\nu^{(i,i)}(z,z')=(-1)^i K_\nu^{(2i)}(z-z')$.
From stationarity we have, with $h = z-z'$,
\begin{equation*}
K_\nu^{(2i)}(z,z)+K_\nu^{(2i)}(z',z')-2K_\nu^{(2i)}(z,z')=2K_\nu^{(2i)}(0)-2K_\nu^{(2i)}(h) .
\end{equation*}

From similar calculations to those performed in the proof of \Cref{thm: Maternmsdifferentiable}, we have that for all $0 \leq i \leq p$,
\begin{align*}
K_\nu^{(2i)}(h) &= \frac{2 \sigma^2 p!}{(2p)!}  \exp\left(-\frac{|h|}{\rho}\right)\sum_{n=0}^{\min(2i,p)} \sum_{k=0}^{p-n}{2i\choose n} \frac{(2p-n-k)!}{(p-n-k)!k!}\frac{(-1)^{2i-n+1}}{\rho^{2i-n}}\left(\frac{2}{\rho} \right)^{n+k}|h|^{k} \\
&= \exp\left(-\frac{|h|}{\rho}\right) \sum_{k=0}^{p}a_k |h|^{k}
\end{align*}
for some real coefficients $a_k$. Therefore:

\begin{align}
|2K_\nu^{(2i)}(0)-2K_\nu^{(2i)}(h)| &= 2\left | a_0-a_0\exp\left(-\frac{|h|}{\rho}\right) -\exp\left(-\frac{|h|}{\rho}\right) \sum_{k=1}^{p}a_k |h|^{k} \right | \nonumber \\
& \leq 2\left | a_0-a_0\exp\left(-\frac{|h|}{\rho}\right) \right | + 2\left | \exp\left(-\frac{|h|}{\rho}\right) \sum_{k=1}^{p}a_k |h|^{k} \right | \label{eq: poly bd}
\end{align}

This final term can be upper bounded by an expression of the form $C|h|^\gamma$ for sufficiently large $C > 0$ and $\gamma=1$.
Indeed, as $h \rightarrow 0$ the behaviour of \eqref{eq: poly bd} is $O(|h|)$. 
As $|h| \rightarrow \infty$ the exponential term dominates and \eqref{eq: poly bd} decays to $a_0$. 
In the region $0<|h|<\infty$, \eqref{eq: poly bd} is smooth. 
Thus we can use \Cref{lemma: 2ndordermomentsonly} to conclude that $\partial_{\textsc{ms}}^{(i)}X$  has a modification that is sample continuous.
\end{proof}

\subsection{Proof of \Cref{thm: Sampledifferentiabilitynthordercriterion}} \label{app: Sampledifferentiabilitynthordercriterion}

For $p=1$ the result is immediate from \Cref{thm: Sampledifferentiabilitycriterion}, so in what follows we concentrate on $p > 1$. 

The main technical challenge of this proof is to deal with modifications, which arise with each application of \Cref{thm: Sampledifferentiabilitycriterion}.
Recall that two stochastic processes $X$ and $\tilde{X}$ are said to be \textit{indistinguishable} if $\mathbb{P}(X(z,\omega) = \tilde{X}(z,\omega) \; \forall z \in I) = 1$.
If $X$ and $\tilde{X}$ are modifications of each other and each is sample continuous, then $X$ and $\tilde{X}$ are indistinguishable \citep[][Section 1.1]{Jeanblanc2009}.

\begin{proof}[Proof of \Cref{thm: Sampledifferentiabilitynthordercriterion}]



We first present a proof for $d=1$, to improve transparency of the argument, then we present the argument for the general case $d \geq 1$.
Note that, since we are considering Gaussian processes, the requirement for a second moment in \Cref{thm: Sampledifferentiabilitycriterion} is automatically satisfied.

For each $0 \leq i < p$, it is assumed that $\partial_{\textsc{ms}}^{i}X$ has mean square derivative $\partial_{\textsc{ms}}^{i+1}X$ that is mean square continuous and sample continuous. 
\Cref{thm: Sampledifferentiabilitycriterion} therefore implies that each $\partial_{\textsc{ms}}^{i}X$ has a modification, denoted $\psi_i$, that is sample continuously differentiable, and satisfying $\partial\psi_i=\partial_{\textsc{ms}}^{i+1}X$ almost surely.
Since $\psi_i$ and $\partial_{\textsc{ms}}^{i}X$ are sample continuous they are indistinguishable; i.e. almost surely $\psi_i = \partial_{\textsc{ms}}^{i}X$.
It follows that, for each $0 \leq i < p$, we have almost surely that $\partial^i \psi_0 = \partial^{i-1} (\partial \psi_0) = \partial^{i-1} \psi_1 = \dots = \psi_i$, while for $i=p$ we have that $\partial^p \psi_0 = \partial \psi_{p-1} = \partial_{\textsc{ms}}^{p}X$.



The case $d \geq 1$ is analogous with more notation is involved; though,
since we assumed $\Sigma \in C^{(p,p)}(I \times I)$ we may employ the shorthand notation $\partial_{\textsc{ms}}^\beta X$ for all $|\beta| \leq p$ (since the order of derivatives can be freely interchanged).
For each $0 \leq i < p$ and $|\beta| = i$, it is assumed that $\partial_{\textsc{ms}}^{\beta}X$ has mean square partial derivatives $\partial_{\textsc{ms}}^{\beta + \gamma}X$, $|\gamma| = 1$, that are mean square continuous and sample continuous. 
\Cref{thm: Sampledifferentiabilitycriterion} therefore implies that each $\partial_{\textsc{ms}}^{\beta}X$ has a modification, denoted $\psi_\beta$, that is sample continuously differentiable, and satisfying $\partial^\gamma \psi_\beta=\partial_{\textsc{ms}}^{\beta + \gamma}X$ almost surely for all $|\gamma| = 1$.
Since $\psi_\beta$ and $\partial_{\textsc{ms}}^\beta X$ are sample continuous they are indistinguishable; i.e. almost surely $\psi_\beta = \partial_{\textsc{ms}}^\beta X$.
It follows that, for each $0 \leq i < p$ and $|\beta| = i$, we have almost surely that $\partial^\beta \psi_0 = \psi_\beta$, while for and $|\beta| = p$ we have that $\partial^\beta \psi_0 =  \partial_{\textsc{ms}}^{\beta}X$.
\end{proof}

\subsection{Proof of \Cref{thm: product of Matern}} \label{app: product of Matern}

The following technical result will be exploited:

\begin{theorem}[Integrated Gaussian Process; Corollary 1 and 2 of \cite{LASINGER1993}]\label{thm: gaussianprocessintegrated}

Let $I=(0,T)$ and let $X : I \times \Omega \rightarrow \mathbb{R}$ be a stationary Gaussian process with mean function $\mu = 0$ and autocovariance function $\Sigma$.
Let $q \in \mathbb{N}$ and recursively define $X^{(-i)}$ for $i = 1,\dots,q$ as follows:
\begin{equation*}
X^{(-i)}(t) = X^{(-i)}(0)+\int_{0}^{t} X^{(-i+1)}(t') \, \mathrm{d}t'
\end{equation*}
Let $A_i$, $i = 1,\dots,q$ be constants and recursively define $\Sigma^{(-i,-i)}$ for $i = 1,\dots,q$ as follows:
\begin{equation*}
\Sigma^{(-i,-i)}(t) = A_i^2 - \int_{0}^{t} \int_{0}^{s} \Sigma^{(- i+1,-i+1)}(s') \, \mathrm{d}s' \mathrm{d}s .
\end{equation*}
Let $\mathcal{H}(\Sigma)$ denote the reproducing kernel Hilbert space associated with the covariance function $\Sigma$, and $\mathcal{H}(\Sigma^{(-i,-i)})$ the reproducing kernel Hilbert space associated with the covariance function $\Sigma^{(-i,-i)}$.
If
\begin{align}
\partial^{q+i-1} \Sigma^{(-q,-q)}(\cdot) & \in \mathcal{H}(\Sigma) \label{eq: cond1}  
\end{align}
and
\begin{align}
A_i^2 >\| \partial \Sigma^{(-i,-i)}(\cdot) \|^2_{\mathcal{H}(\Sigma^{(-(i-1),-(i-1))})} \label{eq: cond2}  
\end{align}
for each $1 \leq i \leq q$, then there exists a stationary Gaussian process $X^{(-q)} : I \times \Omega \rightarrow \mathbb{R}$ with covariance function $\Sigma^{(-q,-q)}(t,s) := \Sigma^{(-q,-q)}(t-s)$ such that $\partial^q \bar{\partial}^q \Sigma^{(-q,-q)}(t,s) = \Sigma(t-s)$.
Moreover, \eqref{eq: cond1} implies $\partial \Sigma^{(-i,-i)}(\cdot) \in \mathcal{H}(\Sigma^{(-(i-1),-(i-1))})$, so that $A_i$ satisfying  \eqref{eq: cond2} exist. 
\end{theorem}

From \Cref{thm: gaussianprocessintegrated} we obtain the following intermediate result:

\begin{proposition}\label{thm: Maternintegralsamplecont}
Let $I=(0,T)$ and let $\mu \in C^p(I)$ be bounded in $I$. 
Let $X \sim \mathcal{GP}(\mu,K_{\nu})$, with $K_\nu$ as in \eqref{eq: Maternkernel}.
Let $q \in \{0,\dots,p\}$.
Then there exists a stationary Gaussian process $X^{(-q)} : I \times \Omega \rightarrow \mathbb{R}$ whose mean function $\mu^{(-q)}$ satisfies $\partial^q \mu^{(-q)}(t)=\mu(t)$ and whose covariance function $K_{\nu}^{(-q,-q)}(\cdot)$ satisfies $\partial^q \bar{\partial}^q K_\nu^{(-q,-q)}(t,s) = K_{\nu}(t,s)$. 
Furthermore, $X^{(-q)}$ is sample continuous.
\end{proposition}

\begin{proof}
Since $\mu \in C^p(I)$ we may define $\mu^{(-q)}$ to be the $q$ times integrated $\mu$ for each $q \in \{0,\dots,p\}$.
Note that $\mu^{(-q)}$ is well-defined since $\mu$ is assumed to be bounded.


Our aim is to use \Cref{thm: gaussianprocessintegrated}, so we must check that the preconditions are satisfied.
\eqref{eq: cond1} is satisfied because $H_K$ for the Mat\'{e}rn $\nu$ kernel is norm equivalent to the Sobolev Space $W^\alpha_2$ where $\alpha=\nu+1/2$ \citep{KanHenSejSri18} and integrals of $W^\alpha_2$ functions are $W^\alpha_2$ functions. 
\eqref{eq: cond2} is automatically satisfied since we may select suitably large $A_i^2$. 
Thus \Cref{thm: gaussianprocessintegrated} is satisfied and there exists a stationary Gaussian process $X^{(-q)} : I \times \Omega \rightarrow \mathbb{R}$ whose mean function $\mu^{(-q)}$ satisfies $\partial^q \mu^{(-q)}(t)=\mu(t)$ and whose covariance function $K_{\nu}^{(-q,-q)}(\cdot)$ satisfies $\partial^q \bar{\partial}^q K_\nu^{(-q,-q)}(t,s) = K_{\nu}(t,s)$. 

It remains only to show that there is a modification of $X^{(-q)}$ is sample continuous.
Since $K_{\nu}(t)$ is of the form $K_{\nu}(t)=\exp(-t/\rho)\sum_{k=0}^{p} a_k t^k$ for some real constants $a_k$, we have
$$K_{\nu}^{(-q,-q)}(t)=\exp \left( -\frac{t}{\rho} \right) \sum_{k=0}^{p} b_k t^k+\sum_{k=0}^{2q-1} c_k t^k$$
for some real constants $b_k$ and $c_k$.
Now we check the preconditions of \Cref{lemma: 2ndordermomentsonly}:
\begin{align*}
& \hspace{-15pt} K_{\nu}^{(-q,-q)}(t,t) + K_{\nu}^{(-q,-q)}(s,s) - 2 K_{\nu}^{(-q,-q)}(t,s)  \\
& = 2K_{\nu}^{(-q,-q)}(0)-2K_{\nu}^{(-q,-q)}(|t-s|) \\
& = 2b_0 - 2b_0\exp \left( -\frac{|t-s|}{\rho} \right) - 2\exp \left( -\frac{|t-s|}{\rho} \right)\sum_{k=1}^{p} b_k |t-s|^k - 2 \sum_{k=1}^{2q-1} c_k |t-s|^k \\
& \leq C|t-s| 
\end{align*}
on a bounded domain, for some constant $C$. 
Therefore from \Cref{lemma: 2ndordermomentsonly} there exists a modification of $X^{(-q)}$  that is sample continuous, as required.
\end{proof}

Finally we present the proof of \Cref{thm: product of Matern}.
The main technical challenge is to smooth the process in such a way that it has an equal number of derivatives in each argument, to allow for a more straightforward application of earlier results.

\begin{proof}[Proof of \Cref{thm: product of Matern}]

If the $a_i$, $b_i$ are not equal to 0 and 1 respectively, we perform a linear re-scaling so that without loss of generality we can consider the hyper-rectangle $I = (0,1)^d$ in the remainder.






Let $p = \max_i \beta_i$ and let $P \in \mathbb{N}_0^d$ denote the vector whose entries are all equal to $p$. 
Consider a stationary Gaussian process $Z : I \times \Omega \rightarrow \mathbb{R}$ with mean function $\mu^{(\beta - P )}(\cdot)$ and autocovariance function 
\begin{align}
\Sigma^{(\beta - P , \beta - P )}(t-s) &=\prod_{i=1}^{d} K_{\nu_i}^{(\beta_i-p ,\beta_i-p )}(t_i-s_i)  ,
\label{eq: Maternproductintegralcovariance}
\end{align}
expressed using the notation of \Cref{thm: gaussianprocessintegrated}.
This construction ensures that $Z$ is order $p$ mean-square differentiable, since $\mu^{(\beta - P )} \in C^{p}(I)$ and since \eqref{eq: Maternproductintegralcovariance}  is a product of functions $K_{\nu_i}^{(\beta_i-p ,\beta_i-p )}(\cdot) \in C^{(p,p)}(I \times I)$, as a consequence of \Cref{thm: Maternmsdifferentiable,thm: Maternintegralsamplecont}.

Next we show that $\partial_{\textsc{ms}}^{\alpha} Z$ has a sample continuous modification, for any $\alpha \in \mathbb{N}_0^d, |\alpha| \leq p$, using \Cref{lemma: 2ndordermomentsonly}.
To apply \Cref{lemma: 2ndordermomentsonly} it suffices to show that 
\begin{align}
\left| \partial^{\alpha} \bar{\partial}^{\alpha} \Sigma^{(\beta - P  , \beta - P )}(0) -  \partial^{\alpha} \bar{\partial}^{\alpha} \Sigma^{(\beta - P  , \beta - P )}(t) \right| \leq C \|t\|^\gamma \label{eq: Kolmog cdn 2}
\end{align}
for some $\gamma, C \in (0,\infty)$ and all $t \in I$.
From \Cref{thm: Maternintegralsamplecont}, the first term on the left hand side satisfies $\partial^{\alpha} \bar{\partial}^{\alpha} \Sigma^{(\beta - P  , \beta - P )}(t,s) = \Sigma^{(\beta+\alpha-P,\beta+\alpha-P)}(t,s)$ and thus \eqref{eq: Kolmog cdn 2} becomes
\begin{align}
\left| \Sigma^{(\beta+\alpha-P,\beta+\alpha-P)}(0) -   \Sigma^{(\beta+\alpha-P,\beta+\alpha-P)}(t) \right| \leq C \|t\|^\gamma . \label{eq: Kolmog cdn 2 again}
\end{align}
Since $\beta_i+\alpha_i-P_i \leq \beta_i$, and from the product form of $\Sigma^{(\beta+\alpha-P,\beta+\alpha-P)}$ and the univariate calculations performed in the proofs of \Cref{thm: Maternderivativesamplecont} and \Cref{thm: Maternintegralsamplecont}, \eqref{eq: Kolmog cdn 2 again} is verified.
Thus, by \Cref{lemma: 2ndordermomentsonly}, $\partial_{\textsc{ms}}^{\alpha} Z$ has a sample continuous modification, for any $\alpha \in \mathbb{N}_0^d, |\alpha| \leq p$.

Finally, we apply \Cref{thm: Sampledifferentiabilitynthordercriterion} to deduce that there exists a modification $\tilde{Z}$ of $Z$ with $\mathbb{P}(\tilde{Z} \in C^{p}(I))$.
It follows that $\tilde{X} = \partial^{P-\beta}\tilde{Z}$ is a modification of $X$ that satisfies $\mathbb{P}(\tilde{X} \in C^\beta(I))$. 
\end{proof}

%% file: acknowledgements.tex
JW was supported by the EPSRC Centre for Doctoral Training in Cloud Computing for Big Data at Newcastle University, UK. 
JC was supported by Wave 1 of the UKRI Strategic Priorities Fund under the EPSRC Grant EP/T001569/1, particularly the ``Digital Twins for Complex Engineering Systems'' theme within that grant, and the Alan Turing Institute, UK.
TJS is supported in part by the Deutsche Forschungsgemeinschaft through Project 415980428.
CJO was supported by the Lloyd’s Register Foundation programme on data-centric engineering at the Alan Turing Institute, UK.
The authors are grateful to Andrew Duncan for feedback on the manuscript.